\documentclass[a4paper,reqno]{amsart}
\usepackage{amssymb,latexsym,amsmath,amscd,amsthm,amsfonts, enumerate}
\usepackage{multirow}
\usepackage{color}
\usepackage[all]{xy}
\usepackage{etex}
\usepackage{caption}
\usepackage{soul}
\usepackage{float}
\usepackage{titletoc}
\usepackage[normalem]{ulem}
\usepackage{graphicx}

\usepackage{tikz,tikz-cd}
\usepackage{mathrsfs}
\usepackage[all]{xy}
\usepackage{tikz}
\usepackage{extarrows}
\usepackage{tikz-cd}
\usetikzlibrary{calc}
\usetikzlibrary{matrix,arrows,decorations.pathmorphing}
\usetikzlibrary{snakes}
\usetikzlibrary{shapes.geometric,positioning}
\usetikzlibrary{arrows,decorations.pathmorphing,decorations.pathreplacing}
\usetikzlibrary{positioning,shapes,shadows,arrows,snakes}

\usepackage[colorlinks=true,hyperindex]{hyperref}
\newcommand\myshade{85}
\colorlet{mylinkcolor}{violet}
\colorlet{mycitecolor}{red}
\colorlet{myurlcolor}{cyan}

\hypersetup{
  linkcolor  = mylinkcolor!\myshade!black,
  citecolor  = mycitecolor!\myshade!black,
  urlcolor   = myurlcolor!\myshade!black,
  colorlinks = true,
}

\numberwithin{equation}{section}

\newtheorem{theorem}{Theorem}[section]
\newtheorem{proposition}[theorem]{Proposition}
\newtheorem{proposition-definition}[theorem]{Proposition-Definition}

\newtheorem{corollary}[theorem]{Corollary}
\newtheorem{lemma}[theorem]{Lemma}
\newtheorem{theoremA}{Theorem}

\theoremstyle{definition}
\newtheorem{remark}[theorem]{Remark}
\newtheorem{example}[theorem]{Example}
\newtheorem{definition}[theorem]{Definition}

\newcommand{\exc}{\operatorname{exp}\nolimits}
\newcommand{\End}{\operatorname{End}\nolimits}
\newcommand{\thick}{\mathsf{thick}}
\newcommand{\Hom}{\mathrm{Hom}}

\newcommand{\add}{\mathsf{add}}

\newcommand{\za}{\alpha}
\newcommand{\zb}{\beta}

\newcommand{\zD}{\Delta}

\newcommand{\zg}{\gamma}
\newcommand{\zG}{\Gamma}

\newcommand{\Z}{\mathbb{Z}}

\newcommand{\calt}{\mathcal{T}}

\newcommand{\calp}{\mathcal{P}}

\DeclareMathOperator{\A}{\mathsf{A}}
\DeclareMathOperator{\ka}{\mathsf{K^b(A)}}
\renewcommand{\P}{P^\bullet}

\newcommand{\xra}{\xrightarrow}

\newcommand{\ra}{\rightarrow}

\def\s{\stackrel}

\definecolor{dark-green}{RGB}{14,150,2}
\definecolor{red}{RGB}{250,0,0}
\newcommand{\gpoint}{\color{dark-green}{\circ}}
\newcommand{\rpoint}{\color{red}{\bullet}}

\UseRawInputEncoding
\begin{document}

\title[Exceptional sequences in the derived category of a gentle algebra]{Exceptional sequences in the derived category of a gentle algebra}

\author{Wen Chang}
\address{School of Mathematics and Statistics, Shaanxi Normal University, Xi'an 710062, China}
\email{changwen161@163.com}

\author{Sibylle Schroll}
\address{Insitut f\"ur Mathematik, Universit\"at zu K\"oln, Weyertal 86-90, K\"oln, Germany and
Institutt for matematiske fag, NTNU, N-7491 Trondheim, Norway}
\email{schroll@math.uni-koeln.de}

\thanks{The first author is supported by Shaanxi Province and Shaanxi Normal University.
During most of the work on this paper, the second author was supported by the EPSRC Early Career Fellowship EP/P016294/1. }

\date{\today}


\begin{abstract}
In this paper, using the correspondence of gentle algebras and dissections of marked surfaces,
we study  full exceptional sequences in the perfect derived category $\ka$ of a gentle algebra $\A$.
We show that full exceptional sequences in  $\ka$ exist if and only if the associated marked surface has no punctures and has at least two marked points on the boundary. Furthermore, by using induction on cuts of surfaces, we characterise when an exceptional sequence can be completed to a full exceptional sequence.
If the genus of the associated surface is zero then we show that the action of the braid group  together with the grading shift on full exceptional sequences in $\ka$ is transitive. For the case of surfaces of higher genus, we reduce the  problem of transitivity to the problem of  the existence of certain sequences of pairs of exceptional objects.
Finally, we interpret the duality of a full exceptional sequence induced by the longest element in the associated symmetric group  using Koszul duality.
\end{abstract}

\maketitle
\setcounter{tocdepth}{1}

\tableofcontents

\section*{Introduction}\label{Introductions}

The notion of an exceptional object was introduced by Drezet and Le Potier in \cite{DL85} to classify vector bundles on the projective plane $\mathbb{P}^2$. The theory fast developed when the group of algebraic geometers around Rudakov in Moscow systematically studied exceptional objects \cite{R90}. They introduced  exceptional  sequences in derived categories and connected these sequences through operations called mutations. An exceptional  sequence is called a full exceptional sequence if its length equals $n$, the rank of the Grothendieck group of the category. General techniques of mutations for  exceptional sequences in triangulated categories were developed by Bondal \cite{B90}. In particular, he shows that mutations give rise to an action of the braid group $B_n$ on the set of full exceptional sequences.

While for some special derived categories the theory of exceptional sequences is well-studied,
for example,  the derived categories of hereditary algebras \cite{C93} and the derived categories of coherent sheaves over some special varieties such as weighted projective lines \cite{M04}, for a general derived category, we know little. Even essential questions such as the  the existence of full exceptional sequences is not known.

In this paper we study  full exceptional sequences in the perfect derived category $\ka$ of a gentle algebra $\A$. In particular, we consider the following questions in this context: When do full exceptional sequences exist? When  can an exceptional sequence  be completed to a full exceptional sequence? Is the braid group action together with the grading shift transitive on the set of full exceptional sequences? And finally, what are  the relations between  Koszul duality, Serre duality and  dualities of full exceptional sequences induced by the longest element of a symmetric group associated to the sequence?

Gentle algebras are classical objects in the representation theory of associative algebras. They were introduced  in the 1980s  as a  generalization of iterated tilted algebras of type $A_n$ \cite{AH81}, and affine type $\widetilde A_n$ \cite{AS87}.
Remarkably, they connect to many other areas of mathematics. For example, they play an important role in the homological mirror symmetry of surfaces, see e.g. \cite{B16, HKK17, BD18}. In particular, it has been shown  that the perfect derived category of a homologically smooth gentle algebra is triangle equivalent to the partially wrapped Fukaya category of a graded oriented smooth surface with marked points on the boundary \cite{HKK17, LP20}.

From the perspective of representation theory, it is shown in \cite{OPS18} that any gentle algebra $\A$ gives rise to a dissection $\zD_{\A}$ of a marked surface $(S,M)$, where $S$ is an oriented surface with boundary, $M$ is a union of $\gpoint$-marked points and $\rpoint$-marked points which alternatively appear on the boundary $\partial S$ as well as a finite set of marked points (called punctures) in the interior of $S$. Conversely, any gentle algebra arises in this way. The objects and morphisms in $\ka$ can also been interpreted on the surface model $(S,M,\zD_{\A})$. In particular, any indecomposable object $\P_{(\zg,f)}$ corresponds to a graded curve $(\zg,f)$ on the surface, where $\zg$ is a closed curve or an arc with $\gpoint$-points as endpoints and $f$ is an integer valued function.

Although gentle algebras have been a constant object of study and much about their representation theory is known, except for some particular cases, such as, for example, those considered in \cite{BPP17} and \cite{LP18}, we know little  in general about  exceptional sequences in their derived categories.
The aim of this paper is to use the surface model established in \cite{OPS18} to determine the existence of full exceptional sequences in $\ka$ and to then study their properties.
In \cite{HKK17, LP20} it is shown that  $\ka$ is a particular realisation of a partially wrapped Fukaya category of a surface with stops and our results imply that full exceptional sequences exist in these partially wrapped Fukaya categories corresponding to $\ka$ except for certain small exceptions.

More precisely, our first result is a characterisation  of  gentle algebras that admit full exceptional sequences in their perfect derived categories in terms of their surface models.
For convenience, we assume the algebras to be connected, so that the associated surface is also connected. We denote by $\mathbb{T}_{(g,1,1)}$ the marked surface of genus $g \geq 1$ with exactly one boundary component and exactly one $\gpoint$-marked point on it.

\begin{theoremA}[Theorem \ref{thm:existence}]
Let $\A$ be a gentle algebra with surface model $(S,M,\zD_{\A})$. The following are equivalent:
 \begin{enumerate}[\rm(1)]
 \item There exists a full exceptional sequence in $\ka$; \item There are at least two $\gpoint$-points in $M$ and there are no $\rpoint$-marked points in $S \setminus \partial S$;
 \item $(S,M)$ is not homeomorphic to $\mathbb{T}_{(g,1,1)}$ and there are no $\rpoint$-marked points in $S \setminus \partial S$.
 \end{enumerate}
\end{theoremA}

Theorem A translates into  the following criterion to determine the existence of full exceptional sequences in $\ka$ by directly looking at the quiver and the relations of the algebra $\A$. For this we denote by $Q_0$ the set of vertices and by $Q_1$ the set of arrows of $Q$.

\begin{theoremA}[Corollary \ref{thm:quiver}]
Let $\A=kQ/I$ be a gentle algebra. Then there exists a full exceptional sequence in $\ka$ if and only if the global dimension of $A$ is finite and the pair $(|Q_0|,|Q_1|)$  is not equal to $(2g,4g-1)$, for any integer $g \geq 1$.
\end{theoremA}

The above theorem also confirms the existence of full exceptional sequences in the derived category of a derived-discrete algebra, see Corollary \ref{cor:existence2}, which has been shown  in \cite{BPP17} by considering a semi-orthogonal decomposition induced by an exceptional object.

By considering cuts of  marked surfaces, we can use  induction techniques to determine when an exceptional sequence can be completed to a full exceptional sequence. More precisely, we show the following.

\begin{theoremA}[Theorem \ref{theorem:completion1}]
Assume that there  exists  a full exceptional sequence in $\ka$ for a gentle algebra $\A$ arising from a marked surface $(S,M)$ and let $\mathbf{\P}=(\P_{(\zg_1,f_{1})},\cdots,\P_{(\zg_m,f_{m})})$
be an exceptional sequence arising from a collection $\zD=\{\zg_1,\cdots,\zg_m\}$ of arcs on $(S,M)$. Then $\mathbf{\P}$ can be completed to a full exceptional sequence if and only if the cut surface of $(S,M)$ along the arcs in $\zD$ has no connected component of the form $\mathbb{T}_{(g,1,1)}$.
\end{theoremA}

For any triangulated category with a full exceptional sequence, it is conjectured by Bondal and Polishchuk in \cite{BP94} that the action of $\mathbb{Z}^n\ltimes B_n$ on the set of full exceptional sequences is transitive, where the group $\mathbb{Z}^n$ acts on full exceptional sequences by shifting the grading of the objects. This has been shown to hold for derived categories of hereditary algebras \cite{C93}, and also for derived categories of coherent sheaves over some special varieties such as del Pezzo surfaces, projective planes, and weighted projective lines, see more details in \cite{M04}. For  the case of the derived category of gentle algebras, we can show the following transitivity result.

\begin{theoremA}[Theorem \ref{theorem:transitivity2}, Corollary \ref{corollary:transitivity}]\label{main theorem4}
Let $\A$ be a gentle algebra arising from a marked surface with genus zero. The action of $\mathbb{Z}^n\ltimes B_n$ on the set of full exceptional sequences in $\ka$ is transitive. In particular, if $\A$ is a derived-discrete algebra, then the action of $\mathbb{Z}^n\ltimes B_n$ on the set of full exceptional sequences in $\ka$ is transitive.
\end{theoremA}

In fact, we show that the braid group action is transitive in general for any gentle algebra with full exceptional sequences if a certain reachability condition, {Condition  \textbf{RCEA} (see Definition \ref{condition:rcea}), holds  for any  two arcs corresponding to indecomposable exceptional objects (see Proposition \ref{proposition:transitivity}).  However,  it seems difficult to establish in general for which surfaces (other than genus zero) the \textbf{RCEA} condition holds.

Finally, we determine  the connection between  Koszul duality and  the right duality ${\bf R}$ and  the left duality ${\bf L}$ induced by the longest element in the associated symmetric group (see Definition \ref{definition:exceptional duality}).
For this,
let $\mathbf X=(X_1,\cdots,X_n)$ be a full exceptional sequence in $\ka$ associated with an ordered exceptional dissection $\zG=(\zg_1,\cdots,\zg_n)$.  Let $S_n$ be the symmetric group on $n$ elements and let $w_0$ be the longest element in $S_n$. We define the right dual of $\mathbf{X}$ by setting  $\mathbf{RX}=\omega_0\mathbf X$ (see Definition \ref{definition:exceptional duality} and Lemma \ref{lemma:exceptional duality} ).  Let $\zG^*=(\zg^*_n,\cdots,\zg^*_1)$ be the dual of $\zG$, where $\zg^*_i, 1 \leq i \leq n,$ is the unique $\rpoint$-arc which intersects $\zg_i$ exactly once and intersects no other arcs in $\zG$.
Denote by $\A(\zG)$ and $(\A(\zG)^!)^{\rm op}=\A(\zG^*)$ the gentle algebra and its Koszul dual  arising from $\zG$ and $\zG^*$ respectively  (we recall the constructions in subsection \ref{subsection:derived categories of gentle algebras}). The next result shows that the right dual $\mathbf{RX}$ of $\mathbf X$ is induced by a twist to the next $\gpoint$-point of every arc $\zg^*$ in the dissection $\zG^*$ corresponding to the Koszul dual $(\A(\zG)^!)^{\rm op}$ of $\A(\zG)$.

\begin{theoremA}[Theorem \ref{theorem:two dualities1}]
Let ${\bf R}\zG$ be the ordered exceptional dissection associated to the right dual ${\bf R}\mathbf X$ of $\mathbf X$. Then  ${\bf R}\zG=D(\zG^*)$, where  $D(\zG^*)$ is the twist of $\zG^*$, which is obtained by rotating both endpoints of each arc in $\zG^*$ anticlockwise to the next respective $\gpoint$-point.
\end{theoremA}

Using our geometric interpretation of the left (and right) duality of a full exceptional sequence, we  recover  a result which Bondal in  \cite{Bon90} proves for a general triangulated category with a full exceptional sequence.

\begin{theoremA}[Theorem \ref{theorem:two dualities}]
Let $\mathbb{S}$ be the Serre functor in $\ka$.
Then $\mathbb{S} (\mathbf X) = {\mathbf L}^2(\mathbf X)$  for any full exceptional sequence  $\mathbf X$ in $\ka$.
\end{theoremA}

The paper is organized as follows. We recall some background material on marked surfaces, gentle algebras and their derived categories, as well as the theory of exceptional sequences in  Section \ref{Preliminaries}.
In Section \ref{subsection:exceptional sequences in terms of surfaces dissections}, we interpret  exceptional sequences in the derived category of a gentle algebra in terms of ordered exceptional collections on the associated marked surface. We prove the existence of full exceptional sequences in Section \ref{Section:The existence}.
Section \ref{subsection:cut surface} is devoted to the study of cutting  surfaces and the characterisation of when an exceptional sequence can be completed to a full exceptional sequence. In Section \ref{subsection:Braid group action on exceptional dissections} we give a geometric realization of the braid group action on the set of full exceptional sequences, and  study the transitivity of this action. Finally, in Section \ref{subsection:three dualities}, we  describe the relation between the left and right duality of exceptional sequences, Koszul duality and Serre duality.

\section*{Acknowledgments}
The first author would like to express his thanks to Feng Luo for his patient explanations of properties of curves on surfaces. The second author would like to thank Lutz Hille for an insightful discussion on the existence of full exceptional sequences.

\section{Preliminaries}\label{Preliminaries}

In this paper, an algebra will be assumed to be basic of finite dimension over a base field $k$. A quiver will be denoted by $Q=(Q_0,Q_1)$, where $Q_0$ is the set of vertices and $Q_1$ is the set of arrows.
Arrows in a quiver are composed from left to right as follows: for arrows $a$ and $b$ we write $ab$ for the path from the source $s(a)$  of $a$ to the target $t(b)$ of $b$.
We adopt the convention that maps are also composed  from left to right, that is if $f: X \to Y$ and $g: Y \to Z$ then $fg : X \to Z$. In general, we consider left modules.
We denote by $\mathbb{Z}$ the set of integer numbers, and by $\mathbb{Z}^*$ the set of non-zero integer numbers.
For a finite set $M$, we denote by $|M|$ its cardinality.

\subsection{Marked surfaces}\label{subsection: marked surfaces}
We recall some concepts about marked surfaces associated to gentle algebras, for which there are many references such as \cite{HKK17,PPP18,LP20}, in this paper we closely follow \cite{OPS18} and \cite{APS19}.

\begin{definition}
\label{definition:marked surface}
A triple~$(S,M,P)$ is called a \emph{marked surface} if
  \begin{enumerate}[\rm(1)]
 \item $S$ is an oriented surface with non-empty boundary with
  connected components $\partial S=\sqcup_{i=1}^{b}\partial_i S$;
 \item $M = M^{\gpoint} \cup M^{\rpoint}$ is a finite set of \emph{marked points} on $\partial S$.
       The elements of~$M^{\gpoint}$ and~$M^{\rpoint}$ will be  respectively represented by symbols~$\gpoint$ and~$\rpoint$.
       Each connected component $\partial_i S$ is required to contain at least one marked point of each colour, where in general the points~$\gpoint$ and~$\rpoint$ are alternating on $\partial_i S$;
 \item $P=P^{\rpoint}$ is a finite set of marked points in the interior of $S$. We refer to these points as \emph{punctures}, and we will  also represent them by the symbol $\rpoint$.
  \end{enumerate}
\end{definition}

Unless otherwise stated, we  always assume that a marked surface is connected
and that in  the case of a disk with no punctures there are at least two $\gpoint$-marked points and two $\rpoint$-marked points on the boundary.

\begin{definition}
\label{definition:arcs}
Let $(S,M,P)$ be a marked surface.
\begin{enumerate}[\rm(1)]
 \item[$\bullet$] An \emph{$\gpoint$-arc} is a non-contractible curve, with endpoints in~$M^{\gpoint}$.
 \item[$\bullet$]  A \emph{loop} is an $\gpoint$-arc whose endpoints coincide.
 \item[$\bullet$] An \emph{$\rpoint$-arc} is a non-contractible curve, with endpoints in $M^{\rpoint}\cup P^{\rpoint}$.
 \item[$\bullet$]  An \emph{infinite-arc} is a non-contractible curve with endpoints in $M^{\gpoint}\cup P^{\rpoint}$, which has at least one endpoint in  $P^{\rpoint}$.
 \item[$\bullet$] An \emph{admissible arc} is an $\gpoint$-arc or an infinite-arc.
 \item[$\bullet$] A \emph{closed curve} is a non-contractible curve in the interior of $S$ whose endpoints coincide. We always assume a closed curve to be primitive, that is, it is not a non-trivial power of a closed curve in the fundamental group of $S$.
 \item[$\bullet$] A \emph{simple closed curve} is a closed curve without self-intersections.
 \end{enumerate}
\end{definition}

We have the following lemma about simple closed curves, which can be found, for example in \cite[Section 1.3.1]{FM12}.

\begin{lemma}\label{lemma:simple closed curves}
If $\za$ and $\zb$ are any two non-separating simple closed curves in a surface $S$, then there is a homeomorphism $\phi: S\mapsto S$ with $\phi(\za)=\zb$.
\end{lemma}

On the surface, all curves are considered up to homotopy, and all  intersections of curves are required to be transversal.
We fix the orientation of the surface, such that when drawing surfaces in the plane, when we follow  the boundary, the interior of the surface is on the right. For simplicity, we denote a marked surface without punctures by $(S,M)$.

\begin{definition}
\label{definition:addmissable dissections}
A collection of $\gpoint$-arcs~$\{\gamma_1, \ldots, \gamma_r\}$ is \emph{admissible} if the only possible intersections of these arcs are at the endpoints, and each subsurface enclosed by the arcs contains at least one $\rpoint$-point from $M^{\rpoint}\cup P^{\rpoint}$. A maximal admissible collection $\zD$ of~$\gpoint$-arcs is called an \emph{admissible $\gpoint$-dissection}. The notion of~\emph{admissible $\rpoint$-dissection} is defined in a similar way. For convenience, we will often simply write admissible dissection instead of admissible $\gpoint$-dissection.
\end{definition}
Note that an admissible dissection cuts the surface into polygons such that each polygon contains exactly one $\rpoint$-marked point. Furthermore, whenever we write $(S,M,P,\zD)$, we implicitly assume that $\zD$ is an admissible dissection on $(S,M,P)$.
\begin{example}\label{example:admissible dissection}
The pictures in Figure \ref{figure:admissible dissection} are two admissible dissections on an annulus.
\begin{figure}[ht]
\begin{center}
\begin{tikzpicture}[scale=0.3]
\begin{scope}
	\draw (0,0) circle (1cm);
	\clip[draw] (0,0) circle (1cm);
	\foreach \x in {-4,-3.5,-3,-2.5,-2,-1.5,-1,-0.5,0,0.5,1,1.5,2,2.5,3,3.5,4}	\draw[xshift=\x cm]  (-5,5)--(5,-5);
\end{scope}
	\draw (0,0) circle (4cm);
    \draw (0,-1) node {$\gpoint$};
    \draw (0,-4) node {$\gpoint$};
    \draw (0,4) node {$\gpoint$};
    \draw (-4,0) node {$\gpoint$};
    \draw (4,0) node {$\gpoint$};

    \draw (0,1) node {$\rpoint$};
    \draw (-3.5,2) node {$\rpoint$};
    \draw (3.5,2) node {$\rpoint$};
    \draw (-2,-3.5) node {$\rpoint$};
    \draw (2,-3.5) node {$\rpoint$};

    \draw (0.5,-2) node {\tiny$\zg_3$};
    \draw[](0,-4)to(0,-1);

    \draw[bend left](0,-1)to(-4,0);
    \draw[](0,-4)to(4,0);
    \draw plot [smooth,tension=1] coordinates {(0,-4) (2,1) (0,4)};
    \draw plot [smooth,tension=1] coordinates {(0,-1) (-1.5,0) (0,4)};

    \draw (2.5,-2.2) node {\tiny$\zg_5$};
    \draw (2.5,0) node {\tiny$\zg_4$};
    \draw (-2,1) node {\tiny$\zg_1$};
    \draw (-2.5,-.5) node {\tiny$\zg_2$};
    \draw (0,-5.5) node {$\zD_1$};
\end{tikzpicture}
\begin{tikzpicture}[scale=0.3]
\begin{scope}
	\draw (0,0) circle (1cm);
	\clip[draw] (0,0) circle (1cm);
	\foreach \x in {-4,-3.5,-3,-2.5,-2,-1.5,-1,-0.5,0,0.5,1,1.5,2,2.5,3,3.5,4}	\draw[xshift=\x cm]  (-5,5)--(5,-5);
\end{scope}
	\draw (0,0) circle (4cm);
    \draw (0,-1) node {$\gpoint$};
    \draw (0,-4) node {$\gpoint$};
    \draw (0,4) node {$\gpoint$};
    \draw (-4,0) node {$\gpoint$};
    \draw (4,0) node {$\gpoint$};

    \draw (0,1) node {$\rpoint$};
    \draw (-3.5,2) node {$\rpoint$};
    \draw (3.5,2) node {$\rpoint$};
    \draw (-2,-3.5) node {$\rpoint$};
    \draw (2,-3.5) node {$\rpoint$};

    \draw (0.5,-2) node {\tiny$\zg_3$};
    \draw[](0,-4)to(0,-1);

    \draw[](0,-4)to(4,0);
    \draw[](0,-4)to(-4,0);
    \draw plot [smooth,tension=1] coordinates {(0,-4) (-2,1) (0,4)};
    \draw plot [smooth,tension=1] coordinates {(0,-4) (2,1) (0,4)};
    \draw (2.5,-2.2) node {\tiny$\zg_5$};
    \draw (2.5,0) node {\tiny$\zg_4$};
    \draw (-2.4,1) node {\tiny$\zg_1$};
    \draw (-2.5,-.5) node {\tiny$\zg_2$};

    \draw (-6,-.5) node {};
    \draw (0,-5.5) node {$\zD_2$};
\end{tikzpicture}
\end{center}
\begin{center}
\caption{Two admissible dissections on an annulus.}\label{figure:admissible dissection}
\end{center}
\end{figure}
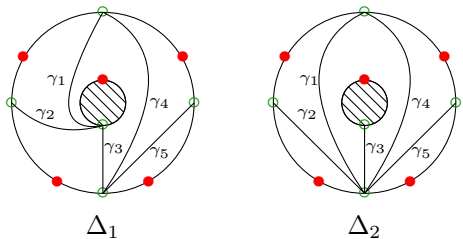
\end{example}

Denote by $g$ the genus of $S$ and by $b$ the number of connected components of $\partial S$. It is shown in \cite{APS19}, see also \cite{PPP18}, that an admissible collection of~$\gpoint$-arcs is an admissible dissection if and only if it is maximal, that is, it contains exactly  \( |M^{\gpoint}| -\chi \) arcs, where $\chi=2-2g-b-|P|$ is the Euler characteristic of the marked surface. Furthermore, let $\Delta$ be an admissible dissection. Then there exists a unique admissible $\rpoint$-dissection $\Delta^*$ (up to homotopy) such that each arc of $\Delta^*$ intersects exactly one arc of $\Delta$. We call $\Delta^*$ the \emph{dual $\rpoint$-dissection} of $\Delta$.
For a fixed admissible dissection $\Delta$ and its dual $\rpoint$-dissection $\Delta^*$, unless otherwise stated, all curves on the surface are always assumed to be in minimal position with respect to both these dissections.

A \emph{graded curve} $(\gamma,f)$ on $(S,M,P,\zD)$ is an admissible arc or a closed curve $\gamma$, together with a function
   \[
    f: \gamma \cap \Delta^* \longrightarrow \mathbb{Z},
   \]
  where~$\gamma \cap \Delta^*$ is the totally ordered set of intersection points of~$\gamma$ with~$\Delta^*$, where the order is induced by the direction of $\gamma$.
  The function $f$ is defined as follows:
  If~$p$ and~$q$ are in~$\gamma \cap \Delta^*$ and~$q$ is the direct successor of~$p$,
  then~$\gamma$ enters a polygon enclosed by~$\rpoint$-arcs of~$\Delta^*$ via~$p$ and leaves it via~$q$.
  If the unique $\gpoint$-point in this polygon is to the left of~$\gamma$, then~$f(q) = f(p)+1$; otherwise,~$f(q) = f(p)-1$.
Note that not all closed curves are gradable.

 For a grading $f$ of $\zg$, the map $f[n]: l \rightarrow f(l)-n$ with $n\in \mathbb{Z}$ is also a grading on $\zg$, and all gradings on $\zg$ are of this form.
We call  $[1]$ the \emph{shift} of the grading $f$.
The grading shift on curves corresponds to the shift functor in the derived category of the associated gentle algebra.

\subsection{Derived categories of gentle algebras}\label{subsection:derived categories of gentle algebras}
We recall in this section the derived category of a gentle algebra, and its geometric realization given in \cite{OPS18}  using a marked surface.

\begin{definition}\cite{AS87}\label{definition:gentle algebras}
We call an algebra $A $ a \emph{gentle algebra}, if $A$ is given by $kQ/I$, where $Q=(Q_0,Q_1)$ is a finite quiver and $I$ an admissible ideal of $kQ$ satisfying the following conditions:
\begin{enumerate}[\rm(1)]
 \item Each vertex in $Q_0$ is the source of at most two arrows and the target of at most two arrows.

 \item For each arrow $a$ in $Q_1$, there is  at most one arrow $b'$ such that $ab'\notin I$; at most one arrow $c'$ such that $c'a\notin I$; at most one arrow $b$ such that $ab$ is a path in $Q$ and $ab\in I$; at most one arrow $c$ such that $ca$ is a path in $Q$ and $ca\in I$.

 \item $I$ is generated by paths of length two.
\end{enumerate}
\end{definition}

For a finite dimensional algebra $A$, it is well known that the bounded derived category is triangle equivalent to the homotopy category $\mathsf{K^{-,b}}(\A)$ of complexes of projective $A$-modules bounded on the right and bounded in homology.  In the following, we will not distinguish these two categories. Denote by $\ka$ the full subcategory of $\mathsf{K^{-,b}}(\A)$, which is the homotopy category of bounded complexes of projective $A$-modules.

 The derived category of a gentle algebra is well-studied. In particular, the authors in \cite{BM03} classified the indecomposable objects in the category in terms of  \emph{(homotopy) string objects} and \emph{(homotopy) band objects}. The morphisms between these indecomposable objects are explicitly described in \cite{ALP16}.

Starting with an admissible collection $\zD$ on $(S,M,P)$, we define an algebra $A(\zD)=kQ(\zD)/I(\zD)$ as follows:
\begin{enumerate}[\rm(1)]
  \item The vertices of $Q(\zD)$ are given by the arcs in $\zD$.
  \item For any two arcs $\za$ and $\zb$, there is an arrow $a: \za\rightarrow \zb$ if they sharing an endpoint $q$ and $\zb$ directly follows $\za$ anticlockwise at $q$.
  \item The ideal $I(\zD)$ is generated by the following relations: whenever $\za$ and $\zb$ intersect at a marked point as above, and the other end of $\zb$ intersects $\zg$ at a marked point as above, then the composition $ab$ of the corresponding arrows $a: \za\rightarrow \zb$ and $b: \zb\rightarrow \zg$ is a relation.
\end{enumerate}
Then $A(\zD)$ is a gentle algebra. In particular, if $\zD$ is an admissible dissection, then this establishes a bijection between the set of homeomorphism classes of  marked surfaces $(S,M,P,\zD)$ and the set of isomorphism classes of gentle algebras $A(\zD)$, see  \cite{OPS18}, and also \cite{BC21}.

\begin{example}\label{example:quiver of admissible dissection}
The pictures in Figure \ref{figure:quiver of admissible dissection} are the quivers with relations associated to the two admissible dissections in Figure \ref{figure:admissible dissection}.
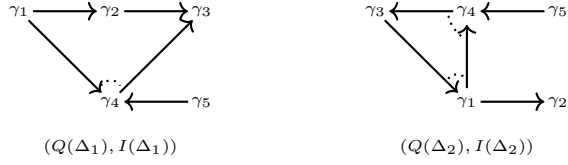
\begin{figure}[ht]
\begin{center}
{\begin{tikzpicture}[scale=0.6]

    \draw (0,0) node {\tiny$\zg_2$};
    \draw (2,0) node {\tiny$\zg_3$};
    \draw (-2,0) node {\tiny$\zg_1$};
    \draw (0,-2) node {\tiny$\zg_4$};
    \draw (2,-2) node {\tiny$\zg_5$};
    \draw [thick,->] (-1.7,0) -- (-.3,0);
    \draw [thick,->] (.3,0) -- (1.7,0);
    \draw [thick,<-] (.3,-2) -- (1.7,-2);
    \draw [thick,->] (-1.8,-0.2) -- (-0.2,-1.8);
    \draw [thick,<-] (1.8,-0.2) -- (0.2,-1.8);

    \draw[dotted,thick](-.2,-1.65) to [out=90,in=90] (.2,-1.65);
    \draw (0,-3) node {\tiny$(Q(\zD_1),I(\zD_1))$};

\end{tikzpicture}}
{\begin{tikzpicture}[scale=0.6]
      \draw (-5,-2) node {};

    \draw (0,0) node {\tiny$\zg_4$};
    \draw (2,0) node {\tiny$\zg_5$};
    \draw (-2,0) node {\tiny$\zg_3$};
    \draw (0,-2) node {\tiny$\zg_1$};
    \draw (2,-2) node {\tiny$\zg_2$};

    \draw [thick,<-] (-1.7,0) -- (-.3,0);
    \draw [thick,<-] (.3,0) -- (1.7,0);
    \draw [thick,->] (.3,-2) -- (1.7,-2);
    \draw [thick,->] (-1.8,-0.2) -- (-0.2,-1.8);
    \draw [thick,<-] (0,-0.3) -- (0,-1.7);

    \draw[dotted,thick](-.4,-1.5) to [out=90,in=90] (-.05,-1.5);
    \draw[dotted,thick](-.4,-.1) to [out=-90,in=-90] (-.05,-.5);
    \draw (0,-3) node {\tiny$(Q(\zD_2),I(\zD_2))$};

\end{tikzpicture}}
\end{center}
\begin{center}
\caption{The quivers with relations associated  to the two admissible dissections in Figure \ref{figure:admissible dissection}, where we denote the relations by dotted lines.}\label{figure:quiver of admissible dissection}
\end{center}
\end{figure}
\end{example}

\begin{remark}\label{remark: dissection and quiver}
 Similarly, one can associate an algebra $A(\zD^*)$ to an $\rpoint$-admissible collection $\zD^*$, which is the Koszul dual of $A(\zD)$ by \cite[Section 1.7]{OPS18}.
\end{remark}

Let $A$ be a gentle algebra associated with a marked surface $(S,M,P,\zD_A)$. It is shown in \cite{OPS18} that there is a correspondence between the indecomposable objects in  $\mathsf{K^{-,b}}(\A)$ and graded curves on the surface.
More precisely, the graded admissible arcs $(\zg,f)$ on $(S,M,P,\Delta_A)$
are in bijection with the isomorphism classes of the indecomposable string objects $\P_{(\zg,f)}$ in $\mathsf{K^{-,b}}(\A)$; the graded closed curves $(\zg,f)$ on $(S,M,P,\Delta_A)$ together with an indecomposable $k[x,x^{-1}]$-module $\mathfrak{m}$ are in bijection with the isomorphism classes of the indecomposable band objects $\P_{(\zg,f,\mathfrak{m})}$ in $\mathsf{K^{-,b}}(\A)$. Furthermore, under the above bijections, the indecomposable objects in $\ka$ correspond to the graded $\gpoint$-arcs and the pairs of graded closed curves and indecomposable $k[x,x^{-1}]$-modules.
For convenience of notation, we often will drop the $\mathfrak{m}$ from the notation, so that $P_{(\gamma, f)}$ denotes both a string or band object.

The morphisms in $\mathsf{K^{-,b}}(\A)$ are described in terms of graded oriented intersections, see  \cite[Theorem 3.3]{OPS18}. Roughly speaking, each intersection of two curves on the boundary (i.e. at a $\gpoint$-marked point) gives rise to exactly one morphism between the corresponding indecomposable objects. Each intersection in the interior (not at a $\rpoint$-puncture) of $S$ gives rise to two morphisms, while an intersection at a $\rpoint$-puncture gives rise to infinitely many   morphisms.
Since the geometric description of the morphisms is a crucial part of many of the proofs, we will now give a more detailed description.
For that, let $\P_{(\zg_1,f_1)}$ and $\P_{(\zg_2,f_2)}$ be two indecomposable objects in $\mathsf{K^{-,b}}(\A)$ corresponding to admissible graded arcs or graded closed curves $(\zg_1, f_1)$ and $(\zg_2, f_2)$.

Assume first that $\zg_1$ and $\zg_2$ intersect at a $\gpoint$-point $p$ on the boundary as depicted in Figure~\ref{boundary intersection}, where for $i\in \{1,2\}$, $q_i$ is the intersection in $\gamma_i \cap \Delta_A^*$ which is nearest to $p$.
If $f_1(q_1) = f_2(q_2)$, then there is a morphism from $\P_{(\zg_1,f_1)}$ to $\P_{(\zg_2,f_2)}$. Note that there is  no morphism from $\P_{(\zg_2,f_2)}$ to $\P_{(\zg_1,f_1)}[j]$ for any $j\in \mathbb{Z}$ arising from this intersection at $p$.

\begin{figure}[H]
	{\begin{tikzpicture}[scale=0.3]
	\draw[green] (0,0) circle [radius=0.1];
	\draw[](0.1,0)to(4.9,0);
	\draw[red](4.3,-0.5)to(4.6,0.5);
	\draw (2.5,.5) node {\tiny$\zg_2$};	
	\draw[](0.1,-0.05)to(4.93,-2.93);
	\draw[red](4.3,-0.5-2.6)to(4.8,0.5-2.9);
	\draw (2.5,-2.2) node {\tiny$\zg_1$};	
	\draw (4.5,-2.1) node {\tiny$q_1$};
	\draw (5.1,.6) node {\tiny$q_2$};	
	\draw (-1,0) node {\tiny$p$};			
	\draw[bend right](-.2,-.6)to(-.2,.6);
	\draw[-]  (-0.1,-.33)to(-.4,-.6);
	\draw[-]  (-.05,-.1)to(-.4,-.4);
	\draw[-]  (-.05,.1)to(-.4,-.2);
	\draw[-]  (-.05,.3)to(-.4,0);
	\draw[-]  (-.1,.45)to(-.4,0.2);
	\end{tikzpicture}}
	\caption{Boundary intersection of $\zg_1$ and $\zg_2$ giving rise to exactly one morphism.}~\label{boundary intersection}
	\end{figure}
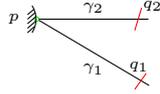

Assume now that $\zg_1$ and $\zg_2$ intersect at some point $p$ in the interior of the surface, which is not a $\rpoint$-puncture.  Suppose that we have the following local configuration in $S$ depicted in  Figure~\ref{interior intersection} below.
If $f_1(q_1) = f_2(q_2)$, then there
is one morphism from $\P_{(\zg_1,f_1)}$ to $\P_{(\zg_2,f_2)}$, and one morphism from $\P_{(\zg_2,f_2)}$ to $\P_{(\zg_1,f_1)}[1]$ (or to $\P_{(\zg_1,f_1)}[-1]$ depending on the position of the corresponding $\gpoint$-marked point).
	
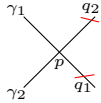
\begin{figure}[H]
		{\begin{tikzpicture}[scale=0.2]
			\draw[-] (0.15,0.15) -- (4.85,4.85);
			\draw[-] (0.15,4.85) -- (4.85,0.15);

			\draw[-,red] (3.5,0.8) -- (4.8,1);
			\draw[-,red] (3.9,4.8) -- (5.1,4.5);
			
			\draw (-.3,5.3) node {\tiny$\zg_1$};
			\draw (-.3,-.3) node {\tiny$\zg_2$};
			
			\draw (4.6,5.4) node {\tiny$q_2$};
			\draw (4.1,0.2) node {\tiny$q_1$};
			\draw (2.5,1.7) node {\tiny$p$};
			\end{tikzpicture}}
			\caption{Interior intersection of $\zg_1$ and $\zg_2$ giving rise to two morphisms in opposite directions.}\label{interior intersection}
			\end{figure}

The Auslander-Reiten translation $\tau$ in $\ka$ is  described in \cite{OPS18} in terms of arcs in the surface. We will rephrase the result  by introducing a slightly different notion of  twisting an arc as follows. The \emph{direct twist} (resp. \emph{inverse twist}) of an $\gpoint$-arc $\zg$ in $(S,M,P)$ is an $\rpoint$-arc $D\zg$ (resp. $D^{-1}\zg$) obtained from $\zg$ by rotating both endpoints anticlockwise (resp. clockwise) to the next  $\rpoint$-point. Similarly, we define the \emph{direct twist} $D\zg^*$ and the \emph{inverse twist} $D^{-1}\zg^*$ of an $\rpoint$-arc $\zg^*$ in $(S,M,P)$. Note that there is a canonical intersection between $\zg$ and $D^2\zg$ arising from the rotation.

\begin{lemma}\cite[Corollary 5.2]{OPS18}\label{lemma:tau and Dehn}
Let $\P_{(\zg,f)}$ be an indecomposable string object in $\ka$ corresponding to a graded $\gpoint$-arc $(\zg,f)$. Then $\tau^{-1}\P_{(\zg,f)}=\P_{(D^{2}\zg,f_{D^{2}\zg})}$, where $f_{D^{2}\zg}$ is a grading on $D^{2}\zg$ which is uniquely determined by the fact that the canonical intersection between $\zg$ and $D^2\zg$ gives rise to a map from $\P_{(D^{2}\zg,f_{D^{2}\zg})}$ to $\P_{(\zg,f)}[1]$.
\end{lemma}

\subsection{Silting objects and exceptional sequences}\label{subsection: silting objects}
We now recall some background on silting theory and the theory of exceptional sequences in a general triangulated category as well as for the particular case of  the bounded derived categories of gentle algebras.

Let $\calt$ be a triangulated category.
We call a full subcategory $\calp$ in $\calt$ \emph{pre-silting} if for all $i>0$, $\Hom_{\calt}(\calp,\calp[i])=0$. It is \emph{silting} if in addition $\calt=\thick\calp$. An object $P$ of $\calt$ is said to be \emph{pre-silting} if $\add P$ is a pre-silting subcategory and \emph{silting} if $\add P$ is a silting subcategory. We always assume that pre-silting objects as well as silting objects are basic.

In the triangulated category $\calt$, we call an object $X\in\calt$ \emph{exceptional} if
$\Hom_{\calt}(X,X[\neq0])=0$ and $\End_{\calt}(X)$ is a division algebra.
We call an (ordered) sequence $(X_1,\cdots,X_n)$ of exceptional objects in $\calt$ an \emph{exceptional sequence} if
$\Hom_{\calt}(X_i,X_j[\mathbb{Z}])=0\ \mbox{ for any }\ 1\le j<i\le n$,
which is said to be \emph{full} if in addition
$\thick_\calt(\bigoplus_{i=1}^nX_i)=\calt$.
Note that the grading does not have any bearing  on whether a sequence is exceptional or not. In general, there is an action of  $\Z^n$  on the set of full exceptional sequences $\exc\calt$ in $\calt$ given as follows
\[(\ell_1,\cdots,\ell_n)(X_1,\cdots,X_n):=(X_1[\ell_1],\cdots,X_n[\ell_n]).\]

The following proposition shows that there are close relations between exceptional sequences and silting objects.
\begin{proposition}[Proposition 3.5 \cite{AI12}]\label{prop:exceptional to silting}
Let $\calt$ be a totally Hom-finite triangulated category, that is, for any $X,Y\in\calt$ we have $\Hom_{\calt}(X,Y[\ell])=0$ for any $|\ell|\gg0$. Let $(X_1,\cdots,X_n)$ be an exceptional sequence in $\calt$. Then

 (1) there exists $a\in\mathbb{Z}$ such that $X_1[\ell_1]\oplus\cdots\oplus X_n[\ell_n]$
is a pre-silting object for any integers $\ell_1\cdots,\ell_n\in\mathbb{Z}$
satisfying $\ell_i+a\le\ell_{i+1}$, for all $1\le i<n$.

(2) the exceptional sequence $(X_1,\cdots,X_n)$ is a  full exceptional sequence if and only if $X_1[\ell_1]\oplus\cdots\oplus X_n[\ell_n]$
is a silting object.
\end{proposition}

The category $\ka$ always has silting objects. The category $\mathsf{K^{-,b}}(\A)$ contains silting objects if and only if the global dimension of $A$ is finite, that is,  $\mathsf{K^{-,b}}(\A)$ is  triangle equivalent to $\ka$. For a gentle algebra $A$, the silting objects in $\ka$ are described in \cite[Theorem 3.2]{APS19}. Namely, any basic silting object $X$ in $\ka$ is of the form
$\P_{(\zD, f)}=\bigoplus_{i=1}^n \P_{(\gamma_i, f_{\zg_i})}$,
where $\zD=\{\gamma_1, \ldots, \gamma_n\}$ is an admissible~dissection of~$(S,M,P)$ and $f$ is a set of gradings $f_{\zg_i}$ on $\zg_i, 1 \leq i \leq n$.

Note that, in general, for an arbitrary admissible~dissection $\zD$, there currently is no characterisation of when there exists a set of gradings $f$ over $\zD$ such that $(\zD,f)$ gives rise to a silting object.

\section{Exceptional sequences in terms of ordered exceptional collections}\label{subsection:exceptional sequences in terms of surfaces dissections}

In this section we describe exceptional sequences in terms of surface dissections. Throughout this section, let $A$ be a gentle algebra associated with a marked surface $(S,M,P,\zD_A)$. Let $\zD^*_A$ be the dual admissible $\rpoint$-dissection of $\zD_A$.

We begin by describing indecomposable exceptional objects by the surface.

\begin{lemma}\label{lemma:braid group}
Let $\P_{(\zg,f)}$ be an indecomposable object in $\ka$ arising from a graded curve $(\zg,f)$.

(1) Then $\P_{(\zg,f)}$ is exceptional if and only if $\zg$ is an $\gpoint$-arc without self-intersections. In particular, $\zg$ is not a loop.

(2) If $\zg$ is an $\gpoint$-arc without self-intersections, then for any grading $f$ of $\zg$, $\P_{(\zg,f)}$ is an indecomposable exceptional object in $\ka$.
\end{lemma}

\begin{proof}
(1) Assume $\P_{(\zg,f)}$ is exceptional. Note that if $\zg$ is a closed curve, then $\P_{(\zg,f)}$ is a band object which has self-extension and so cannot be exceptional.
If $\zg$ is an infinite arc, that is at least one end of $\zg$ wraps infinitely many times around a puncture,   then following \cite[Remark 3.8]{OPS18} there exist infinitely many maps from $\P_{(\zg,f)}$ to  positive shifts of $\P_{(\zg,f)}$. So $\zg$ can only be an $\gpoint$-arc. Moreover, if $\zg$ has an interior self-intersection, then we have two linearly independent maps arising from this intersection, and $\Hom(\P_{(\zg,f)},\P_{(\zg,f)}[m])\neq 0$ for some $m\in \mathbb{Z}^*$. Therefore $\P_{(\zg,f)}$ is not exceptional. On the other hand, if $\zg$ is a loop without interior self-intersection, then we have two cases. If the two gradings near the endpoint are equal, then we have two linearly independent maps from $\P_{(\zg,f)}$ to itself: the identity and the map corresponding to the intersection. Then the endomorphism algebra of $\P_{(\zg,f)}$ is isomorphic to $k[x]/(x^2)$, which is not a division algebra. Otherwise, if the gradings are different, then there is some $m \in \mathbb{Z}^*$ such that there is a non-zero element in $\Hom(\P_{(\zg,f)},\P_{(\zg,f)}[m])$ arising from the intersection. In both cases $\P_{(\zg,f)}$ is not exceptional.

Conversely, it is clear that if $\zg$ is an $\gpoint$-arc without self-intersections, then $\Hom(\P_{(\zg,f)},\P_{(\zg,f)}[\neq0])=0$ and $\End(\P_{(\zg,f)})\cong k$ is a division ring, so that $\P_{(\zg,f)}$ is exceptional.

(2) This is clear by the first statement and an observation that shifting the grading of objects does not change the fact whether an object is exceptional or not.
\end{proof}

Following the above lemma, we will call an $\gpoint$-arc an \emph{exceptional arc} if it has no self-intersections.

\begin{definition}\label{definition:tree type}
 Let $\zD$ be an admissible $\gpoint$-collection on $(S,M,P)$ consisting of exceptional arcs. We call $\zD$ an {\it exceptional collection}, if it has no \emph{non-exceptional cycle}, that is a subgraph consisting of arcs in $\zD$ as depicted in Figure \ref{figure:Non exceptional}, where $\zg_{i+1}$ (not necessarily directly) follows $\zg_i$ in the anticlockwise order at the endpoint $q_i$ and where the index $1\leq i \leq m$ is considered modulo $m$. We call an exceptional collection an \emph{exceptional dissection} if it is an admissible dissection.
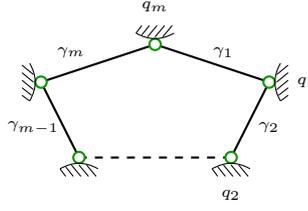
\begin{figure}
 \[\scalebox{1}{
\begin{tikzpicture}[>=stealth,scale=0.5]
\draw[thick] (0,0)--(-1,2)--(2,3)--(5,2)--(4,0);
\draw[thick,dashed](4,0)--(0,0);
\draw[dark-green,thick,fill=white] (0,0) circle (0.15);
\draw[dark-green,thick,fill=white] (-1,2) circle (0.15);
\draw[dark-green,thick,fill=white] (2,3) circle (0.15);
\draw[dark-green,thick,fill=white] (5,2) circle (0.15);
\draw[dark-green,thick,fill=white] (4,0) circle (0.15);
\node at (3.8,2.8) {\tiny$\zg_{1}$};
\node at (-.2,2.8) {\tiny$\zg_{m}$};
\node at (5,.8) {\tiny$\zg_{2}$};
\node at (-1.2,.8) {\tiny$\zg_{m-1}$};
\node at (6,2) {\tiny$q_1$};
\node at (4,-1) {\tiny$q_2$};
\node at (2,4) {\tiny$q_{m}$};

	\draw[bend left](5.3,1.5)to(5.3,2.5);
	\draw[-]  (5.5,1.65)to(5.3,1.5);
	\draw[-]  (5.5,1.85)to(5.25,1.65);
	\draw[-]  (5.5,2.05)to(5.2,1.8);
	\draw[-]  (5.5,2.25)to(5.2,2);
	\draw[-]  (5.5,2.45)to(5.2,2.2);

	\draw[bend right](-5.3+4,1.5)to(-5.3+4,2.5);
	\draw[-]  (-5.5+4,1.65)to(-5.3+4,1.5);
	\draw[-]  (-5.5+4,1.85)to(-5.25+4,1.65);
	\draw[-]  (-5.5+4,2.05)to(-5.2+4,1.8);
	\draw[-]  (-5.5+4,2.25)to(-5.2+4,2);
	\draw[-]  (-5.5+4,2.45)to(-5.2+4,2.2);

	\draw[bend right](1.5, 3.3)to(2.5, 3.3);
	\draw[-]  (1.65, 3.5)to(1.5, 3.3);
	\draw[-]  (1.85, 3.5)to(1.65, 3.25);
	\draw[-]  (2.05, 3.5)to(1.8, 3.2);
	\draw[-]  (2.25, 3.5)to(2, 3.2);
	\draw[-]  (2.45, 3.5)to(2.2, 3.2);

	\draw[bend left](3.5, -.3)to(4.5, -.3);
	\draw[-]  (3.65, -.5)to(3.5, -.3);
	\draw[-]  (3.85, -.5)to(3.65, -.25);
	\draw[-]  (4.05, -.5)to(3.8, -.2);
	\draw[-]  (4.25, -.5)to(4, -.2);
	\draw[-]  (4.45, -.5)to(4.2, -.2);

	\draw[bend left](-.5, -.3)to(0.5, -.3);
	\draw[-]  (3.65-4, -.5)to(3.5-4, -.3);
	\draw[-]  (3.85-4, -.5)to(3.65-4, -.25);
	\draw[-]  (4.05-4, -.5)to(3.8-4, -.2);
	\draw[-]  (4.25-4, -.5)to(4-4, -.2);
	\draw[-]  (4.45-4, -.5)to(4.2-4, -.2);
\end{tikzpicture}}\]
\begin{center}
\caption{A non-exceptional cycle}\label{figure:Non exceptional}
\end{center}
\end{figure}
\end{definition}

The quiver $Q(\zD)$ associated to an admissible collection $\zD$ gives an easy way of checking whether this collection is exceptional or not.

\begin{lemma}\label{lemma:quvier of exceptional dissection}
An admissible collection $\zD$ is an exceptional collection if and only if there exists no oriented cycle in the associated quiver $Q(\zD)$.
\end{lemma}

\begin{proof}
Clearly every non-exceptional cycle in $\zD$ gives rise to an oriented cycle in $Q(\zD)$. On the other hand suppose that we have an oriented cycle $a_1 a_2 \ldots a_m$ in $Q(\zD)$. Then the subgraph of $\zD$ corresponding to all vertices $t(a_i)$ of $Q(\zD)$ such that $a_i a_{i+1} \in I$,  considering the indexing set  modulo $m$, form a non-exceptional cycle in $\zD$. Thus $\zD$ is not an exceptional collection.
\end{proof}

\begin{example}\label{example:exceptional dissection}
In Figure \ref{figure:admissible dissection}, $\zD_1$ is an exceptional dissection, while $\zD_2$ is not. This can be seen from the  quivers of the two algebras in Figure \ref{figure:quiver of admissible dissection}: the quiver $Q(\zD_1)$ does not have an oriented cycle whereas the quiver $Q(\zD_2)$ has an oriented cycle.
\end{example}

Let $\zD=\{\zg_i, 1 \leq i \leq n\}$ be an exceptional collection. We define a partial order on $\zD$ as follows. For two arcs $\zg_i, \zg_j \in \zD$, set $\zg_i\preceq \zg_j$ if $\zg_i$ and $\zg_j$  share an endpoint $q$ and if $\zg_i$ follows $\zg_j$ in the clockwise order at $q$ or if $\zg_i=\zg_j$.
We consider the transitive closure of this relation which we will still denote by $\preceq$.  Then it is clear that $\preceq$ is reflexive.
On the other hand, it follows from the definition of an exceptional collection that if $\zg_i \preceq \zg_j$ and $\zg_i \neq \zg_j$, then $\zg_j \npreceq \zg_i$, since there are no `oriented' subgraphs as depicted above. Thus this order relation is antisymmetric, that is $\zg_i \preceq \zg_j$ and $\zg_j \preceq \zg_i$ implies $\zg_i = \zg_j$.
Therefore $\preceq$ is a well-defined partial order on $\zD$.

\begin{definition}\label{definition:tree type}
We call an ordered set of arcs $(\zg_1,\cdots,\zg_m)$ in $(S, M, P)$ an \emph{ordered exceptional collection} if it is an exceptional collection and the order of the arcs are compatible with the partial order $\preceq$ introduced above, that is, $\zg_i\preceq \zg_j$ implies $i \leq j$. We call an ordered exceptional collection an \emph{ordered exceptional dissection} if it is an admissible dissection. We will use $\zG$ to denote an ordered exceptional dissection. We define $\zG^*=(\zg^*_m,\cdots,\zg^*_1)$ as the \emph{dual ordered exceptional dissection} of $\zG$. Note that we reverse the indices of arcs in $\zG^*$.
\end{definition}

While there is no geometric criterion to determine which admissible dissections give rise to silting objects, we can give such a criterion to determine the admissible dissections which correspond to full exceptional sequences as follows.

\begin{proposition}\label{prop:tree type}
(1) If $(X_1,\cdots,X_n)$ is a full exceptional sequence in $\ka$, then $X_i=\P_{(\zg_i,f_i)}$ for some graded $\gpoint$-arcs $(\zg_i,f_i)$ such that $\Delta=(\zg_1,\cdots,\zg_n)$ is an ordered exceptional dissection of $(S,M,P)$.

(2) Let $\Delta=(\zg_1,\cdots,\zg_n)$ be an ordered exceptional dissection of $(S,M,P)$. Then for any grading $f_i$ over $\zg_i, 1\leq i \leq n$, $(\P_{(\zg_1,f_1)},\cdots,\P_{(\zg_n,f_n)})$ is a full exceptional sequence in $\ka$.
\end{proposition}
\begin{proof}
(1) Since each $X_i$, $1\leq i\leq n$, is an indecomposable exceptional object, by Lemma \ref{lemma:braid group}, we may assume that $X_i=\P_{(\zg_i,f_i)}$ for some graded arc $(\zg_i,f_i)$, where $\zg_i$ is an $\gpoint$-arc without self-intersections.
On the other hand, since $\ka$ is totally Hom-finite, there exists integers $\ell_1\cdots,\ell_n\in\mathbb{Z}$ such that $X_1[\ell_1]\oplus\cdots\oplus X_n[\ell_n]$
is a silting object by Proposition \ref{prop:exceptional to silting}.
Thus $\zD=\{\zg_1,\cdots,\zg_n\}$ is an admissible dissection of $(S,M,P)$.
Suppose $\zD$ is not an exceptional dissection. Then there is a cycle formed by arcs $\zg_{i_1}, \zg_{i_2}, \cdots, \zg_{i_m}$ in $\zD$. Without loss of generality, we may assume that this is a cycle as follows.

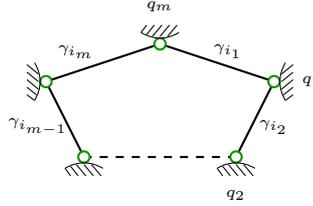
\[\scalebox{1}{
\begin{tikzpicture}[>=stealth,scale=0.5]
\draw[thick] (0,0)--(-1,2)--(2,3)--(5,2)--(4,0);
\draw[thick,dashed](4,0)--(0,0);
\draw[dark-green,thick,fill=white] (0,0) circle (0.15);
\draw[dark-green,thick,fill=white] (-1,2) circle (0.15);
\draw[dark-green,thick,fill=white] (2,3) circle (0.15);
\draw[dark-green,thick,fill=white] (5,2) circle (0.15);
\draw[dark-green,thick,fill=white] (4,0) circle (0.15);
\node at (3.8,2.8) {\tiny$\zg_{i_1}$};
\node at (-.2,2.8) {\tiny$\zg_{i_m}$};
\node at (5,.8) {\tiny$\zg_{i_2}$};
\node at (-1.2,.8) {\tiny$\zg_{i_{m-1}}$};
\node at (6,2) {\tiny$q_1$};
\node at (4,-1) {\tiny$q_2$};
\node at (2,4) {\tiny$q_{m}$};

	\draw[bend left](5.3,1.5)to(5.3,2.5);
	\draw[-]  (5.5,1.65)to(5.3,1.5);
	\draw[-]  (5.5,1.85)to(5.25,1.65);
	\draw[-]  (5.5,2.05)to(5.2,1.8);
	\draw[-]  (5.5,2.25)to(5.2,2);
	\draw[-]  (5.5,2.45)to(5.2,2.2);

	\draw[bend right](-5.3+4,1.5)to(-5.3+4,2.5);
	\draw[-]  (-5.5+4,1.65)to(-5.3+4,1.5);
	\draw[-]  (-5.5+4,1.85)to(-5.25+4,1.65);
	\draw[-]  (-5.5+4,2.05)to(-5.2+4,1.8);
	\draw[-]  (-5.5+4,2.25)to(-5.2+4,2);
	\draw[-]  (-5.5+4,2.45)to(-5.2+4,2.2);

	\draw[bend right](1.5, 3.3)to(2.5, 3.3);
	\draw[-]  (1.65, 3.5)to(1.5, 3.3);
	\draw[-]  (1.85, 3.5)to(1.65, 3.25);
	\draw[-]  (2.05, 3.5)to(1.8, 3.2);
	\draw[-]  (2.25, 3.5)to(2, 3.2);
	\draw[-]  (2.45, 3.5)to(2.2, 3.2);

	\draw[bend left](3.5, -.3)to(4.5, -.3);
	\draw[-]  (3.65, -.5)to(3.5, -.3);
	\draw[-]  (3.85, -.5)to(3.65, -.25);
	\draw[-]  (4.05, -.5)to(3.8, -.2);
	\draw[-]  (4.25, -.5)to(4, -.2);
	\draw[-]  (4.45, -.5)to(4.2, -.2);

	\draw[bend left](-.5, -.3)to(0.5, -.3);
	\draw[-]  (3.65-4, -.5)to(3.5-4, -.3);
	\draw[-]  (3.85-4, -.5)to(3.65-4, -.25);
	\draw[-]  (4.05-4, -.5)to(3.8-4, -.2);
	\draw[-]  (4.25-4, -.5)to(4-4, -.2);
	\draw[-]  (4.45-4, -.5)to(4.2-4, -.2);
\end{tikzpicture}}\]
Then  we have
$$\Hom(\P_{(\zg_{i_j},f_{{i_j}})},
\P_{(\zg_{i_{j+1}},f_{{i_{j+1}}})}[\Z]) \neq 0,$$
for $1 \leq j \leq m-1$,
 and
$$\Hom(\P_{(\zg_{i_m},f_{{i_m}})},
\P_{(\zg_{i_{1}},f_{{i_{1}}})}[\Z]) \neq 0.$$
This contradicts that $(X_1, \ldots, X_n)$ is an exceptional sequence.
Therefore $\zD$ is an exceptional collection, and thus an exceptional dissection.

(2) Let $\Delta=(\zg_1,\cdots,\zg_n)$ be an ordered exceptional dissection of $(S,M,P)$. Since the order of $\Delta$ is compatible with the partial order $\preceq$ of the arcs on $(S,M,P)$, $i < j$ implies that $\zg_i\preceq \zg_j$ or $\zg_i$ and $\zg_j$ are not comparable by $\preceq$. If $\zg_i\preceq \zg_j$, then $\Hom(\P_{(\zg_{j},f_{{j}})},
\P_{(\zg_{i},f_{{i}})}[\Z])=0$ for any grading $f_{i}$ and $f_{j}$. On the other hand, if $\zg_i$ and $\zg_j$ are not comparable by $\preceq$, then in particular, there exists no intersections between $\zg_i$ and $\zg_j$, and thus $\Hom(\P_{(\zg_{j},f_{{j}})},
\P_{(\zg_{i},f_{{i}})}[\Z])=0$.
Therefore $\Hom(\P_{(\zg_{j},f_{{j}})},
\P_{(\zg_{i},f_{{i}})}[\Z])=0$ when $i < j$.
Furthermore, since  $\zg_i$ has no self-intersection, $\P_{(\zg_i,f_i)}$ is an exceptional object for all $1 \leq i \leq n$ by Lemma \ref{lemma:braid group}.
Thus $(\P_{(\zg_{i_1},f_{{i_1}})},\P_{(\zg_{i_2},f_{{i_2}})},
\cdots,\P_{(\zg_{i_n},f_{{i_n}})})$ is an exceptional sequence, which is full since $\zD$ is maximal.
\end{proof}

We have the following immediate corollary.

\begin{corollary}\label{coro:tree type}
(1) There exists a full exceptional sequence in $\ka$ if and only if there exists an exceptional dissection on $(S,M,P)$.

(2) Let $A_1$ and $A_2$ be two gentle algebras with the same surface model, then there exists a full exceptional sequence in $\mathsf{K^b(A_1)}$ if and only if there exists a full exceptional sequence in $\mathsf{K^b(A_2)}$.
\end{corollary}

\section{The existence of full exceptional sequences}\label{Section:The existence}
In this section, we consider the existence of full exceptional sequences in $\ka$ by considering the existence of exceptional dissections on the associated surface. We begin by considering the case of surfaces containing one or two boundary components and one or two marked points.

(1) $\mathbb{T}_{(g,1,1)}$ is the marked surface of genus  $g \geq 1$ with only one boundary component and  exactly one $\gpoint$-marked point.

(2) $\mathbb{T}_{(g,1,2)}$ is the marked surface of genus  $g \geq 1$ with only one boundary component and  exactly two $\gpoint$-marked point.

(3) $\mathbb{T}_{(g,2,2)}$ is the marked surface of genus  $g \geq 1$ with  two boundary components each of which has one $\gpoint$-marked point.

\begin{lemma}\label{lemma:existence0}
There exists no exceptional dissection in $\mathbb{T}_{(g,1,1)}$.
\end{lemma}
\begin{proof}
Note that if there exists only one $\gpoint$-marked point on the surface, then any $\gpoint$-arc is a loop, and thus can not be exceptional. Therefore there exists no exceptional dissections in $\mathbb{T}_{(g,1,1)}$.
\end{proof}

\begin{lemma}\label{lemma:existence1}
There exist exceptional dissections in both $\mathbb{T}_{(g,1,2)}$ and $\mathbb{T}_{(g,2,2)}$.
\end{lemma}
\begin{proof}
On  $\mathbb{T}_{(g,1,2)}$ and $\mathbb{T}_{(g,2,2)}$ there are admissible dissections such that the associated gentle algebras are of the form as in Figure \ref{figure:Tg12} and Figure \ref{figure:Tg22}, respectively.
\begin{figure}[ht]
\begin{center}
{\begin{tikzpicture}[scale=0.35]
\def \radius {4cm}
    \draw [thick,->] (0,.2) -- (2,.2);
    \draw [thick,->] (.5+3,.2) -- (2.5+3,.2);
    \draw [thick,->] (.5+9,.2) -- (2.5+9,.2);
    \draw [thick,->] (.5+12,.2) -- (2.5+12,.2);
    \draw [thick,->] (.5+15,.2) -- (2.5+15,.2);

    \draw [thick,->] (0,-.2) -- (2,-.2);
    \draw [thick,->] (.5+3,-.2) -- (2.5+3,-.2);
    \draw [thick,->] (.5+9,-.2) -- (2.5+9,-.2);
    \draw [thick,->] (.5+12,-.2) -- (2.5+12,-.2);
    \draw [thick,->] (.5+15,-.2) -- (2.5+15,-.2);

    \draw [thick,dotted] (.5+6,0) -- (2.5+6,0);

    \draw[dotted,thick](1.5,.5) to [out=90,in=90] (3.5,.5);
    \draw[dotted,thick](1.5,-.5) to [out=-90,in=-90] (3.5,-.5);
    \draw[dotted,thick](1.5+3.4,.5) to [out=90,in=90] (3.5+3.4,.5);
    \draw[dotted,thick](1.5+3.4,-.5) to [out=-90,in=-90] (3.5+3.4,-.5);
    \draw[dotted,thick](1.5+3.4+3,.5) to [out=90,in=90] (3.5+3.4+3,.5);
    \draw[dotted,thick](1.5+3.4+3,-.5) to [out=-90,in=-90] (3.5+3.4+3,-.5);
    \draw[dotted,thick](1.5+3.4+3+3,.5) to [out=90,in=90] (3.5+3.4+3+3,.5);
    \draw[dotted,thick](1.5+3.4+3+3,-.5) to [out=-90,in=-90] (3.5+3.4+3+3,-.5);
    \draw[dotted,thick](1.5+3.4+3+3+3,.5) to [out=90,in=90] (3.5+3.4+3+3+3,.5);
    \draw[dotted,thick](1.5+3.4+3+3+3,-.5) to [out=-90,in=-90] (3.5+3.4+3+3+3,-.5);

    \draw (-1.2,0) node {\tiny$2g+1$};
    \draw (2.7,0) node {\tiny$2g$};
    \draw (12,0) node {\tiny$3$};
    \draw (15,0) node {\tiny$2$};
    \draw (18,0) node {\tiny$1$};
\end{tikzpicture}}
\end{center}
\begin{center}
\caption{The quiver of an exceptional dissection on $\mathbb{T}_{(g,1,2)}$.}\label{figure:Tg12}
\end{center}
\end{figure}
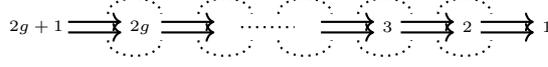

\begin{figure}[ht]
\begin{center}
{\begin{tikzpicture}[scale=0.35]
\def \radius {4cm}
    \draw [thick,->] (-1,.2) -- (1,.2);
    \draw [thick,->] (.5+3,.2) -- (2.5+3,.2);
    \draw [thick,->] (.5+9,.2) -- (2.5+9,.2);
    \draw [thick,->] (.5+12,.2) -- (2.5+12,.2);
    \draw [thick,->] (.5+15,.2) -- (2.5+15,.2);

    \draw [thick,->] (-1,-.2) -- (1,-.2);
    \draw [thick,->] (.5+3,-.2) -- (2.5+3,-.2);
    \draw [thick,->] (.5+9,-.2) -- (2.5+9,-.2);
    \draw [thick,->] (.5+12,-.2) -- (2.5+12,-.2);
    \draw [thick,->] (.5+15,-.2) -- (2.5+15,-.2);

    \draw [thick,dotted] (.5+6,0) -- (2.5+6,0);

    \draw[dotted,thick](.6,.5) to [out=90,in=90] (3.5,.5);
    \draw[dotted,thick](.6,-.5) to [out=-90,in=-90] (3.5,-.5);
    \draw[dotted,thick](1.5+3.4,.5) to [out=90,in=90] (3.5+3.4,.5);
    \draw[dotted,thick](1.5+3.4,-.5) to [out=-90,in=-90] (3.5+3.4,-.5);
    \draw[dotted,thick](1.5+3.4+3,.5) to [out=90,in=90] (3.5+3.4+3,.5);
    \draw[dotted,thick](1.5+3.4+3,-.5) to [out=-90,in=-90] (3.5+3.4+3,-.5);
    \draw[dotted,thick](1.5+3.4+3+3,.5) to [out=90,in=90] (3.5+3.4+3+3,.5);
    \draw[dotted,thick](1.5+3.4+3+3,-.5) to [out=-90,in=-90] (3.5+3.4+3+3,-.5);
    \draw[dotted,thick](1.5+3.4+3+3+3,.5) to [out=90,in=90] (3.5+3.4+3+3+3,.5);
    \draw[dotted,thick](1.5+3.4+3+3+3,-.5) to [out=-90,in=-90] (3.5+3.4+3+3+3,-.5);

    \draw (-2.5,0) node {\tiny$2g+2$};
    \draw (2.3,0) node {\tiny$2g+1$};
    \draw (12,0) node {\tiny$3$};
    \draw (15,0) node {\tiny$2$};
    \draw (18,0) node {\tiny$1$};
\end{tikzpicture}}
\end{center}
\begin{center}
\caption{The quiver of an exceptional dissection on $\mathbb{T}_{(g,2,2)}$.}\label{figure:Tg22}
\end{center}
\end{figure}
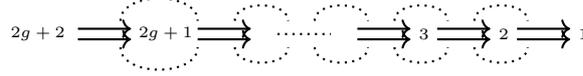

The existence of such admissible dissections follows from the  computation of the ribbon graphs of these quivers with relations, see   \cite{S15, OPS18}, and from showing that the surfaces given by the ribbon graphs are exactly the surfaces $\mathbb{T}_{(g,1,2)}$ and $\mathbb{T}_{(g,2,2)}$, respectively. On the other hand, since the quivers have no oriented cycles, by Lemma \ref{lemma:quvier of exceptional dissection}, the associated admissible dissections are exceptional dissections on $\mathbb{T}_{(g,1,2)}$ and $\mathbb{T}_{(g,2,2)}$, respectively.
\end{proof}

In Proposition \ref{prop:quivers of two surfaces}, we will show that any quiver with relations arising from an exceptional dissection on $\mathbb{T}_{(g,1,2)}$ and $\mathbb{T}_{(g,2,2)}$ is isomorphic to the quiver with relations depicted in Figure \ref{figure:Tg12} and Figure \ref{figure:Tg22} respectively.
The following lemma describes a numerical property of the quiver arising from these three special surfaces.

\begin{lemma}\label{lemma:surface and quiver}
Let $A=kQ/I$ be a gentle algebra arising from $(S,M,P,\zD_A)$. Then

(1) $(S,M)$ is of the form $\mathbb{T}_{(g,1,1)}$   if and only if $(|Q_0|,|Q_1|)=(2n,4n-1)$, for some $n \geq 1$.

(2) $(S,M)$ is of the form $\mathbb{T}_{(g,1,2)}$ if and only if $(|Q_0|,|Q_1|)=(2n+1,4n)$ for some $n \geq 1$.

(3) $(S,M)$ is of the form $\mathbb{T}_{(g,2,2)}$ if and only if $(|Q_0|,|Q_1|)=(2n+2,4n+2)$ for some $n \geq 1$.

Furthermore, if $(1), (2)$ or $(3)$ holds then $n=g$.
\end{lemma}
\begin{proof}
For a marked surface $(S,M)$, denote by $\chi=2-2g-b$ the Euler characteristic of $S$, where $g$ is the genus of $S$ and $b$ is the number of boundary components of $S$. Then we have the following equalities:
\begin{equation}\label{equation:1}
	|Q_0|=|\zD|=|M^{\gpoint}|-\chi
\end{equation}
\begin{equation}\label{equation:2}
	\chi=|Q_0|-|Q_1|,
\end{equation}
see for example in \cite{APS19,LP20}.

It follows directly from these equalities that if $(S,M)=\mathbb{T}_{(g,1,1)}$, then $(|Q_0|,|Q_1|)=(2g,4g-1)$; if $(S,M)=\mathbb{T}_{(g,1,2)}$, then $(|Q_0|,|Q_1|)=(2g+1,4g)$; if $(S,M)=\mathbb{T}_{(g,2,2)}$, then $(|Q_0|,|Q_1|)=(2g+2,4g+2)$.

Now we show the converse.

(1) If $(|Q_0|,|Q_1|)=(2n,4n-1)$ for some $n\geq 1$, then by  equality \eqref{equation:2}, the Euler characteristic of $S$ equals $1-2n$. Thus by  equality \eqref{equation:1} there exists only one $\gpoint$-marked point on $(S,M)$. Furthermore, there is only one boundary component on $S$, since there exists at least one $\gpoint$-marked point on each boundary component. Therefore $(S,M)$ is of the form $\mathbb{T}_{(g,1,1)}$, and we have $g=n$.

(2) If $(|Q_0|,|Q_1|)=(2n+1,4n)$ for some $n\geq 1$, then by  equality \eqref{equation:2}, the Euler characteristic of $S$ equals $1-2n$. Thus by  equality \eqref{equation:1} there are exactly two $\gpoint$-marked points on $(S,M)$. Hence $S$ has only one or two boundary components. On the other hand, note that since $2-2g-b=\chi=1-2n$, we have that  $b$ is odd. Therefore there is only one boundary component on $S$. So $(S,M)$ is of the form $\mathbb{T}_{(g,1,1)}$, and $g=n$.

(3) The proof is similar to the proof of the second part.
\end{proof}

Before giving equivalent conditions for the existence of full exceptional dissections on a marked surface, we show the following technical lemma.

\begin{lemma}\label{lemma:adding} Let $(S,M,P)$ be a marked surface. Let $(S',M',P)$ be a marked surface obtained from $(S,M,P)$ by adding a $\gpoint$-point and a $\rpoint$-point on some boundary component or by adding a new boundary component with exactly one $\gpoint$-point and one $\rpoint$-point. If there exists an exceptional dissection on $(S,M,P)$ then there  exists an exceptional dissection on $(S',M',P)$.
\end{lemma}

\begin{proof} Let $\zg$ be a  minimal element in an exceptional dissection $\zD$ on $(S,M,P)$  with respect to the partial order introduced in Section \ref{subsection:exceptional sequences in terms of surfaces dissections}.  Denote the new  $\gpoint$-marked point by $p$ and the new $\rpoint$-marked point by $q$.
The  two generic  possible ways to add $p$ and $q$ is as in Figure \ref{figure:add marked points}.
\begin{figure}[ht]
\begin{center}
\begin{tikzpicture}[scale=0.22]
\begin{scope}
	\draw (0,0) circle (3cm);
	\clip[draw] (0,0) circle (3cm);
	\foreach \x in {-4,-3.5,-3,-2.5,-2,-1.5,-1,-0.5,0,0.5,1,1.5,2,2.5,3,3.5,4}	\draw[xshift=\x cm]  (-5,5)--(5,-5);
\end{scope}
    \draw (0,3) node {$\gpoint$};
    \draw (3,0) node {$\gpoint$};
    \draw (2.1,2.1) node {$\rpoint$};
    \draw (2.5,2.6) node {\tiny$q$};
        \draw[bend right](0,3)to(-4,4);
        \draw (-3,3.5) node {\tiny$\zg$};
        \draw (3.7,0) node {\tiny$p$};
        \draw (3.5,3.6) node {\tiny$\za$};
    \draw plot [smooth,tension=1] coordinates {(0,3) (3,3) (3,0)};
        \draw (0,-4) node {};
\end{tikzpicture}
\begin{tikzpicture}[scale=0.2]
\begin{scope}
	\draw (0,0) circle (1cm);
	\clip[draw] (0,0) circle (1cm);
	\foreach \x in {-2.5,-2,-1.5,-1,-0.5,0,0.5,1,1.5,2,2.5}	\draw[xshift=\x cm]  (-3,3)--(3,-3);
\end{scope}
    \draw (0,1) node {$\rpoint$};
    \draw (0,1.9) node {\tiny$q$};
    \draw (0,-1) node {$\gpoint$};

\begin{scope}
	\draw (0,-6) circle (1cm);
	\clip[draw] (0,-6) circle (1cm);
	\foreach \x in {-2.5,-2,-1.5,-1,-0.5,0,0.5,1,1.5,2,2.5}	\draw[xshift=\x cm]  (-3,3-6)--(3,-3-6);
\end{scope}
    \draw (0,-5) node {$\gpoint$};
    \draw[-]  (0,-5)to(0,-1);
    \draw[bend right](0,-5)to(-4,-4);
    \draw plot [smooth,tension=1] coordinates {(0,-5) (-2,-.5) (0,2.7) (2,1) (1.5,-1) (0,-1)};
        \draw (-3,-4.5) node {\tiny$\zg$};
        \draw (-3,0) node {\tiny$\za_2$};
        \draw (1,-3) node {\tiny$\za_1$};
        \draw (-.5,-1.7) node {\tiny$p$};

        \draw (-9,-1.5) node {};
\end{tikzpicture}
\end{center}
\begin{center}
\caption{Local configuration of the new surface $(S',M',P)$ obtained from a surface $(S,M,P)$ by adding marked points or a boundary with marked points. In the left case we add a $\gpoint$-marked point $p$ and a $\rpoint$-marked point $q$ on a boundary component of $S$, while in the second case we add a new boundary component (upper one) with one $\gpoint$-marked point $p$ and one $\rpoint$-marked point $q$ on $(S,M,P)$.}\label{figure:add marked points}
\end{center}
\end{figure}
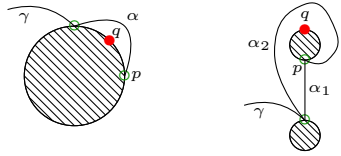

We now construct an admissible collection of arcs $\zD'$ on $(S',M',P)$ by adding new arcs to $\zD$ as shown in Figure \ref{figure:add marked points}, where we add an arc $\za$ as in  the case of the left of Figure \ref{figure:add marked points} and we add two arcs $\za_1$ and $\za_2$ as in  the case of the right  of  Figure \ref{figure:add marked points}.
 Then  $\zD'$ is  an admissible dissection on $(S',M',P)$, since it is maximal, that is it has  precisely $|\zD|+1$ arcs in the first case and $|\zD|+2$ arcs in the second case.

Furthermore, note that the quiver $Q(\zD')$ is obtained from $Q(\zD)$ by adding an arrow from $\za$ to $\zg$ in the first case, and by adding two arrows from $\za_1$ to $\za_2$ and one arrow from $\za_2$ to $\zg$ in the second case. On the other hand, since $\zD$ is exceptional, there is no oriented cycle in $Q(\zD)$ by Lemma \ref{lemma:quvier of exceptional dissection}. Thus there is no oriented cycle in $Q(\zD')$. So $\zD'$ is exceptional again by Lemma \ref{lemma:quvier of exceptional dissection}.
\end{proof}

\begin{proposition}\label{thm:existence of full exceptional sequence}
Let $(S,M,P)$ be a marked surface. Then the following are equivalent:
\begin{enumerate}[\rm(1)]
\item There exists an exceptional dissection on $(S,M,P)$;
\item $P=\emptyset$ and $|M^{\gpoint}|\geq 2$;
 \item  $P=\emptyset$ and $(S,M,P)$ is not  homeomorphic to  $\mathbb{T}_{(g,1,1)}$.
\end{enumerate}

\end{proposition}
\begin{proof}
Clearly (2) and (3) are equivalent.  We now show that  (1) implies (3).
Suppose $P\neq \emptyset$ in $(S,M,P)$.
Let $\zD$ be any admissible dissection of $(S,M,P)$. Then there  exists a polygon  in $\zD$ containing a $\rpoint$-puncture and this polygon is a non-exceptional cycle in $\zD$. Thus $\zD$ can not be an exceptional dissection.
Furthermore, by Lemma \ref{lemma:existence0}, a surface  of the form $\mathbb{T}_{(g,1,1)}$  has no exceptional dissection.

We now show that (2) implies (1).
So assume that $P=\emptyset$ and $|M^{\gpoint}|\geq 2$. If the genus of $S$ is greater than one then the surfaces  $\mathbb{T}_{(g,1,2)}$ and $\mathbb{T}_{(g,2,2)}$ are the marked surfaces with minimal number of $\gpoint$-marked points containing at least two marked points. Lemma  \ref{lemma:existence1} shows that there exists at least one exceptional dissection on both $\mathbb{T}_{(g,1,2)}$ and $\mathbb{T}_{(g,2,2)}$. Then the result follows from Lemma~\ref{lemma:adding} by induction on the number of marked points of the surface. If the genus of $S$ is zero, we can also use an inductive argument   based on iteratively using Lemma \ref{lemma:adding}, and noticing that the disk with two $\gpoint$-marked points is the minimal surface with genus zero, which  has an exceptional dissection.
\end{proof}

The following result shows that if there is an exceptional dissection on $(S,M)$ then there are infinitely many, unless $S$ is a disk.

\begin{proposition}\label{prop:number}
Let $(S,M)$ be a marked surface with an exceptional dissection. Then the number of exceptional dissections is finite if and only if $S$ is a disk.
\end{proposition}
\begin{proof}
Recall that the twist of an $\gpoint$-arc $\zg$ in $(S,M)$ is an $\rpoint$-arc $D\zg$ obtained from $\zg$ by  rotating anticlockwise both endpoints to the next respective $\rpoint$-point.
Then the result follows from the following two observations: any square of a twist of an exceptional dissection is again an exceptional dissection; the group generated by homeomorphisms of twisting each boundary component of $(S,M)$ which maps a $\gpoint$-marked point to the (anticlockwise) next $\gpoint$-marked point is finite if and only $S$ is a disk. Furthermore, if $S$ is not a disk then the action of this group is faithful. \end{proof}

From Proposition~\ref{thm:existence of full exceptional sequence} and  Proposition \ref{prop:tree type} (2), we now obtain a complete description of the existence of full exceptional sequences in $\ka$ for a gentle algebra $A$.

\begin{theorem}\label{thm:existence}
Let $A$ be a gentle algebra with marked surface $(S,M,P,\zD_A)$. The following are equivalent:
 \begin{enumerate}[\rm(1)]
 \item There exists a full exceptional sequence in $\ka$; \item $P=\emptyset$ and $|M^{\gpoint}|\geq 2$;
 \item $P=\emptyset$ and $(S,M,P)$ is not homeomorhpic to $\mathbb{T}_{(g,1,1)}$.
 \end{enumerate}
\end{theorem}

By the above Theorem and  the fact that the global dimension of a gentle algebra is finite if and only if the associated marked surface has no punctures, we have the following.
\begin{corollary}\label{cor:existence1}
Let $A$ be a gentle algebra. If the global dimension of $A$ is infinite then there exists no full exceptional sequence in $\ka$. \end{corollary}

 In particular, noting that the surface model of a non-hereditary discrete-derived algebra is an annulus, Theorem \ref{thm:existence} gives a new proof of the following result on the existence of full exceptional sequences for the derived category of derived-discrete algebras which was first proved in \cite[Proposition 6.6]{BPP17}. We refer the reader to \cite{V01,BPP17} for details on derived-discrete algebras.

\begin{corollary}\label{cor:existence2}
Let $A$ be a non-hereditary derived-discrete algebra. Then there exists a full exceptional sequence in $\ka$.
\end{corollary}

We now give a criterion to determine the existence of full exceptional sequences in $\ka$ of a gentle algebra $A$ by directly looking at the quiver and the relations of the algebra.

\begin{corollary}\label{thm:quiver}
Let $A=kQ/I$ be a gentle algebra, then there exists a full exceptional sequence in $\ka$ if and only if the following two conditions are satisfied
\begin{enumerate}[\rm(1)]
 \item  the pair $(|Q_0|,|Q_1|)$ is not of the form $(2g,4g-1)$ for any $g \geq 1$;
 \item there exist no oriented cycles in $Q$ such that the composition of any two neighboring arrows belongs to $I$.
\end{enumerate}
\end{corollary}
\begin{proof}
 Lemma \ref{lemma:surface and quiver} (1) says that a gentle algebra $A=kQ/I$ arising from an admissible dissection on a surface of the form $\mathbb{T}_{(g,1,1)}$ if and only if $(|Q_0|,|Q_1|)=(2g,4g-1)$.  Furthermore, the global dimension of a gentle algebra is finite if and only if there exist no oriented cycles with relations at each vertex. This is further equivalent to the surfaces having not punctures, that is $P=\emptyset$.
Then the result follows from Theorem \ref{thm:existence of full exceptional sequence}.
\end{proof}

\section{Cutting  surfaces and completing exceptional sequences}\label{subsection:cut surface}

We now define the notion of a cut surface which will play an important role in the inductive arguments later in the paper.
\begin{definition}\label{definition:cut surface1}
Let $\zg$ be an exceptional $\gpoint$-arc on a marked surface $(S,M)$ which has two distinct endpoints $p$ and $q$.
We define $(S_\zg , M_\zg)$ to be the marked surface obtained from $(S,M)$ by cutting along $\zg$
and then contracting along the cut, and we set $M_\zg=M\setminus \{p,q\}\sqcup \{pq,p'q'\}$, where $pq$ and $p'q'$ are new marked points obtained from $p$ and $q$ after cutting and  contracting along $\zg$.
 We call $(S_\zg , M_\zg)$ the \emph{cut surface} of $(S,M)$ along $\zg$.
\end{definition}

Note that whenever we consider a cut surface $(S_\gamma, M_\gamma)$, we implicitly assume that $\gamma$ is an exceptional $\gpoint$-arc. When cutting the surface, we always ignore components corresponding to disks with only one $\gpoint$-point and one $\rpoint$-point. If the arc $\zg$ is a separating arc, then the cut surface $(S_\zg , M_\zg)$ may not be  connected anymore and $(S_\zg , M_\zg)$ is in fact a union of marked surfaces.

\begin{example}\label{example:cut surface}
Figure \ref{figure:cut surface} gives an example of a cut surface. Namely, we cut an annulus along an arc $\zg$ connecting the two boundary components and then contract along the cut. From this we obtain the cut surface which is a disk.

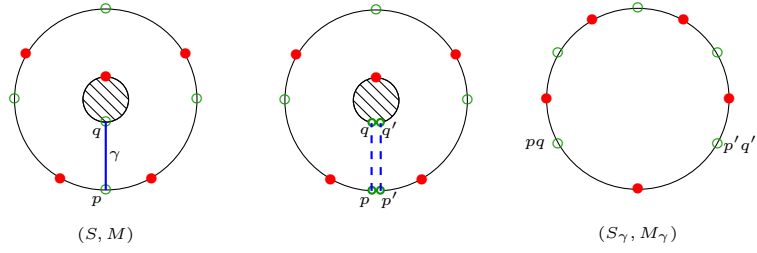
\begin{figure}[ht]
\begin{center}
\begin{tikzpicture}[scale=0.3]
\begin{scope}
	\draw (0,0) circle (1cm);
	\clip[draw] (0,0) circle (1cm);
	\foreach \x in {-4,-3.5,-3,-2.5,-2,-1.5,-1,-0.5,0,0.5,1,1.5,2,2.5,3,3.5,4}	\draw[xshift=\x cm]  (-5,5)--(5,-5);
\end{scope}
	\draw (0,0) circle (4cm);
    \draw (0,-1) node {$\gpoint$};
    \draw (0,-4) node {$\gpoint$};
    \draw (0,4) node {$\gpoint$};
    \draw (-4,0) node {$\gpoint$};
    \draw (4,0) node {$\gpoint$};

    \draw (0,1) node {$\rpoint$};
    \draw (-3.5,2) node {$\rpoint$};
    \draw (3.5,2) node {$\rpoint$};
    \draw (-2,-3.5) node {$\rpoint$};
    \draw (2,-3.5) node {$\rpoint$};

    \draw (0.4,-2.5) node {\tiny$\zg$};
    \draw[blue,thick](0,-4)to(0,-1);
%

    \draw (-.4,-4.5) node {\tiny$p$};
    \draw (-.4,-1.5) node {\tiny$q$};

    \draw (0,-6) node {\tiny$(S,M)$};
\end{tikzpicture}
\begin{tikzpicture}[scale=0.3]
\begin{scope}
	\draw (0,0) circle (1cm);
	\clip[draw] (0,0) circle (1cm);
	\foreach \x in {-4,-3.5,-3,-2.5,-2,-1.5,-1,-0.5,0,0.5,1,1.5,2,2.5,3,3.5,4}	\draw[xshift=\x cm]  (-5,5)--(5,-5);
\end{scope}
	\draw (0,0) circle (4cm);
    \draw[dark-green,thick,fill=white] (-.185,-.97) circle (0.15);
    \draw[dark-green,thick,fill=white] (.185,-.97) circle (0.15);
    \draw[dark-green,thick,fill=white] (-.185,-3.97) circle (0.15);
    \draw[dark-green,thick,fill=white] (.185,-3.97) circle (0.15);

    \draw (0,4) node {$\gpoint$};
    \draw (-4,0) node {$\gpoint$};
    \draw (4,0) node {$\gpoint$};

    \draw (0,1) node {$\rpoint$};
    \draw (-3.5,2) node {$\rpoint$};
    \draw (3.5,2) node {$\rpoint$};
    \draw (-2,-3.5) node {$\rpoint$};
    \draw (2,-3.5) node {$\rpoint$};

    \draw[blue,thick,dashed](-.2,-4)to(-.2,-1);
    \draw[blue,thick,dashed](.2,-4)to(.2,-1);

    \draw (-.5,-4.5) node {\tiny$p$};
    \draw (-.5,-1.4) node {\tiny$q$};
    \draw (.6,-4.4) node {\tiny$p'$};
    \draw (.6,-1.3) node {\tiny$q'$};

    \draw (-6.4,-6.3) node {};
\end{tikzpicture}
\begin{tikzpicture}[scale=0.3]
    \draw (-6,-.5) node {};
	\draw (0,0) circle (4cm);
    \draw (0,-4) node {$\rpoint$};
    \draw (0,4) node {$\gpoint$};
    \draw (-4,0) node {$\rpoint$};
    \draw (4,0) node {$\rpoint$};
    \draw (-2,3.5) node {$\rpoint$};
    \draw (2,3.5) node {$\rpoint$};

    \draw (3.5,-2) node {$\gpoint$};
    \draw (3.5,2) node {$\gpoint$};
    \draw (-3.5,-2) node {$\gpoint$};
    \draw (-3.5,2) node {$\gpoint$};

    \draw (-4.5,-2) node {\tiny$pq$};
    \draw (4.5,-2) node {\tiny$p'q'$};
%
%
%
    \draw (0,-6) node {\tiny$(S_\zg,M_\zg)$};
\end{tikzpicture}
\end{center}
\begin{center}
\caption{An example of cut surface.}\label{figure:cut surface}
\end{center}
\end{figure}
\end{example}

Let $\za$ be an $\gpoint$-arc in $(S,M)$. It is easy to see that
if $\za$ does not intersect $\zg$, then $\za$ is still an $\gpoint$-arc in $S_\zg$ and there are no other arcs in $S$ which are identified with $\za$ in $(S_\zg,M_\zg)$. While if $\za$ and $\zg$ intersect in the interior of $S$ then $\za$ disappears after cutting, and it does not give rise to a curve in $(S_\zg,M_\zg)$ anymore.
When $\za$ and $\zg$ only intersect at endpoints, the situation becomes more complicated. In particular, in that case, two distinct arcs in $S$ might be identified in $S_\gamma$,  see for example in Figure \ref{figure:possible positions between a and g}, $\za_1$ and $\za_2$ are identified in $(S_\zg,M_\zg)$.

More precisely, for an $\gpoint$-arc $\za$ in $(S,M)$, we denote by $\overline{\za}$ the $\gpoint$-arc in $(S_\zg,M_\zg)$ induced by $\za$. For example $\overline{\za_1}=\overline{\za_2}$ for $\za_1$ and $\za_2$ in Figure \ref{figure:possible positions between a and g}. In particular, if $\za$ intersects $\zg$ in the interior, then we view $\overline{\za}$ as an empty element. For example, $\overline{\zg}$ is empty.  We refer to \cite[Proposition 4.6]{CS20} for a complete description of when two  $\gpoint$-arcs are identified in the cut surface.

\begin{figure}[H]
\begin{center}
{\begin{tikzpicture}[scale=0.3]
\draw[green] (0,0) circle [radius=0.2];
\draw[green] (0,5) circle [radius=0.2];
\draw[green] (5,0) circle [radius=0.2];

\draw[-] (0,0.2) -- (0,4.8);

\draw[-,blue] (0.2,0) -- (4.8,0);

\draw[-] (0.15,4.85) -- (4.85,0.15);

\draw (2.3,-.5) node {$\za_1$};

\draw (-1.2,2.5) node {$\zg$};
\draw (3.5,2.5) node {$\za_2$};

\draw (-.6,5.5) node {$q_2$};
\draw (-.5,-.5) node {$q_1$};
\draw (5.5,-.5) node {$q_3$};
\end{tikzpicture}}
\end{center}
\begin{center}
\caption{An example of how arcs are identified in the cut surface, where on the cut surface $(S_\zg,M_\zg)$, we have $\overline{\za_1}=\overline{\za_2}$.}
\label{figure:possible positions between a and g}
\end{center}
\end{figure}
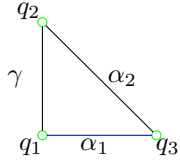

The following two lemmas show that any admissible (resp. exceptional) collection on $(S,M)$ induces an admissible (resp. exceptional) collection on $(S_\zg,M_\zg)$.
\begin{lemma}\label{lemma:completion}
Let $\zD$ be an admissible collection on $(S,M)$ and  $\zg$ an exceptional arc in $\zD$. Denote respectively by $Q$ and $I$ the quiver and the set of relations associated to $\zD$. Then $\overline{\zD}:=\{\overline{\za}, \za \in \zD\}$ is an admissible collection on $(S_\zg,M_\zg)$, whose quiver $\overline{Q}$ and the set of relations $\overline{I}$ are obtained from $Q$ and $I$ as follows:
 \begin{enumerate}[\rm(1)]
 \item The vertex set $\overline{Q}_0$ is given by  $Q_{0}\backslash \{ \zg\}$. The set $\overline{Q}_1$ of arrows is the union
 $$\{ c: s\ra t \in Q_{1}\mid s\not=\zg \text{ and } t\not=\zg\}~\sqcup~\{ [c_1c_2]=j \xra{c_1} \zg\xra{c
 _2} k \mid c_1 c_2\in I\}.$$

 \item $\overline{I}$ is generated by
\[ \{c_1c_2\mid c_1c_2\in I, c_1, c_2 \in Q_1\cap\overline{Q}_1\}~\sqcup~\{c_0[c_1c_2]\mid c_0c_1\in I \}~\sqcup~\{[c_1c_2]c_3\mid c_2c_3\in I \}.\]
 \end{enumerate}
Furthermore, if $\zD$ is an admissible dissection, then $\overline{\zD}$ is also an admissible dissection.
\end{lemma}

\begin{proof}
The proof is straightforward.
\end{proof}

\begin{lemma}\label{lemma:completione}
 Assuming the notations in Lemma \ref{lemma:completion}, we have that  if $\zD$ is an exceptional collection, then $\overline{\zD}$ is an exceptional collection. Furthermore, if $\zD$ is an exceptional dissection, then $\overline{\zD}$ also is  an exceptional dissection.
\end{lemma}

\begin{proof}
By Lemma \ref{lemma:completion}, we know that $\overline{\zD}$ is an admissible collection on $(S_\zg,M_\zg)$.
Now we show that $\overline{\zD}$ is also exceptional.
By Lemma \ref{lemma:quvier of exceptional dissection}, we only need to show that there is no oriented cycle in $\overline{Q}$. Assume for contradiction that $\omega$ is an oriented cycle in $\overline{Q}$. Then $\omega$ must contain a new arrow of the form $[c_1c_2]$ described in Lemma \ref{lemma:completion}, otherwise $\omega$ would already be an oriented cycle in $Q$, which contradicts to the fact that $\zD$ is exceptional. But then after replacing all the new arrows $[c_1c_2]$ in $\omega$ by the path $c_1c_2$, we have an oriented cycle in $Q$, which is again a contradiction. Thus there is no oriented cycle in $\overline{Q}$, and thus $\overline{\zD}$ is exceptional.

 Now we prove the second statement. Let $\za\neq \zb$ be two arcs in $\zD\setminus \{\zg\}$, then $\overline{\za}\neq \overline{\zb}$ since otherwise $\za$, $\zb$ and $\zg$ form a triangle on the surface which has no $\rpoint$-point, which contradicts the fact that $\zD$ is admissible. On the other hand, note that if $\overline{\za}$ is empty for an arc $\za$ in $\zD$, then $\za=\zg$. Thus $\overline{\zD}$ has exact $|\zD|-1$ arcs. Therefore by \cite[Proposition 1.11]{APS19}, $\overline{\zD}$ is a maximal admissible collection on $(S_\zg,M_\zg)$, that is, $\overline{\zD}$ is an exceptional dissection.
\end{proof}

Now we discuss how to lift an admissible (resp. exceptional) collection on $(S_\zg,M_\zg)$ to an admissible (resp. exceptional) collection on $(S,M)$. Let $\zD_\zg$ be a collection of arcs on $(S_\zg,M_\zg)$. Each arc on $(S_\zg,M_\zg)$ can be viewed as an arc induced by an arc (which may not be unique) on $(S,M)$. Thus there always exists a collection of arcs $\zD=\{\zg_1,\cdots,\zg_n\}$ on $(S,M)$ such that $\zD_\zg=\overline{\zD}$. So we can always denote a collection of arcs on $(S_\zg,M_\zg)$ by $\zD_\zg=\{\overline{\za}, \za\in \zD\}$ for some collection of arcs $\zD$ on $(S,M)$, and we will always assume that $\zD$ contains $\zg$, noticing that $\overline{\zg}$ is empty. We call $\zD$ a \emph{lift} of $\zD_\zg$.

Let $\zD_\zg=\{\overline{\zg_i}, 1 \leq i \leq n, \zg_i\in \zD\}$
be an admissible collection on $(S_\zg,M_\zg)$ where $\zD$ is a lift of $\zD_\zg$.
Note that the collection $\zD$ is not necessarily admissible, see Figure \ref{figure:non-admisssible lift} for an example.
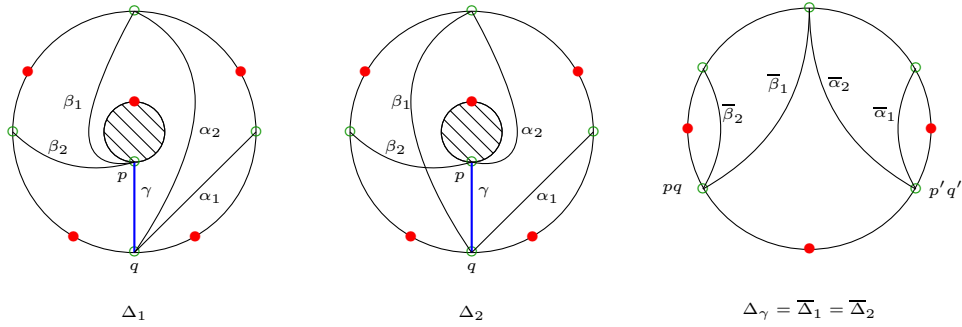
\begin{figure}[ht]
\begin{center}
\begin{tikzpicture}[scale=0.4]
\begin{scope}
	\draw (0,0) circle (1cm);
	\clip[draw] (0,0) circle (1cm);
	\foreach \x in {-4,-3.5,-3,-2.5,-2,-1.5,-1,-0.5,0,0.5,1,1.5,2,2.5,3,3.5,4}	\draw[xshift=\x cm]  (-5,5)--(5,-5);
\end{scope}
	\draw (0,0) circle (4cm);
    \draw (0,-1) node {$\gpoint$};
    \draw (0,-4) node {$\gpoint$};
    \draw (0,4) node {$\gpoint$};
    \draw (-4,0) node {$\gpoint$};
    \draw (4,0) node {$\gpoint$};

    \draw (0,1) node {$\rpoint$};
    \draw (-3.5,2) node {$\rpoint$};
    \draw (3.5,2) node {$\rpoint$};
    \draw (-2,-3.5) node {$\rpoint$};
    \draw (2,-3.5) node {$\rpoint$};

    \draw (0.4,-2) node {\tiny$\zg$};
    \draw[blue,thick](0,-4)to(0,-1);

    \draw[bend left](0,-1)to(-4,0);
    \draw[](0,-4)to(4,0);
    \draw plot [smooth,tension=1] coordinates {(0,-4) (2,1) (0,4)};
    \draw plot [smooth,tension=1] coordinates {(0,-1) (-1.5,0) (0,4)};

    \draw (0,-4.5) node {\tiny$q$};
    \draw (-.4,-1.5) node {\tiny$p$};
    \draw (2.5,-2.2) node {\tiny$\za_1$};
    \draw (2.5,0) node {\tiny$\za_2$};
    \draw (-2,1) node {\tiny$\zb_1$};
    \draw (-2.5,-.5) node {\tiny$\zb_2$};
    \draw (0,-6) node {\tiny$\zD_1$};
\end{tikzpicture}
\begin{tikzpicture}[scale=0.4]
\begin{scope}
	\draw (0,0) circle (1cm);
	\clip[draw] (0,0) circle (1cm);
	\foreach \x in {-4,-3.5,-3,-2.5,-2,-1.5,-1,-0.5,0,0.5,1,1.5,2,2.5,3,3.5,4}	\draw[xshift=\x cm]  (-5,5)--(5,-5);
\end{scope}
	\draw (0,0) circle (4cm);
    \draw (0,-1) node {$\gpoint$};
    \draw (0,-4) node {$\gpoint$};
    \draw (0,4) node {$\gpoint$};
    \draw (-4,0) node {$\gpoint$};
    \draw (4,0) node {$\gpoint$};

    \draw (0,1) node {$\rpoint$};
    \draw (-3.5,2) node {$\rpoint$};
    \draw (3.5,2) node {$\rpoint$};
    \draw (-2,-3.5) node {$\rpoint$};
    \draw (2,-3.5) node {$\rpoint$};

    \draw (0.4,-2) node {\tiny$\zg$};
    \draw[blue,thick](0,-4)to(0,-1);

    \draw[bend left](0,-1)to(-4,0);
    \draw[](0,-4)to(4,0);
    \draw plot [smooth,tension=1] coordinates {(0,-4) (-2,1) (0,4)};
    \draw plot [smooth,tension=1] coordinates {(0,-1) (1.5,0) (0,4)};

    \draw (0,-4.5) node {\tiny$q$};
    \draw (-.4,-1.5) node {\tiny$p$};
    \draw (2.5,-2.2) node {\tiny$\za_1$};
    \draw (2,0) node {\tiny$\za_2$};
    \draw (-2.3,1) node {\tiny$\zb_1$};
    \draw (-2.5,-.5) node {\tiny$\zb_2$};

    \draw (-6,-.5) node {};
    \draw (0,-6) node {\tiny$\zD_2$};
\end{tikzpicture}
\begin{tikzpicture}[scale=0.4]
    \draw (-6,-.5) node {};
	\draw (0,0) circle (4cm);
    \draw (0,-4) node {$\rpoint$};
    \draw (0,4) node {$\gpoint$};
    \draw (-4,0) node {$\rpoint$};
    \draw (4,0) node {$\rpoint$};

    \draw (3.5,-2) node {$\gpoint$};
    \draw (3.5,2) node {$\gpoint$};
    \draw (-3.5,-2) node {$\gpoint$};
    \draw (-3.5,2) node {$\gpoint$};

    \draw (-4.5,-2) node {\tiny$pq$};
    \draw (4.5,-2) node {\tiny$p'q'$};

    \draw (2.5,.5) node {\tiny$\overline{\za}_1$};
    \draw (1,1.5) node {\tiny$\overline{\za}_2$};
    \draw (-1,1.5) node {\tiny$\overline{\zb}_1$};
    \draw (-2.5,.5) node {\tiny$\overline{\zb}_2$};

    \draw[bend left](0,4)to(-3.5,-2);
    \draw[bend right](-3.5,-2)to(-3.5,2);

    \draw[bend right](0,4)to(3.5,-2);
    \draw[bend left](3.5,-2)to(3.5,2);
    \draw (0,-6) node {\tiny$\zD_\zg=\overline{\zD}_1=\overline{\zD}_2$};
\end{tikzpicture}
\end{center}
\begin{center}
\caption{Two lifts $\zD_1$ and $\zD_2$ of an admissible dissection $\zD_\zg$ on $(S_\zg,M_\zg)$, where $\zD_1$ is an admissible dissection on $(S,M)$ while $\zD_2$ is not.}\label{figure:non-admisssible lift}
\end{center}
\end{figure}
We will give conditions which guarantee that $\zD$ is admissible.
For this we denote by $(\overline{\zb}_1,\cdots,\overline{\zb}_v)$ and $(\overline{\za}_1,\cdots,\overline{\za}_t)$, respectively,  the ordered subsets of $\zD_\zg$  of all arcs with endpoints at  $pq$, respectively at $p'q'$.
Here by our conventions the orders of $(\overline{\zb}_1,\cdots,\overline{\zb}_v)$ and $(\overline{\za}_1,\cdots,\overline{\za}_t)$ are compatible with the clockwise orders at the endpoints, that is, $\overline{\zb}_1 \preceq \cdots \preceq \overline{\zb}_v$ and $\overline{\za}_1 \preceq \cdots \preceq \overline{\za}_t$, see the left picture in Figure \ref{figure:lift}.

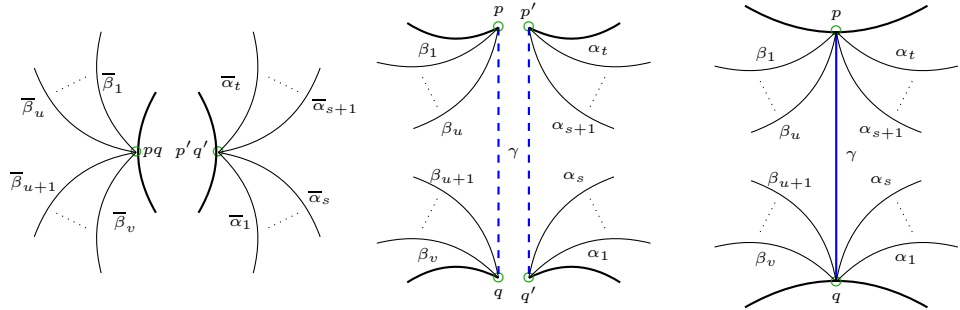
\begin{figure}[ht]
\begin{center}
\begin{tikzpicture}[scale=0.4]
    \draw (0,-5) node {};
    \draw[bend left,thick](4.8,-2)to(4.8,2);
    \draw (4.7,0) node {\tiny$pq$};
    \draw (6,0.1) node {\tiny$p'q'$};
    \draw (4.15,0) node {$\gpoint$};

    \draw[bend right](4.15,0)to(3,-4);
    \draw[bend right](4.15,0)to(.8,-2.8);
    \draw[dotted](2.5,-2.5)to(1.5,-2);
    \draw (3.8,-2.3) node {\tiny$\overline{\zb}_{v}$};
    \draw (.8,-1) node {\tiny$\overline{\zb}_{u+1}$};

    \draw (-4.15+11,0) node {$\gpoint$};
    \draw[bend left](4.15,0)to(3,4);
    \draw[bend left](4.15,0)to(.8,2.8);
    \draw[dotted](2.5,2.5)to(1.5,2);
    \draw (3.4,2.3) node {\tiny$\overline{\zb}_{1}$};
    \draw (.8,1.5) node {\tiny$\overline{\zb}_{u}$};

\draw[bend right,thick](-4.8+11,-2)to(-4.8+11,2);

    \draw[bend left](-4.15+11,0)to(-3+11,-4);
    \draw[bend left](-4.15+11,0)to(-.8+11,-2.8);
    \draw[dotted](-2.5+11,-2.5)to(-1.5+11,-2);
    \draw (-3.4+11,-2.3) node {\tiny$\overline{\za}_1$};
    \draw (-.8+11,-1.5) node {\tiny$\overline{\za}_s$};

    \draw[bend right](-4.15+11,0)to(-3+11,4);
    \draw[bend right](-4.15+11,0)to(-.8+11,2.8);
    \draw[dotted](-2.5+11,2.5)to(-1.5+11,2);
    \draw (-3.7+11,2.3) node {\tiny$\overline{\za}_{t}$};
    \draw (10.7,1.5) node {\tiny$\overline{\za}_{s+1}$};

\end{tikzpicture}
\begin{tikzpicture}[scale=0.4]

    \draw[bend right,thick](-3.5,4.3)to(-.5,4.15);
    \draw[bend left,thick](3.5,4.3)to(.5,4.15);
    \draw[bend right,thick](3.5,-4.3)to(.5,-4.15);
    \draw[bend left,thick](-3.5,-4.3)to(-.5,-4.15);

    \draw (-.5,4.15) node {$\gpoint$};
    \draw (-.5,-4.15) node {$\gpoint$};
    \draw (.5,4.15) node {$\gpoint$};
    \draw (.5,-4.15) node {$\gpoint$};
    \draw[blue,thick,dashed] (-.5,4.15)to(-.5,-4.15);
    \draw[blue,thick,dashed] (.5,4.15)to(.5,-4.15);
    \draw (-.5,4.7) node {\tiny$p$};
    \draw (-.5,-4.7) node {\tiny$q$};
    \draw (.5,4.7) node {\tiny$p'$};
    \draw (.5,-4.7) node {\tiny$q'$};
    \draw (0,0) node {\tiny$\zg$};

    \draw[bend left](-.5,4.15)to(-.5-4,3);
    \draw[bend left](-.5,4.15)to(-.5-2.8,.8);
    \draw[dotted](-.5-2.5,2.5)to(-.5-2,1.5);
    \draw (-.5-2.3,3.4) node {\tiny$\zb_1$};
    \draw (-.5-1.5,.8) node {\tiny$\zb_u$};

    \draw[bend right](-.5,-4.15)to(-.5-4,-3);
    \draw[bend right](-.5,-4.15)to(-.5-2.8,-.8);
    \draw[dotted](-.5-2.5,-2.5)to(-.5-2,-1.5);
    \draw (-.5-2.3,-3.4) node {\tiny$\zb_{v}$};
    \draw (-.5-1.5,-.8) node {\tiny$\zb_{u+1}$};

    \draw[bend right](.5,4.15)to(.5+4,3);
    \draw[bend right](.5+0,4.15)to(.5+2.8,.8);
    \draw[dotted](.5+2.5,2.5)to(.5+2,1.5);
    \draw (.5+2.3,3.4) node {\tiny$\za_{t}$};
    \draw (.5+1.5,.8) node {\tiny$\za_{s+1}$};

    \draw[bend left](.5+0,-4.15)to(.5+4,-3);
    \draw[bend left](.5+0,-4.15)to(.5+2.8,-.8);
    \draw[dotted](.5+2.5,-2.5)to(.5+2,-1.5);
    \draw (.5+2.3,-3.4) node {\tiny$\za_{1}$};
    \draw (.5+1.5,-.8) node {\tiny$\za_{s}$};

    \draw (6,0) node {};
\end{tikzpicture}
\begin{tikzpicture}[scale=0.4]

    \draw[bend right,thick](-3,5)to(3,5);
    \draw[bend left,thick](-3,-5)to(3,-5);

    \draw (0,4.15) node {$\gpoint$};
    \draw (0,-4.15) node {$\gpoint$};
    \draw[blue,thick] (0,4.15)to(0,-4.15);
    \draw (0,4.7) node {\tiny$p$};
    \draw (0,-4.7) node {\tiny$q$};
    \draw (.5,0) node {\tiny$\zg$};

    \draw[bend left](0,4.15)to(-4,3);
    \draw[bend left](0,4.15)to(-2.8,.8);
    \draw[dotted](-2.5,2.5)to(-2,1.5);
    \draw (-2.3,3.4) node {\tiny$\zb_1$};
    \draw (-1.5,.8) node {\tiny$\zb_u$};

    \draw[bend right](0,-4.15)to(-4,-3);
    \draw[bend right](0,-4.15)to(-2.8,-.8);
    \draw[dotted](-2.5,-2.5)to(-2,-1.5);
    \draw (-2.3,-3.4) node {\tiny$\zb_{v}$};
    \draw (-1.5,-.8) node {\tiny$\zb_{u+1}$};

    \draw[bend right](0,4.15)to(4,3);
    \draw[bend right](0,4.15)to(2.8,.8);
    \draw[dotted](2.5,2.5)to(2,1.5);
    \draw (2.3,3.4) node {\tiny$\za_{t}$};
    \draw (1.5,.8) node {\tiny$\za_{s+1}$};

    \draw[bend left](0,-4.15)to(4,-3);
    \draw[bend left](0,-4.15)to(2.8,-.8);
    \draw[dotted](2.5,-2.5)to(2,-1.5);
    \draw (2.3,-3.4) node {\tiny$\za_{1}$};
    \draw (1.5,-.8) node {\tiny$\za_{s}$};

\end{tikzpicture}
\end{center}
\begin{center}
\caption{Local configuration of lifting of arcs with endpoints $pq$ or $p'q'$ in an admissible collection $\zD_\zg$ of the cut surface. For different $i$ and $j$, $\overline{\zb}_i$ may coincide with $\overline{\zb}_j$.  Similar for the arcs $\overline{\za}_i$ and $\overline{\za}_j$. There may also exist arcs $\overline{\zb}_i$ and $\overline{\za}_j$ which coincide.}\label{figure:lift}
\end{center}
\end{figure}

\begin{lemma}\label{lemma:completion2}
Let $\zD$ be a collection of arcs on $(S,M)$ such that $\zD_\zg=\{\overline{\zg_i}, 1 \leq i \leq n, \zg_i\in \zD\}$
is an admissible collection on $(S_\zg,M_\zg)$. With the  above notations, the following statements hold.

(1) The collection $\zD$ is an admissible collection if and only if there are integers $0 \leq u \leq v$ and $0 \leq s \leq t$ such that the lifts of $pq$  and $p'q'$ respectively are such that, $p$ is an endpoint of $\zb_i$ for $1 \leq i \leq u$,  $p'$ is an endpoint of $\za_i$ for $1 \leq i \leq s$, $q$ is an endpoint of $\zb_j$ for $u+1 \leq j \leq v$ and $q'$ is an endpoint of $\za_j$ for $s+1 \leq j\leq t$, see the two rightmost pictures in Figure \ref{figure:lift}.  Furthermore, if $\zD_\zg$ is an admissible dissection, then $\zD$ is an admissible dissection if and only if the above conditions are satisfied.

(2) The collection $\zD$ is an exceptional collection (resp. dissection) if and only if the following conditions are satisfied:

(2.1) $\zD_\zg$ is an exceptional collection (resp. dissection);

(2.2) $\zD$ is an admissible collection consisting of exceptional arcs;

(2.3) for any two lifts $\zb_i$ and $\za_j$ with common endpoint $p$ (resp. $q$), either $\overline{\zb}_i$ and $\overline{\za}_j$ are not comparable or $\overline{\zb}_i \preceq \overline{\za}_j$ (resp. $\overline{\za}_j \preceq \overline{\zb}_i$).
\end{lemma}

We will denote  by $\zD_{(u,s)}$ the admissible dissection in $(S,M)$ as in Lemma~\ref{lemma:completion2} (1) determined by a pair of integers $(u,s)$.

\begin{proof}
(1) If there exist integers $0 \leq u \leq v$ and $0 \leq s \leq t$ satisfying the conditions in the first statement, then the lifting of arcs with endpoints $pq$ or $p'q'$ will not give rise to new interior intersections.
Thus $\zD$ is admissible since $\zD_\zg$ is. Conversely, assume that there exists no integer $0 \leq u \leq v$ or $0 \leq s \leq t$ satisfying the conditions in the first statement.
Without loss of generality assume that such an integer $u$ does not exist and that  $\zb_u$ has endpoint $q$ and $\zb_{u'}$ has endpoint $p$ with $u<u'$. Then $\zb_u$ and $\zb_{u'}$ intersect in the interior and $\zD$ is not admissible, see for example the right and middle pictures of Figure~\ref{figure:non-admisssible lift}.

(2) If $\zD$ is an exceptional collection, then in particular it is an admissible collection and any arc in $\zD$ is an exceptional arc with two distinct endpoints. By Lemma \ref{lemma:completione}, $\zD_\zg$ is an exceptional collection. So  conditions (2.1) and (2.2) hold. Furthermore, note that the condition $\overline{\za}_j \preceq \overline{\zb}_i$ implies that ${\za}_j \preceq {\zb}_i$ in $\zD$. On the other hand, the fact that $\zb_i$ and $\za_j$ share a common endpoint $p$ shows that $\zb_i \preceq \za_j$ in $\zD$. So $(\zb_i, \za_j)$ cannot be an exceptional pair, and $\zD$ will not be exceptional. A similar discussion works for the case that $\zb_i$ and $\za_j$ share a common endpoint $q$.

Conversely, assume that the three conditions in the second statement hold. We prove that $\zD$ is exceptional. Suppose for contradiction that there is a non-exceptional cycle in $\zD$, or equivalently, that there is an oriented cycle $\mathbf{c}$ in the quiver $Q$ associated to $\zD$. Since $Q$ is the quiver of a gentle algebra, there are at most four arrows at the vertex in $Q$ corresponding to  $\zg$, depicted as follows:

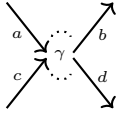
\begin{center}
{\begin{tikzpicture}[scale=0.7]
    \draw (0,0) node {\tiny$\zg$};
    \draw (-.8,.4) node {\tiny$a$};
    \draw (.8,.4) node {\tiny$b$};
    \draw (-.8,-.4) node {\tiny$c$};
    \draw (.8,-.4) node {\tiny$d$};
    \draw [thick,->] (-1,-1) -- (-0.25,-0.05);
    \draw [thick,<-] (1,-1) -- (0.25,-0.05);
    \draw [thick,->] (-1,1) -- (-0.25,0.05);
    \draw [thick,<-] (1,1) -- (0.25,0.05);

    \draw[dotted,thick](-.15,0.35) to [out=90,in=90] (.2,0.35);
    \draw[dotted,thick](-.15,-0.35) to [out=-90,in=-90] (.2,-0.35);
\end{tikzpicture}}
\end{center}
Note that since $\zD_\zg$ is exceptional, the associated quiver $\overline{Q}$ has no oriented cycles. On the other hand, $\overline{Q}$ is obtained from $Q$ by deleting the arrows attached to $\zg$ and then adding new arrows $[ab]$ and $[cd]$ (if $a, b$ and $c, d$, respectively, exist in $Q$). Thus the cycle $\mathbf{c}$ goes through $\zg$ and contains a deleted arrow in $Q$.

Assume $\mathbf{c}$ contains $a$.  Then $\mathbf{c}$ also contains $b$ or $d$, since it is a cycle. In the following we always assume that  $\zg$ appears in $\mathbf{c}$ once, and each deleted arrow appears at most once in $\mathbf{c}$ and we call $\mathbf{c}$ a minimal cycle.

Case 1. If $b$ belongs to $\mathbf{c}$, since $\mathbf{c}$ is minimal, $a$ and $b$ appear in $\mathbf{c}$ just once and no other deleted arrow appears in $\mathbf{c}$. After replacing $ab$ in $\mathbf{c}$ by $[ab]$, we obtain a path in $\overline{Q}$, which is still an oriented cycle. This contradicts  the assumption that $\zD_\zg$ is exceptional.

Case 2. If $d$ belongs to $\mathbf{c}$, then we have two subcases. If $ad=\mathbf{c}$, then the source of $a$ and the target of $d$ coincide and we denote this vertex by $\za$. Then the corresponding arc $\za$ on $(S,M)$ is a loop in $\zD$, which contradicts  the assumption that any arc in $\zD$ is exceptional.
Now assume $ad \neq \mathbf{c}$. We denote by $\zb$ the source of $a$ and by $\za$ the target of $d$. Since the composition of $a$ and $d$ is nonzero, $\zb$ and $\za$ share an endpoint $p$ or $q$. Without loss of generality, we assume that they share the endpoint $p$.  On the other hand, under the assumption that $\mathbf{c}$ is minimal, the path obtained from $\mathbf{c}$ by deleting the arrows $a$ and $d$ is a path in $\overline{Q}$, which is still connected and gives rise to the relation $\overline{\za} \preceq \overline{\zb}$. This contradicts condition $(2.3)$ in the second statement.

The proof for the case that $\mathbf{c}$ contains $c$ is similar.
To sum up, if the three conditions in the second statement hold, then $\zD$ is exceptional.
 \end{proof}

The following corollary is a direct consequence of Lemma~\ref{lemma:completion2}.
\begin{corollary}\label{corollary:completion2}
For any exceptional collection $\zD_\zg$ on $(S_\zg,M_\zg)$, the admissible collection $\zD_{(v,t)}$ is the unique exceptional collection on $(S,M)$ such that $\zD_\zg=\overline{\zD}_{(v,t)}$ and such that $\zg$ is a maximal element in $\zD_{(v,t)}$. Furthermore, if $\zD_\zg$ is an exceptional dissection, then $\zD_{(v,t)}$ is also an exceptional dissection.
\end{corollary}

We apply iterated cuts of surfaces to prove the following.

\begin{proposition}\label{prop:quivers of two surfaces}
The quiver with relations associated to any exceptional dissection on the surface $\mathbb{T}_{(g,1,2)}$ is of the form

\begin{equation}\label{eq:quiver1}
{\begin{tikzpicture}[scale=0.4]
\def \radius {4cm}
    \draw [thick,->] (0,.2) -- (2,.2);
    \draw [thick,->] (.5+3,.2) -- (2.5+3,.2);
    \draw [thick,->] (.5+9,.2) -- (2.5+9,.2);
    \draw [thick,->] (.5+12,.2) -- (2.5+12,.2);
    \draw [thick,->] (.5+15,.2) -- (2.5+15,.2);

    \draw [thick,->] (0,-.2) -- (2,-.2);
    \draw [thick,->] (.5+3,-.2) -- (2.5+3,-.2);
    \draw [thick,->] (.5+9,-.2) -- (2.5+9,-.2);
    \draw [thick,->] (.5+12,-.2) -- (2.5+12,-.2);
    \draw [thick,->] (.5+15,-.2) -- (2.5+15,-.2);

    \draw [thick,dotted] (.5+6,0) -- (2.5+6,0);

    \draw[dotted,thick](1.5,.5) to [out=90,in=90] (3.5,.5);
    \draw[dotted,thick](1.5,-.5) to [out=-90,in=-90] (3.5,-.5);
    \draw[dotted,thick](1.5+3.4,.5) to [out=90,in=90] (3.5+3.4,.5);
    \draw[dotted,thick](1.5+3.4,-.5) to [out=-90,in=-90] (3.5+3.4,-.5);
    \draw[dotted,thick](1.5+3.4+3,.5) to [out=90,in=90] (3.5+3.4+3,.5);
    \draw[dotted,thick](1.5+3.4+3,-.5) to [out=-90,in=-90] (3.5+3.4+3,-.5);
    \draw[dotted,thick](1.5+3.4+3+3,.5) to [out=90,in=90] (3.5+3.4+3+3,.5);
    \draw[dotted,thick](1.5+3.4+3+3,-.5) to [out=-90,in=-90] (3.5+3.4+3+3,-.5);
    \draw[dotted,thick](1.5+3.4+3+3+3,.5) to [out=90,in=90] (3.5+3.4+3+3+3,.5);
    \draw[dotted,thick](1.5+3.4+3+3+3,-.5) to [out=-90,in=-90] (3.5+3.4+3+3+3,-.5);

    \draw (-1.2,0) node {\tiny$2g+1$};
    \draw (2.7,0) node {\tiny$2g$};
    \draw (12,0) node {\tiny$3$};
    \draw (15,0) node {\tiny$2$};
    \draw (18,0) node {\tiny$1$};
\end{tikzpicture}}
\end{equation}

and  on the surface $\mathbb{T}_{(g,2,2)}$ it is of the form

\begin{equation}\label{eq:quiver2}
{\begin{tikzpicture}[scale=0.4]
\def \radius {4cm}
    \draw [thick,->] (-1,.2) -- (1,.2);
    \draw [thick,->] (.5+3,.2) -- (2.5+3,.2);
    \draw [thick,->] (.5+9,.2) -- (2.5+9,.2);
    \draw [thick,->] (.5+12,.2) -- (2.5+12,.2);
    \draw [thick,->] (.5+15,.2) -- (2.5+15,.2);

    \draw [thick,->] (-1,-.2) -- (1,-.2);
    \draw [thick,->] (.5+3,-.2) -- (2.5+3,-.2);
    \draw [thick,->] (.5+9,-.2) -- (2.5+9,-.2);
    \draw [thick,->] (.5+12,-.2) -- (2.5+12,-.2);
    \draw [thick,->] (.5+15,-.2) -- (2.5+15,-.2);

    \draw [thick,dotted] (.5+6,0) -- (2.5+6,0);

    \draw[dotted,thick](.6,.5) to [out=90,in=90] (3.5,.5);
    \draw[dotted,thick](.6,-.5) to [out=-90,in=-90] (3.5,-.5);
    \draw[dotted,thick](1.5+3.4,.5) to [out=90,in=90] (3.5+3.4,.5);
    \draw[dotted,thick](1.5+3.4,-.5) to [out=-90,in=-90] (3.5+3.4,-.5);
    \draw[dotted,thick](1.5+3.4+3,.5) to [out=90,in=90] (3.5+3.4+3,.5);
    \draw[dotted,thick](1.5+3.4+3,-.5) to [out=-90,in=-90] (3.5+3.4+3,-.5);
    \draw[dotted,thick](1.5+3.4+3+3,.5) to [out=90,in=90] (3.5+3.4+3+3,.5);
    \draw[dotted,thick](1.5+3.4+3+3,-.5) to [out=-90,in=-90] (3.5+3.4+3+3,-.5);
    \draw[dotted,thick](1.5+3.4+3+3+3,.5) to [out=90,in=90] (3.5+3.4+3+3+3,.5);
    \draw[dotted,thick](1.5+3.4+3+3+3,-.5) to [out=-90,in=-90] (3.5+3.4+3+3+3,-.5);

    \draw (-2.5,0) node {\tiny$2g+2$};
    \draw (2.3,0) node {\tiny$2g+1$};
    \draw (12,0) node {\tiny$3$};
    \draw (15,0) node {\tiny$2$};
    \draw (18,0) node {\tiny$1.$};
\end{tikzpicture}}
\end{equation}
\end{proposition}

\begin{proof}
At first we show that any quiver $Q$ associated to an exceptional dissection on the surface $\mathbb{T}_{(1,1,2)}$ is isomorphic to a quiver as in \eqref{eq:quiver1}.  
By Lemma \ref{lemma:surface and quiver} there are three vertices and four arrows in $Q$ and  there must exist two vertices, for example vertices  $3$ and $2$, such that there are two arrows between them. On the other hand, since there is no oriented cycle in $Q$, these two arrows have the same direction. So we can assume that there are two arrows from $3$ to $2$ in $Q$. Furthermore, since $Q$ is the quiver of a gentle algebra, there are at most two arrows starting and at most to arrow ending at a given vertex,  so the other two arrows must either start at $2$ and end at $1$, or start at $1$ and end at $3$. So $Q$ is  isomorphic to a quiver as \eqref{eq:quiver1}.
Furthermore, for such a quiver, up to isomorphism there is a unique choice of relations, depicted in \eqref{eq:quiver1}, such that the resulting quotient algebra is gentle.

Let $\zD$ be an exceptional dissection on the surface $(S,M)=\mathbb{T}_{(1,2,2)}$. Let $\zg$ be any arc in $\zD$, then the cut surface $(S_\zg,M_\zg)$ is homeomorphic to the surface $\mathbb{T}_{(1,1,2)}$ with an exceptional dissection $\overline{\zD}$, whose quiver $\overline{Q}$ is then necessarily  as in \eqref{eq:quiver1}. Furthermore, we may assume that $\zg$ is maximal in $\zD$ with respect to the partial order of $\zD$. Then by Lemma \ref{lemma:completion}, the quiver $Q$ of $\zD$ is obtained from $\overline{Q}$ by adding a vertex corresponding to $\zg$ and two new arrows. Note that the latter follows from the fact that by Lemma \ref{lemma:surface and quiver}, there are 6 arrows in $Q$. By the maximality of $\zg$ the two new arrows start at $\zg$. Since $\zD$ is an exceptional dissection, there are no oriented cycles in $Q$ and thus in order for $Q$ to give rise to a gentle algebra, it must be  as in \eqref{eq:quiver2}.
Thus the statement holds for surfaces of genus one.

For the general case, the following two claims guarantee that we can use inductions on the genus of the surface, where the first claim is clear and we will prove the second one.

Claim 1. Let $\zg$ be any arc in an exceptional dissection $\zD$ of the surface $(S,M)=\mathbb{T}_{(g,2,2)}$. Then the cut surface $(S_\zg,M_\zg)$ is homeomorphic to  $\mathbb{T}_{(g,1,2)}$.

Claim 2. Let $\zg$ be a maximal arc in an exceptional dissection $\zD$ of the surface $(S,M)=\mathbb{T}_{(g,1,2)}$. Then the cut surface $(S_\zg,M_\zg)$ is homeomorphic to $\mathbb{T}_{(g-1,2,2)}$.

We now show Claim 2. Since $\zg$ is a  minimal arc in $\zD$, by Lemma \ref{lemma:completion}, the quiver $\overline{Q}$ is obtained from $Q$ by deleting the vertex corresponding to $\zg$ (which is a  source) and the arrows incident with it. Note that since the surface $(S,M)$ is connected, the quiver $Q$ is connected. So the quiver $\overline{Q}$ is connected, and thus the cut surface $(S_\zg,M_\zg)$ is also connected. Therefore $\zg$ is non-separating. Denote by $\tilde{S}$ the compact surface obtained from $S$ by contracting the boundary to a point. Then $\zg$ becomes a simple closed curve on $\tilde{S}$, which is denoted by $\tilde{\zg}$. By Lemma \ref{lemma:simple closed curves}, there is a homeomorphism $\tilde{\phi}:\tilde{S}\mapsto \tilde{S}$ which maps $\tilde{\zg}$ to a simple closed curve $\tilde{\phi}(\tilde{\zg})$ which is a curve going around a handle of $\tilde{S}$ once. In fact, $\tilde{\phi}$ can be lifted to a homeomorphism $\phi$ of $(S,M)$ by keeping the local configuration of the boundary, where $\phi(\zg)$  is a curve going around a handle of $(S,M)$ once.  Thus the cut surface of $(S,M)$ along $\phi(\zg)$ is homeomorphic to $\mathbb{T}_{(g-1,2,2)}$. Furthermore, the cut surface of $(S,M)$ along $\zg$ is also homeomorphic to $\mathbb{T}_{(g-1,2,2)}$, since $\phi$ is a homeomorphism of $(S,M)$.
\end{proof}

\begin{corollary}\label{cor:quivers of two surfaces}
(1) Any arc in an exceptional dissection on  the surfaces $\mathbb{T}_{(g,1,2)}$ and on $\mathbb{T}_{(g,2,2)}$ is a non-separating arc.

(2) Let $\zD_1$ and $\zD_2$ be two exceptional dissections on $\mathbb{T}_{(g,1,2)}$ (on $\mathbb{T}_{(g,2,2)}$ resp.), then there is homeomorphism $\phi$ of $\mathbb{T}_{(g,1,2)}$ (of $\mathbb{T}_{(g,2,2)}$ resp.) such that $\phi(\zD_1)=\zD_2$.
\end{corollary}

\begin{proof}
(1) Let $\zg$ be an arc in an exceptional dissection $\zD$ on the surface $\mathbb{T}_{(g,1,2)}$ or on $\mathbb{T}_{(g,2,2)}$.
By Proposition \ref{prop:quivers of two surfaces}, the quiver with relations associated to $\zD$ is as in \eqref{eq:quiver1} or \eqref{eq:quiver2} respectively. On the other hand, by the construction  in Lemma~\ref{lemma:completion} of the quiver $\overline{Q}$ associated to the exceptional dissection $\overline{\zD}$, $\overline{Q}$ is still connected. Thus $\zg$ is non-separating.

 (2) This follows from Proposition \ref{prop:quivers of two surfaces} and the bijective correspondence between gentle quivers and admissible dissections of the surfaces.
\end{proof}

As another application of cutting marked surfaces, we will describe when an exceptional sequence in $\ka$ can be completed to a full exceptional sequence. We begin with the following straightforward observation.

\begin{lemma}\label{lemma:completion1}
An exceptional sequence $(\P_{(\zg_1,f_{1})},\cdots,\P_{(\zg_m,f_{m})})$ in $\ka$ can be completed to a full exceptional sequence if and only if the associated exceptional collection $\{\zg_1,\cdots,\zg_m\}$ on $(S,M)$ can be completed to an exceptional dissection.
\end{lemma}

We now give an example showing that not every exceptional sequence can be completed to a full exceptional sequence.

\begin{example}\label{example:completion1}
Consider the marked surface $\mathbb{T}_{(g,1,2)}$ of genus $g \geq 1$ with one boundary component and  two $\gpoint$-marked point $p$ and $q$. Let $\zg$ be an arc connecting $p$ and $q$, and such that $\zg$ is homotopic to a boundary segment, see Figure \ref{figure:completion1}. Then for any grading $f_\zg$ on $\zg$, the associated object $\P_{(\zg,f_\zg)}$ is exceptional. However, it cannot be completed to a full exceptional sequence: Let $\za$ be an $\gpoint$-arc which has no self-intersections. So in particular, $\za$ must be have  endpoints $p$ and $q$. If $\za$ is not homotopic to $\zg$, then $\zg$ and $\za$ form a non-exceptional cycle, and $\{\zg,\za\}$ cannot be an exceptional collection. Thus by Lemma \ref{lemma:completion1}, $\P_{(\zg,f_\zg)}$ can not be completed to a full exceptional sequence.
\begin{figure}[ht]
\begin{center}
\begin{tikzpicture}[scale=0.3]
\begin{scope}
	\draw (0,0) circle (3cm);
	\clip[draw] (0,0) circle (3cm);
	\foreach \x in {-4,-3.5,-3,-2.5,-2,-1.5,-1,-0.5,0,0.5,1,1.5,2,2.5,3,3.5,4}	\draw[xshift=\x cm]  (-5,5)--(5,-5);
\end{scope}
    \draw (0,3) node {$\gpoint$};
    \draw (3,0) node {$\gpoint$};
    \draw (2.1,2.1) node {$\rpoint$};
    \draw (-2.1,-2.1) node {$\rpoint$};
    \draw (0,3.7) node {$p$};
        \draw (3.6,0) node {$q$};
        \draw (3.5,3.5) node {$\zg$};
    \draw plot [smooth,tension=1] coordinates {(0,3) (3,3) (3,0)};

\end{tikzpicture}
\end{center}
\begin{center}
\caption{The local configuration of an arc $\zg$ on $\mathbb{T}_{(1,1,2)}$.}\label{figure:completion1}
\end{center}
\end{figure}
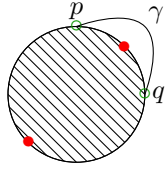
\end{example}

Let $\{\zg_1,\zg_2\}$ be an exceptional collection, it is easy to see that
$$((S_{\zg_1})_{\zg_2},(M_{\zg_1})_{\zg_2})
=((S_{\zg_2})_{\zg_1},(M_{\zg_2})_{\zg_1}).$$ This ensures that for any exceptional collection $\zD$, we may define the cut surface $(S_\zD,M_\zD)$ by successively cutting the arcs in $\zD$ in any order.

\begin{proposition}\label{proposition:completion1}
Let $\zD=\{\zg_1,\cdots,\zg_m\}$ be an exceptional collection on $(S,M)$. Then $\zD$ can be completed to an exceptional dissection if and only if $(S_\zD,M_\zD)$ has no connected component which is homeomorphic to $\mathbb{T}_{(g,1,1)}$.
\end{proposition}
\begin{proof}
We prove this by induction on $m$. We begin with the case $m=1$. In this case $\zD=\{\zg_1\}$.
 Suppose that there is no connected component of  $(S_\zD,M_\zD)$ homeomorphic to $\mathbb{T}_{(g,1,1)}$. Then by Proposition \ref{thm:existence of full exceptional sequence}, there exists an exceptional dissection $\zD_{\zg_1}$ on $(S_\zD,M_\zD)$. Furthermore, by Corollary \ref{corollary:completion2}, there exists an exceptional dissection $\zD$ of $(S,M)$ containing $\zg_1$, that is, $\zD$ can be completed to an exceptional dissection on $(S,M)$.

For the general case, let $\zD=\{\zg_1,\cdots,\zg_m\}$ be an exceptional collection on $(S,M)$ and suppose that $(S_\zD,M_\zD)$ has no connected component homeomorphic to $\mathbb{T}_{(g,1,1)}$. Then $\{\zg_2,\cdots,\zg_m\}$ is an exceptional collection on $(S_{\zg_1},M_{\zg_1})$ by Lemma \ref{lemma:completione}. By induction, it can be completed to an exceptional dissection on $(S_{\zg_1},M_{\zg_1})$. Then again by Corollary \ref{corollary:completion2}, there exists an exceptional dissection on $(S,M)$ containing $\{\zg_1,\cdots,\zg_m\}$.

The converse directly follows from the fact that if there is a connected component of the form $\mathbb{T}_{(g,1,1)}$, then there are no exceptional arcs on $\mathbb{T}_{(g,1,1)}$ since this  connected component has a single marked point on its unique boundary.
\end{proof}

Now we have the following geometric characterisation of when an exceptional sequence can be completed to a full exceptional sequence.

\begin{theorem}\label{theorem:completion1}
Let $A$ be a gentle algebra associated with a marked surface $(S,M,\zD_A)$ with associated admissible dissection $\zD_A$. Assume there exist full exceptional sequences in $\ka$ and let   $(\P_{(\zg_1,f_{1})},\cdots,\P_{(\zg_m,f_{m})})$
 be an exceptional sequence arising from an exceptional collection $\zD=\{\zg_1,\cdots,\zg_m\}$. Then it can be completed to a full exceptional sequence if and only if the cut surface $(S_\zD,M_\zD)$ has no subsurface of the form $\mathbb{T}_{(g,1,1)}$.
\end{theorem}

\begin{proof}
This follows directly from Lemma \ref{lemma:completion1} and Proposition \ref{proposition:completion1}.
\end{proof}

\section{Transitivity of the braid group action}\label{subsection:Braid group action on exceptional dissections}

 In this section we will use the geometric model to show that the induced action of  the Artin braid group is transitive on the set of full exceptional sequences in the perfect derived category of a gentle algebra corresponding to a genus zero surface. More precisely, denote by $B_n$ the Artin braid group generated by $\sigma_1,\cdots,\sigma_{n-1}$
with relations
$\sigma_i\sigma_j=\sigma_j\sigma_i, |i-j|>1$ and
$\sigma_i\sigma_{i+1}\sigma_i=\sigma_{i+1}\sigma_i\sigma_{i+1}$. Then we show the following.

\begin{theorem}\label{theorem:transitivity2}
Let $\A$ be a gentle algebra arising from a marked surface of genus zero, let $n$ be the rank of $K_0(\ka)$ and $B_n$ the Artin braid group on $n$ strands. Then  the action of  $\mathbb{Z}^n\ltimes B_n$ on the set of full exceptional sequences in $\ka$ is transitive.
\end{theorem}

The proof is based on  the following strategy:
In the first half of this section we introduce an action of $B_n$ on the set of ordered exceptional dissections of a marked surface $(S,M)$. We then discuss the transitivity of this action and show that it holds for surfaces satisfying a certain condition ({\bf RCEA} condition). In the second half of this section we prove the theorem by showing how this geometric interpretation translates to the braid group action on full exceptional sequences in the associated derived category $\ka$.

\subsection{Braid group action on exceptional surface dissections}

We start with the following observation.

\begin{lemma}
Let $(\zg_1,\zg_2)$ be an ordered exceptional pair in $(S,M)$. Then $\zg_1$ and $\zg_2$ either do not intersect or they intersect at one or two common endpoints, as shown  in Figure~\ref{figure:mutation of exceptional dissection}.
\end{lemma}

\begin{definition}\label{definition:mutations}
(1) Let $(\zg_1,\zg_2)$ be an ordered exceptional pair in $(S,M)$. We now define \emph{the left (resp. right) mutation of $\zg_2$ at $\zg_1$} denoted by $L_{\zg_1}\zg_2$ (resp. $R_{\zg_1}\zg_2$) as follows:
\begin{itemize}
\item If $\zg_1$ and $\zg_2$ do not intersect, then we define $L_{\zg_1}\zg_2$ as $\zg_2$ and define $R_{\zg_2}\zg_1$ as $\zg_1$.
\item  If $\zg_1$ and $\zg_2$ share one endpoint $q_1$ then $L_{\zg_1}\zg_2$ is the smoothing of the crossing of $\zg_1$ and $\zg_2$  at $q_1$.
\item If $\zg_1$ and $\zg_2$ share both endpoints, then   $L_{\zg_1}\zg_2$ is obtained by first smoothing the crossing of  $\zg_1$ and $\zg_2$ at one of the endpoints and then smoothing the crossing of the resulting arc with $\zg_1$ at the other endpoint of $\zg_1$, see  Figure \ref{figure:mutation of exceptional dissection}.
\end{itemize}
Dually, we define $R_{\zg_2}\zg_1$ when $\zg_1$ and $\zg_2$ intersect in one or both endpoints, see  Figure \ref{figure:mutation of exceptional dissection}.

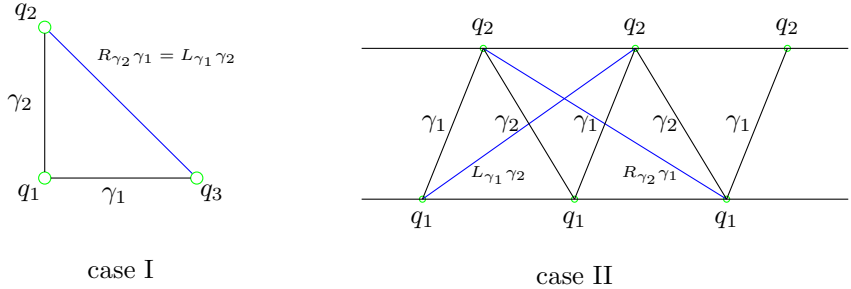
\begin{figure}[H]
\begin{center}
{\begin{tikzpicture}[scale=0.4]
\draw[green] (0,0) circle [radius=0.2];
\draw[green] (0,5) circle [radius=0.2];
\draw[green] (5,0) circle [radius=0.2];

\draw[-] (0,0.2) -- (0,4.8);

\draw[-] (0.2,0) -- (4.8,0);

\draw[-,blue] (0.15,4.85) -- (4.85,0.15);

\draw (2.3,-.6) node {$\zg_1$};

\draw (-.8,2.5) node {$\zg_{2}$};
\draw (4,4) node {\tiny$R_{\zg_2}\zg_1=L_{\zg_1}\zg_2$};
\draw (-.6,5.5) node {$q_2$};
\draw (-.5,-.5) node {$q_1$};
\draw (5.5,-.5) node {$q_3$};
%

\draw (2.5,-3) node {case I};
\draw (2.5,-3.5) node {};

\end{tikzpicture}}
\qquad\qquad
{\begin{tikzpicture}[scale=0.4]
\draw[green] (0,-2.5) circle [radius=0.1];
\draw[green] (5,-2.5) circle [radius=0.1];
\draw[green] (-5,-2.5) circle [radius=0.1];

\draw[green] (7,2.5) circle [radius=0.1];
\draw[green] (2,2.5) circle [radius=0.1];
\draw[green] (-3,2.5) circle [radius=0.1];

\draw[-,blue]  (-3,2.5)to(5,-2.5);
\draw[-,blue]  (2,2.5)to(-5,-2.5);

\draw[-]  (-7,-2.5)to(9,-2.5);
\draw[-]  (-7,2.5)to(9,2.5);

\draw (-4.6,0) node {$\zg_1$};
\draw (0.4,0) node {$\zg_1$};
\draw (5.4,0) node {$\zg_1$};

\draw (-2.2,0) node {$\zg_2$};
\draw (3,0) node {$\zg_2$};

\draw (-5,-3.2) node {$q_1$};
\draw (0,-3.2) node {$q_1$};
\draw (5,-3.2) node {$q_1$};

\draw[-]  (-5,-2.5)to(-3,2.5);
\draw[-]  (0,-2.5)to(-3,2.5);
\draw[-]  (0,-2.5)to(2,2.5);
\draw[-]  (5,-2.5)to(2,2.5);
\draw[-]  (5,-2.5)to(7,2.5);

\draw (-3,3.2) node {$q_2$};
\draw (2,3.2) node {$q_2$};
\draw (7,3.2) node {$q_2$};

\draw (-2.5,-1.7) node {\tiny$L_{\zg_1}\zg_2$};
\draw (2.5,-1.7) node {\tiny$R_{\zg_2}\zg_1$};

\draw (0,-5) node {case II};
\draw (2.5,-3.5) node {};

\end{tikzpicture}}
\end{center}

\begin{center}
\caption{Possible intersections of  $\zg_{1}$ and $\zg_{2}$ in an ordered exceptional pair $(\zg_1,\zg_2)$, where case II is depicted in the universal covering. This figure also illustrates the mutations of the ordered exceptional pair $(\zg_1,\zg_2)$ (in blue).}\label{figure:mutation of exceptional dissection}
\end{center}
\end{figure}
(2) Let $\zG=(\zg_1,\cdots,\zg_m)$ be an ordered exceptional collection of arcs in $(S,M)$, for any $1 \leq i \leq m-1$, we define
\vspace{-.2cm}
$$\begin{array} {l}
\sigma_i \zG=(\zg_1,\cdots,\zg_{i-1},\zg_{i+1}
,R_{\zg_{i+1}}\zg_i,\zg_{i+2},\cdots,\zg_m),
\\
\\
\sigma_i^{-1} \zG=(\zg_1,\cdots,\zg_{i-1},L_{\zg_{i}}\zg_{i+1},\zg_{i}
,\zg_{i+2},\cdots,\zg_m).\\
\end{array}$$
\end{definition}

\medskip

\begin{theorem}\label{theorem:mutation}
Let $\zG=(\zg_1,\cdots,\zg_m)$ be an ordered exceptional collection on $(S,M)$. Then  $\sigma_i \zG$ and $\sigma^{-1}_i \zG$ are  ordered exceptional collections on $(S,M)$ of length $m$, for all $1 \leq i \leq m-1$.
Furthermore, we have
 $$\begin{array}{l} \sigma^{-1}_i\sigma_i =\sigma_i\sigma^{-1}_i =id, \\
\sigma_{i}\sigma_{i+1}\sigma_{i}=\sigma_{i+1}\sigma_{i}\sigma_{i+1},\\
\sigma_{i}\sigma_{j}=\sigma_{j}\sigma_{i}, ~\rm{if}~|i-j| > 1,\
\end{array}$$
and this induces an action of the  braid group $B_m$ on the set of ordered exceptional collections on $(S,M)$ of length $m$.
\end{theorem}
\begin{proof}
The statement that $\sigma_i \zG$ and $\sigma^{-1}_i \zG$ are ordered exceptional collections on $(S,M)$ of length $m$ can be checked directly using the definitions. It is also straightforward to check the equalities  $\sigma^{-1}_i\sigma_i =\sigma_i\sigma^{-1}_i =id$ and $\sigma_{i}\sigma_{j}=\sigma_{j}\sigma_{i}, ~\rm{if}~|i-j| > 1$. Now we consider the braid relations.
Without loss of generality, we may assume $m=3$. Then
\[\sigma_1\sigma_2\sigma_1(\zg_1,\zg_2,\zg_3)=
(\zg_3,R_{\zg_3}\zg_2,R_{\zg_3}R_{\zg_2}\zg_1),\]
\[\sigma_2\sigma_1\sigma_2(\zg_1,\zg_2,\zg_3)=
(\zg_3,R_{\zg_3}\zg_2,R_{R_{\zg_3}\zg_2}R_{\zg_3}\zg_1).\]
So we have to show that $R_{\zg_3}R_{\zg_2}\zg_1=R_{R_{\zg_3}\zg_2}R_{\zg_3}\zg_1.$

If there exists an arc $\zg_j, 1\leq j\leq 3$ which shares no endpoints with the other two arcs, then it can be easily checked that $R_{\zg_3}R_{\zg_2}\zg_1=R_{R_{\zg_3}\zg_2}R_{\zg_3}\zg_1=\zg_1$ if $j=1$, $R_{\zg_3}R_{\zg_2}\zg_1=R_{R_{\zg_3}\zg_2}R_{\zg_3}\zg_1=R_{\zg_3}\zg_1$ if $j=2$, $R_{\zg_3}R_{\zg_2}\zg_1=R_{R_{\zg_3}\zg_2}R_{\zg_3}\zg_1=R_{\zg_2}\zg_1$
if $j=3$.

If any two arcs intersect at one or both of the endpoints then there are several cases to consider depending on the number of the intersections as well as the positions of the different arcs. The result then follows from a straightforward verification of the definitions. We prove  one of the cases in the following figure, the other cases being similar.


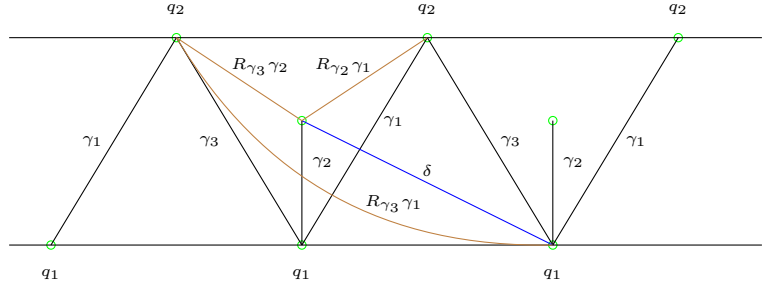
\begin{figure}[H]
\begin{center}
{\begin{tikzpicture}[scale=0.55]
\draw[green] (0,-2.5) circle [radius=0.1];
\draw[green] (6,-2.5) circle [radius=0.1];
\draw[green] (-6,-2.5) circle [radius=0.1];

\draw[green] (9,2.5) circle [radius=0.1];
\draw[green] (3,2.5) circle [radius=0.1];
\draw[green] (-3,2.5) circle [radius=0.1];

\draw[green] (0,0.5) circle [radius=0.1];
\draw[green] (6,0.5) circle [radius=0.1];

\draw[-]  (-7,-2.5)to(11,-2.5);
\draw[-]  (-7,2.5)to(11,2.5);

\draw (-5,0) node {\tiny$\zg_1$};
\draw (2.2,.5) node {\tiny$\zg_1$};
\draw (8,0) node {\tiny$\zg_1$};

\draw (-2.2,0) node {\tiny$\zg_3$};
\draw (5,0) node {\tiny$\zg_3$};

\draw (6.5,-.5) node {\tiny$\zg_2$};
\draw (.5,-.5) node {\tiny$\zg_2$};
\draw (3,-.7) node {\tiny$\delta$};

\draw (2.2,-1.5) node {\tiny$R_{\zg_3}\zg_1$};
\draw (1,1.8) node {\tiny$R_{\zg_2}\zg_1$};
\draw (-1,1.8) node {\tiny$R_{\zg_3}\zg_2$};


\draw (-6,-3.2) node {\tiny$q_1$};
\draw (0,-3.2) node {\tiny$q_1$};
\draw (6,-3.2) node {\tiny$q_1$};

\draw[-]  (-6,-2.5)to(-3,2.5);
\draw[-]  (0,-2.5)to(-3,2.5);
\draw[-]  (0,-2.5)to(3,2.5);
\draw[-]  (6,-2.5)to(3,2.5);
\draw[-]  (6,-2.5)to(9,2.5);

\draw[-,blue]  (0,0.5)to(6,-2.5);
\draw[-]  (0,0.5)to(0,-2.5);
\draw[-]  (6,0.5)to(6,-2.5);
\draw[-,brown]  (0,0.5)to(3,2.5);
\draw[-,brown]  (0,0.5)to(-3,2.5);
\draw[-,bend left,brown]  (6,-2.5)to(-3,2.5);

\draw (-3,3.2) node {\tiny$q_2$};
\draw (3,3.2) node {\tiny$q_2$};
\draw (9,3.2) node {\tiny$q_2$};
%
\end{tikzpicture}}
\end{center}

\begin{center}
\label{figure:braid group action of exceptional dissections3}
\end{center}
\end{figure}
In this situation it directly follows from  the figure that  $R_{\zg_3}R_{\zg_2}\zg_1=R_{R_{\zg_3}\zg_2}R_{\zg_3}\zg_1$, which corresponds to the blue arc $\delta$ for the case when $\zg_2$ shares one endpoint both with $\zg_1$ and $\zg_3$, and $\zg_1$, $\zg_3$ share two endpoints.
\end{proof}

In the following we consider the transitivity of the braid group action on the set of ordered exceptional dissections on $(S,M)$.
For this we first show some technical lemmas.

We observe that different choices of order for a given exceptional collection on $(S,M)$ gives rise to  distinct ordered exceptional collections. However, we show that they are all related through the braid group action.

\begin{lemma}\label{lemma:braid mutate same collection}
Let $\zG=(\zg_1,\cdots,\zg_m)$ and $\zG'=(\zg'_1,\cdots,\zg'_m)$ be two ordered exceptional collections on $(S,M)$ arising from the same exceptional collection $\zD$. Then there is an element $\sigma$ in $B_m$ such that $\sigma \zG=\zG'$.
\end{lemma}
\begin{proof}
Suppose that $\zg_i=\zg'_m$ for some $1\leq i\leq m$. Then note that $\zg_i$ is a maximal element with respect to the partial order on $\zD$. Thus there are no intersections between $\zg_i$ and $\zg_j$, for  $i+1 \leq j \leq m$ and $R_{\zg_j}\zg_i=\zg_i$. So $\sigma_{m-1}\cdots\sigma_i \zG=(\zg_1,\cdots,\zg_{i-1},\zg_{i+1},\cdots,\zg_{m},\zg_i=\zg'_m)$ which is obtained from $\zG$ by moving $\zg'_m$ to the last position and shifting the arcs $\zg_{i+1}, \cdots, \zg_m$ one position to the left.  Then we consider the position of $\zg'_{m-1}$ in $\sigma_{m-1}\cdots\sigma_i \zG$, and similarly mutate it to the position before $\zg'_m$. Continuing in this way, we can  mutate $\zG$ to $\zG'$, that is, we have inductively constructed an element $\sigma$ such that $\sigma\zG=\zG'$.
\end{proof}

For the rest of this section we fix the following conventions.
Let $(S,M)$ be a marked surface with exceptional dissections.
Let $\zg$ be an exceptional $\gpoint$-arc such that there exist exceptional dissections on the cut surface $(S_\zg,M_\zg)$.
We extend the definition of the braid group action to a non-connected marked surface. Note that if the braid group $B_m$ acts on the set of ordered exceptional dissections of a marked surface $(S,M)$, then the braid group $B_{m-1}$ acts on the cut surface $(S_\zg,M_\zg)$, independently from whether  $(S_\zg,M_\zg)$ is connected or not.

Let $\zD_\zg$ be an admissible dissection on $(S_\zg,M_\zg)$.
Recall from Lemma \ref{lemma:completion2} that $\zD_{(u,s)}$ is the unique lift of $\zD_\zg$ to $(S,M)$ which is an admissible dissection determined by the pair $(u,s)$ such that $0 \leq u \leq v, 0 \leq s \leq t$, where $v$ and $t$ are as in the paragraph preceding Lemma~\ref{lemma:completion2}. For two pairs $(u,s)$ and $(u',s')$, we say $(u,s)<(u',s')$ if $u \leq u', s \leq s'$ and at least one of the inequalities is strict. The following lemma allows us to later use induction.

\begin{lemma}\label{lemma:braid mutate same underlying set}
Let $\zD_{(u,s)}$ be a lift of $\zD_\zg$ which is exceptional.
If $(u,s)<(v,t)$, then there exists an ordered exceptional dissection $\zG_{(u,s)}$ arising from $\zD_{(u,s)}$ and an element $\sigma$ in $B_n$ such that $\sigma\zG_{(u,s)}$ is an ordered exceptional dissection with underlying set $\zD_{(u',s')}$ for some $(u',s')$, which is a lift of $\zD_\zg$ with $(u,s)<(u',s')$.
\end{lemma}
\begin{proof}
We first note the following  useful fact: Using the notations in Figure~\ref{figure:lift}, $\za_{s+1}$ and $\zb_{u+1}$ are the only two minimal elements in $\zD_{(u,s)}$ which are greater than $\zg$, that is, if $\zg \preceq \za$ for some $\za\in \zD_{(u,s)}$, then at least one of $\za_{s+1} \preceq \za$ or $\zb_{u+1} \preceq \za$ holds.

Case 1. $\za_{s+1}=\zb_{u+1}$. Let $\zG_{(u,s)}$ be an ordered exceptional dissection arising from $\zD_{(u,s)}$ such that $\za_{s+1}$ directly follows $\zg$. This is possible, since for any $\za\in \zD_{(u,s)}$ with $\zg \preceq \za$, we have $\za_{s+1}=\zb_{u+1} \preceq \za$ by the above fact. Assume $\zg$ is at the $i$-th position in $\zG_{(u,s)}$, then $\sigma^{-1}_i \zG_{(u,s)}$ is an ordered exceptional dissection obtained from $\zG_{(u,s)}$ by exchanging the positions of $\zg$ and $\za_{s+1}$, and then replacing $\za_{s+1}$ by $L_{\zg}\za_{s+1}$. Then $\zG_{(u',s')}:=\sigma^{-1}_i \zG_{(u,s)}$ is what we want, noticing that on the one hand $(u',s')=(u+1,s+1)$, on the other hand $\overline{L_{\zg}\za_{s+1}}=\overline{\za_{s+1}}$ and thus the underlying set of $\zG_{(u',s')}$ is a lift of  $\zD_\zg$.

Case 2. $\za_{s+1}\neq\zb_{u+1}$. Then at most one of $\za_{s+1}\preceq \zb_{u+1}$ and $\zb_{u+1}\preceq\za_{s+1}$ is satisfied. Otherwise, we have $\overline{\zb}_{u+1}\preceq\overline{\za}_{s+1}\preceq\overline{\zb}_{u+1}$, which contradicts  the fact that $\zD_\zg$ is exceptional. Without loss of generality, assume that $\zb_{u+1}\preceq\za_{s+1}$ does not hold.
Then for $\za\in \zD$ such that $\zg\preceq\za$, we have $\za_{s+1}\preceq \za$.
Otherwise, if $\za\preceq\za_{s+1}$, then by the above fact, we have $\zb_{u+1}\preceq \za$, thus $\zb_{u+1}\preceq \za_{s+1}$. A contradiction. Then similar to the first case, we define $\zG_{(u,s)}$ as an ordered exceptional dissection arising from $\zD_{(u,s)}$ such that $\za_{s+1}$ directly follows $\zg$. And then define $\zG_{(u',s')}$ as $\sigma^{-1}_i \zG_{(u,s)}$, where $(u',s')=(u,s+1)$.
\end{proof}

The next corollary follows from the two previous lemmas.

\begin{corollary}\label{corollary:braid mutate same underlying set2}
Let $\zG_1$ and $\zG_2$ be two ordered exceptional dissections on $(S,M)$ containing $\zg$. Let $\overline{\zG}_1$ and $\overline{\zG}_2$  be the  ordered exceptional dissections on $(S_\zg, M_\zg)$ induced by  $\zG_1$ and $\zG_2$, respectively. Suppose $\overline{\zG}_1$ and $\overline{\zG}_2$  have the same underlying set $\zD_\zg$.   Then there is an element $\sigma$ in $B_n$ such that $\sigma \zG_1=\zG_2$.
\end{corollary}
\begin{proof}
Denote by $\zD_{(v,t)}$ the unique exceptional dissection on $(S,M)$ such that $\zD_\zg=\overline{\zD}_{(v,t)}$ and $\zg$ is a maximal element in $\zD_{(v,t)}$. Let $\zG_{(v,t)}$ be any ordered exceptional dissection arising from $\zD_{(v,t)}$.
Denote by $\zD_{(u_1,s_1)}$ and $\zD_{(u_2,s_2)}$ the underlying sets of $\zG_1$ and $\zG_2$ respectively.
Then by alternatively using Lemma \ref{lemma:braid mutate same underlying set} and Lemma \ref{lemma:braid mutate same collection}, for $1 \leq i \leq 2,$ we have an element $\sigma'_i$ in $B_n$ and an ordered exceptional dissection
$\zG_i'$ on $(S,M)$ with underlying set $\zD_{(u_i,s_i)}$, such that $\sigma'_i\zG_i'=\zG_{(v,t)}$.
On the other hand, by Lemma \ref{lemma:braid mutate same collection}, for $1 \leq i \leq 2,$ we have an element $\sigma_i$ in $B_n$ such that $\sigma_i\zG_i=\zG'_i$. Then by setting $\sigma=\sigma_2^{-1}{\sigma'}_2^{-1}\sigma_1'\sigma_1$, we have $\sigma\zG_1=\zG_2$.
\end{proof}

The next lemma sets up the induction step for the proof of Theorem~\ref{theorem:transitivity2}.

\begin{lemma}\label{lemma:braid mutate share gamma}
Let $\zG$ and $\zG'$ be two ordered exceptional dissections on $(S,M)$ containing $\zg$. If $B_{n-1}$ is transitive on the set of ordered exceptional dissections on $(S_\zg,M_\zg)$, then there is an element $\sigma$ in $B_n$ such that $\sigma \zG=\zG'$.
\end{lemma}
\begin{proof}
At first, we assume that $\zg$ is a maximal element in both $\zG$ and $\zG'$. Then due to Lemma \ref{lemma:braid mutate same collection}, without loss of generality, we assume that $\zg$ is the last element in $\zG$ and $\zG'$. Note that   by Lemma \ref{lemma:completione} the induced ordered sets $\overline{\zG}$ and $\overline{\zG}'$ are both ordered exceptional dissections on $(S_\zg,M_\zg)$. So by  assumption, there is an element $\sigma=\delta_{m}\cdots\delta_{1}$ in $B_{n-1}$ such that $\sigma\overline{\zG}=\overline{\zG}'$, where each $\delta_j$ is some generator $\sigma_i$ or $\sigma^{-1}_i$ of $B_{n-1}$ for some $1 \leq i \leq n-2$. By Corollary \ref{corollary:completion2}, there is an unique ordered exceptional dissection $\zG_1$ on $(S,M)$ such that $\overline{\zG}_1=\delta_{1}\overline{\zG}$ and such that $\zg$ is its last element. By viewing $\delta_1$ as an element in $B_n$ by the natural embedding of braid groups $B_{n-1}\subseteq B_n$, the set $\delta_{1}\zG$ is an ordered exceptional dissection on $(S,M)$. Note that $\zG_1$ and $\delta_{1}\zG$ are both  lifts of $\delta_{1}\overline{\zG}$ containing $\zg$ as their last element. Thus
 we have $\zG_1=\delta_{1}\zG$.
Similarly, we construct ordered exceptional dissections $\zG_j, 2 \leq j \leq m,$ on $(S,M)$ such that $\overline{\zG}_j=\delta_{j}\cdots\delta_{1}\overline{\zG}$ and $\zG_j=\delta_{j}\cdots\delta_{1}\zG$. Then we have
$$\zG'=\delta_{m}\zG_{m-1}=\delta_{m}\delta_{m-1}\zG_{m-2}=\cdots=
\delta_{m}\cdots\delta_{1}\zG=\sigma\zG,$$
where $\delta_{j}, 1 \leq j \leq m$ and $\sigma$ are viewed as elements in $B_n$.

Finally, for the general case, that is, the case when $\zg$ is not necessary  a maximal element in $\zG$ or $\zG'$, the proof follows from the above  and Corollary \ref{corollary:braid mutate same underlying set2}.
\end{proof}

Next we formulate a condition on a marked surface such that for surfaces satisfying this condition, the action of the braid group is transitive on the set of ordered exceptional dissections.
We say that a marked surfaces $(S,M)$ satisfies the `reachable condition for exceptional arcs' if the following holds.

\begin{definition}[\bf
Condition RCEA]\label{condition:rcea}  We say that a marked surface $(S,M)$ satisfies the {\it {\bf RCEA} condition} if for any two exceptional arcs $\zg$ and $\zg'$ on $(S,M)$, there is a sequence of exceptional arcs $\zg_0,\cdots,\zg_{m+1}$ on $(S,M)$, such that each pair $\{\zg_i,\zg_{i+1}\}$ is an exceptional collection, where  $0 \leq i \leq m$, $\zg_0=\zg$, and  $\zg_{m+1}=\zg'$.
\end{definition}

\begin{theorem}\label{proposition:transitivity}
Suppose that the  {\bf RCEA} condition holds for all marked surfaces.  Then for any marked surface with exceptional dissections, the braid group action on the set of ordered exceptional dissections is transitive.
\end{theorem}
\begin{proof}
Let $(S,M)$ be a connected marked surface. We prove the result by using induction on the number $n=|M^{\gpoint}|+2g+b-2$ of elements in an exceptional dissection on $(S,M)$.  The base case of the induction corresponds to marked surfaces which are unions of non-punctured-disks each with two $\gpoint$-marked points, where the statement is clearly true.
We now show the induction step. Assume the result holds  for  surfaces with less than $n-1$ elements in an exceptional dissection.

Let $\zG$ and $\zG'$ be any two ordered exceptional dissections on $(S,M)$. We will show that there is an element $\sigma$ in $B_n$ such that $\sigma\zG=\zG'$.
Let $\zg$ and $\zg'$ be any two arcs in $\zG$ and $\zG'$ respectively. Then by the {\bf RCEA} condition, there is a sequence of exceptional arcs $\zg_0,\cdots,\zg_{m+1}$ on $(S,M)$, such that each pair $\{\zg_i,\zg_{i+1}\}$ is an exceptional collection, where  $0 \leq i \leq m$ and $\zg_0=\zg$, $\zg_{m+1}=\zg'$. We have two cases.

Case 1. $(S,M)$ is not a surface of the form $\mathbb{T}_{(g,1,2)}$. Note that in this case, for each $1 \leq i \leq m$, the cut surface $(S_{\zg_i},M_{\zg_i})$ does not have a connected component of the form $\mathbb{T}_{(g,1,1)}$. Then by Proposition \ref{proposition:completion1}, there exists an exceptional dissection $\zG_i$ on $(S,M)$ which contains $\zg_i$ and $\zg_{i-1}$ for each $1 \leq i \leq m$. Furthermore, by Lemma \ref{lemma:braid mutate share gamma}, there is an element $\sigma_i, 0\leq i \leq m$ in $B_n$ such that $\sigma_i\zG_i=\zG_{i+1}$, where we denote $\zG$ by $\zG_0$ and denote $\zG'$ by $\zG_{m+1}$ for simplicity. Then we have $\sigma\zG=\zG'$, where $\sigma=\sigma_m\cdots\sigma_0\in B_n$.

Case 2. $(S,M)$ is a surface of the form $\mathbb{T}_{(g,1,2)}$. Note that in this case, since each $\zg_i, 1 \leq i \leq m$ is contained in an exceptional pair, it can not be an arc homotopic to a boundary segment which contains a $\rpoint$-marked point, see in Figure \ref{figure:completion1}. Thus there exists an exceptional dissection $\zG_i$ which contains $\zg_i$ for each $1 \leq i \leq m$. Then similar to the proof of the first case, we have an element $\sigma$ in $B_n$ such that $\sigma\zG=\zG'$.

\end{proof}

In general, it seems a difficult problem to determine, for which surfaces the \textbf{RCEA} condition holds.  In the following, we show that the \textbf{RCEA} condition does hold for surfaces of genus zero.

\begin{theorem}\label{theorem:transitivity}
The {\bf RCEA} condition is satisfied for any marked surface of zero genus. \end{theorem}
\begin{proof}
Let $\zg$ and $\zg'$ be two exceptional arcs on $(S,M)$. Assume that $(\zg,\zg')$ is not an exceptional collection, otherwise there is nothing to prove. Then $\zg$ and $\zg'$ share two endpoints or they have at least one intersection in the interior of $S$. Suppose first that  $\zg$ and $\zg'$ share two endpoints and that they do not intersect in the interior. Then they surround a bigon which contains at least one boundary component, see in Figure \ref{figure:configuration of intersections0}, noticing that the genus of the surface is zero. Then there exists an exceptional arc  $\zg_1$  which connects an endpoint of $\zg$ and a $\gpoint$-marked point on a boundary component in the bigon, as depicted in Figure \ref{figure:configuration of intersections0}. Then $(\zg,\zg_1)$ and $(\zg_1,\zg')$ are both exceptional collections.

\begin{figure}
\begin{center}
{\begin{tikzpicture}[scale=0.8]

\begin{scope}
	\draw (2.5,0) circle (.3cm);
	\clip[draw] (2.5,0) circle (.3cm);
    \draw[]  (1.5,1-1)--(3.5,1-1);
    \draw[]  (1.5,.85-1)--(3.5,.85-1);
    \draw[]  (1.5,1.15-1)--(3.5,1.15-1);
\end{scope}

	\draw[green] (0,0) circle [radius=0.07];
	\draw[green] (5,0) circle [radius=0.07];
	\draw[green] (2.2,0) circle [radius=0.07];
    \draw[](0,0)..controls (1.5,1) and (3.5,1)..(5,0);
	\draw (2.5,1.1) node {\tiny$\zg$};		
    \draw (2.5,-1.1) node {\tiny$\zg'$};
	\draw (1.2,0.2) node {\tiny$\zg_1$};		

	\draw[](0,0)..controls (1.5,-1) and (3.5,-1)..(5,0);
    \draw[]  (0,0)--(2.2,0);
	\draw[bend right](-.2,-.6)to(-.2,.6);
	\draw[-]  (-0.1,-.33)to(-.4,-.6);
	\draw[-]  (-.05,-.1)to(-.4,-.4);
	\draw[-]  (-.05,.1)to(-.4,-.2);
	\draw[-]  (-.05,.3)to(-.4,0);
	\draw[-]  (-.1,.45)to(-.4,0.2);

	\draw[bend left](5.2,-.6)to(5.2,.6);
	\draw[-]  (5.35,-.35)to(5.15,-.5);
	\draw[-]  (5.35,-.15)to(5.1,-.35);
	\draw[-]  (5.35,.05)to(5.05,-.2);
	\draw[-]  (5.35,.25)to(5.05,0);
	\draw[-]  (5.35,.45)to(5.05,0.2);

\draw (-.7,0) node {\tiny$p$};
	\end{tikzpicture}}
\end{center}
\begin{center}
\caption{The exceptional arcs $\zg$ and $\zg'$ share two endpoints and $(\zg,\zg')$ is a non-exceptional cycle, where the  left and right boundary components may coincide.}\label{figure:configuration of intersections0}
\end{center}
\end{figure}
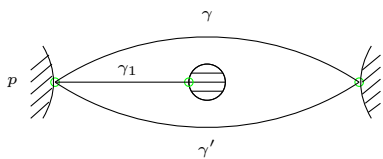

In the following, we assume that $\zg$ and $\zg'$ intersect in the interior of the surface.  Note that since the genus of the surface is zero there are three possible cases as shown in  Figures \ref{figure:configuration of intersections1}, \ref{figure:configuration of intersections2}, and \ref{figure:configuration of intersections3}.

\emph{Case 1. All the endpoints of $\zg$ and $\zg'$ are on the same boundary component.}
Then the local configuration of the intersection nearest to an endpoint $p$ of $\zg'$ is depicted in Figure \ref{figure:configuration of intersections1}. We choose $\zg_1$ as showed in the figure, then $(\zg,\zg_1)$ and $(\zg_1,\zg')$ are both exceptional pairs.

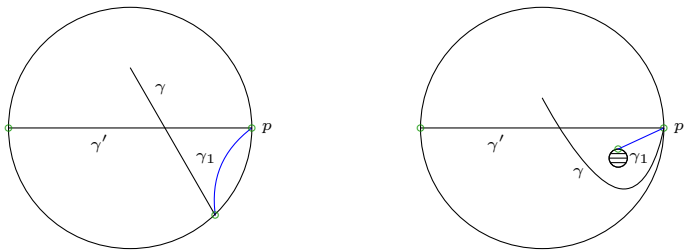
\begin{figure}[ht]
\begin{center}
\begin{tikzpicture}[scale=0.4]

    \draw[]  (-4,0)--(4,0);
    \draw[]  (2.8,-2.9)--(0,2);
    \draw[bend left,blue]  (2.8,-2.9)to(4,0);
	\draw (0,0) circle (4cm);
    \draw (4,0) node {\tiny$\gpoint$};
    \draw (2.8,-2.9) node {\tiny$\gpoint$};
    \draw (-4,0) node {\tiny$\gpoint$};
    \draw (-1,-.5) node {\tiny$\zg'$};
    \draw (1,1.3) node {\tiny$\zg$};
    \draw (2.5,-1) node {\tiny$\zg_1$};
    \draw (4.5,0) node {\tiny$p$};

\end{tikzpicture}
\begin{tikzpicture}[scale=0.4]
    \draw (-8,0) node {};
\begin{scope}
	\draw (2.5,-1) circle (.3cm);
	\clip[draw] (2.5,-1) circle (.3cm);
    \draw[]  (1.5,-1)--(3.5,-1);
    \draw[]  (1.5,-.85)--(3.5,-.85);
    \draw[]  (1.5,-1.15)--(3.5,-1.15);
\end{scope}

    \draw[]  (-4,0)--(4,0);
	\draw (0,0) circle (4cm);
    \draw (4,0) node {\tiny$\gpoint$};
    \draw (-4,0) node {\tiny$\gpoint$};
    \draw (2.5,-.7) node {\tiny$\gpoint$};
    \draw[blue]  (2.5,-.7)--(4,0);
    \draw (-1.5,-.5) node {\tiny$\zg'$};
    \draw (1.2,-1.5) node {\tiny$\zg$};
    \draw (3.2,-1) node {\tiny$\zg_1$};
    \draw plot [smooth,tension=1] coordinates {(4,0)(2.5,-2)(0,1)};
    \draw (4.5,0) node {\tiny$p$};

\end{tikzpicture}
\end{center}
\begin{center}
\caption{Local configurations of the intersection nearest to an endpoint $p$ of $\zg'$ on a surface of genus zero, where the endpoints of the two exceptional arcs $\zg$ and $\zg'$ are on the same boundary component.}\label{figure:configuration of intersections1}
\end{center}
\end{figure}

\emph{Case 2. Both endpoints of $\zg$ (resp. $\zg'$) are on the same boundary component $B$ (resp. $B'$), where $B$ and $B'$ are different.}
By assumption  $\zg$ and $\zg'$ have an interior intersection which implies that there are at least three boundary components on $(S,M)$. Then up to homotopy, the local configuration of the intersection  nearest to one of the endpoints  of $\zg$, which will denote by  $p$, is as in Figure \ref{figure:configuration of intersections2}.
We choose $\zg_1$ and $\zg_2$ as shown in  Figure  \ref{figure:configuration of intersections2}. Then $(\zg,\zg_1)$, $(\zg_1,\zg_2)$ and $(\zg_2,\zg')$ are exceptional pairs.

\begin{figure}[ht]
\begin{center}
\begin{tikzpicture}[scale=0.4]
\begin{scope}
	\draw (2.5-3.5,-1) circle (.3cm);
	\clip[draw] (2.5-3.5,-1) circle (.3cm);
    \draw[]  (1.5-3.5,-1)--(3.5-3.5,-1);
    \draw[]  (1.5-3.5,-.85)--(3.5-3.5,-.85);
    \draw[]  (1.5-3.5,-1.15)--(3.5-3.5,-1.15);
\end{scope}
\draw (-1,-.7) node {\tiny$\gpoint$};

\begin{scope}
	\draw (-2.5,1) circle (.3cm);
	\clip[draw] (-2.5,1) circle (.3cm);
    \draw[]  (-1.5,1)--(-3.5,1);
    \draw[]  (-1.5,.85)--(-3.5,.85);
    \draw[]  (-1.5,1.15)--(-3.5,1.15);
\end{scope}
    \draw (-2.5,.7) node {\tiny$\gpoint$};

    \draw[]  (-4,0)--(4,0);
	\draw (0,0) circle (4cm);
    \draw (4,0) node {\tiny$\gpoint$};
    \draw (-4,0) node {\tiny$\gpoint$};
    \draw (-4.5,0) node {\tiny$p$};
    \draw[blue]  (-2.5,.7)--(-1,-.7);
    \draw[blue]  (4,0)--(-1,-.7);
    \draw (1.5,.5) node {\tiny$\zg'$};
    \draw (1.5,-.8) node {\tiny$\zg_2$};
    \draw (-1.7,.5) node {\tiny$\zg_1$};
    \draw (-2.5,-.5) node {\tiny$\zg$};
    \draw plot [smooth,tension=1] coordinates {(-2.5,.7)(-1,-2)(0,1)};
    \draw (0,3.5) node {\tiny$B'$};
    \draw (-2.5,1.8) node {\tiny$B$};
\end{tikzpicture}
\end{center}
\begin{center}
\caption{Local configuration of the intersection of two exceptional arcs $\zg$ and $\zg'$ nearest to an endpoint $p$ of $\zg$ on a surface of genus zero, where the endpoints of $\zg$ and $\zg'$ are on different boundary components $B$ and $B'$ respectively.}\label{figure:configuration of intersections2}
\end{center}
\end{figure}
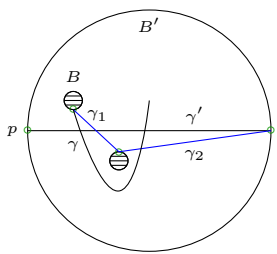

\emph{Case 3. There is at least one arc, for example $\zg$, with endpoints on two different boundary components. Denote one of these by $B$.}
Assume that locally  $\zg$ and $\zg'$ form a bigon containing at least one boundary component as shown in the left picture of Figure \ref{figure:configuration of intersections3}. We construct a new arc $\zg_1$ obtained from $\zg$ by smoothing the intersections of $\zg$ and $\zg'$  at this bigon, as shown in the left picture of  Figure \ref{figure:configuration of intersections3}. Then there are no interior intersections of $\zg$ and $\zg_1$, and $|\zg_1\cap\zg'|=|\zg\cap\zg'|-2$. But $\zg_1$ and $\zg$ do not form an exceptional pair, since $\zg$ and $\zg_1$ intersect each other at both endpoints in such a way that they form a non-exceptional cycle.

Now we  rotate $\zg_1$ anti-clockwise one  time around the boundary component $B$. Then $\zg_1$ and $\zg$ form an exceptional pair,  noticing that the endpoints of $\zg$ are on two different boundary components and that in this case this process does not give rise to   interior intersections of these arcs. On the other hand, such a rotation may produce new intersections between $\zg_1$ and $\zg'$. However, we may assume that  at most one endpoint of $\zg'$ lies on $B$, otherwise we choose the boundary component of the other endpoint of $\zg$. Under such assumption, there may appear at most one new intersection between $\zg_1$ and $\zg'$. So that we have $|\zg_1\cap\zg'|\leq|\zg\cap\zg'|-1$, and $(\zg,\zg_1)$ is an exceptional pair.

Then we replace $\zg$ by $\zg_1$, and continue this process until there are no bigons surrounded by $\zg'$ and  the successive replacements of  $\zg$. By abuse of notation, we will denote this new arc replacing $\zg$ again by $\zg$.  Then there are two cases. Either $\zg$ and $\zg'$ share at most one endpoint and an argument as in case 1 finishes the proof or $\zg$ and $\zg'$ share both endpoints. But then, since the genus of the surface is zero, the only possible intersections between $\zg$ and $\zg'$ come from rotations of one of the arcs, for example $\zg$, around a boundary component which is incident  to the other arc $\zg'$,  see the right picture in Figure \ref{figure:configuration of intersections3}. Then define $\zg_1$ as the new arc which is obtained from $\zg$ by a local surgery, where locally we replace the shown segment of $\gamma$ by the blue segment as in the figure. Note that we have
$|\zg_1\cap\zg'|=|\zg\cap\zg'|-1$, and $(\zg,\zg_1)$ is an exceptional pair. Continue this process until there are no intersections between the new arc and $\zg'$. Then we have a sequence of exceptional arcs as desired.

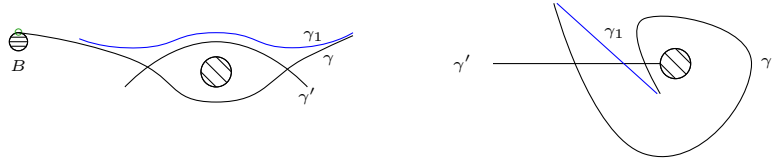
\begin{figure}[ht]
\begin{center}
\begin{tikzpicture}[scale=0.4]
\begin{scope}
	\draw (-6.5,1) circle (.3cm);
	\clip[draw] (-2.5-4,1) circle (.3cm);
    \draw[]  (-1.5-4,1)--(-3.5-4,1);
    \draw[]  (-1.5-4,.85)--(-3.5-4,.85);
    \draw[]  (-1.5-4,1.15)--(-3.5-4,1.15);
\end{scope}
    \draw (-2.5-4,1.3) node {\tiny$\gpoint$};
    \draw (-2.5-4,0.2) node {\tiny$B$};
\begin{scope}
	\draw (0,0) circle (.5cm);
	\clip[draw] (0,0) circle (.5cm);
	\foreach \x in {-4,-3.5,-3,-2.5,-2,-1.5,-1,-0.5,0,0.5,1,1.5,2,2.5,3,3.5,4}	\draw[xshift=\x cm]  (-5,5)--(5,-5);
\end{scope}
     \draw (0,-3) node {};
    \draw (3,-.8) node {\tiny$\zg'$};
    \draw (3.7,.4) node {\tiny$\zg$};
     \draw (3.2,1.2) node {\tiny$\zg_1$};
    \draw plot [smooth,tension=1] coordinates {(-2.5-4,1.3)(-3,.5)(0,-1)(3,.5)(4.5,1.2)};
    \draw plot [smooth,tension=1] coordinates {(-3,-.5)(0,1)(3,-.5)};
    \draw [blue] plot [smooth,tension=1] coordinates {(-4.5,1.1)(-2.5,.8)(0,1.3)(2.5,.8)(4.5,1.3)};

\end{tikzpicture}
\begin{tikzpicture}[scale=0.4]
\begin{scope}
	\draw (0,0) circle (.5cm);
	\clip[draw] (0,0) circle (.5cm);
	\foreach \x in {-4,-3.5,-3,-2.5,-2,-1.5,-1,-0.5,0,0.5,1,1.5,2,2.5,3,3.5,4}	\draw[xshift=\x cm]  (-5,5)--(5,-5);
\end{scope}

    \draw (-10,0) node {};
    \draw (-7,0) node {\tiny$\zg'$};
    \draw (3,0) node {\tiny$\zg$};
    \draw (-2,1) node {\tiny$\zg_1$};
    \draw[-] (-.5,0)to(-6,0);
    \draw plot [smooth,tension=1] coordinates {(-4,2)(-1,-3)(2.5,0)(-1,1.5)(-.5,-1)};
    \draw[-,blue] (-3.9,2)to(-.6,-1);
\end{tikzpicture}
\end{center}
\begin{center}
\caption{Local configurations of intersections between two exceptional arcs on a surface of genus zero, where in the right picture the arc $\zg$ may also surround the boundary component in another direction.}\label{figure:configuration of intersections3}
\end{center}
\end{figure}

\end{proof}

Combining Theorems~\ref{proposition:transitivity} and ~\ref{theorem:transitivity}, we obtain the following.

\begin{corollary}
The action of the braid group on the set of ordered exceptional dissections on a surface of genus zero is transitive.

\end{corollary}

\subsection{Braid group action on full exceptional sequences in $\ka$} In the following, by using the geometric model established above, we will show how the braid group $B_n$ acts on the full exceptional sequences in the perfect derived category $\ka$ of a gentle algebra $A$. This allows us to show the transitivity of the braid group action when the associated surface is of genus zero.

We start by recalling the braid group action on a triangulated category $\mathcal{T}$.
Let $(X,Y)$ be an exceptional pair in $\calt$.
Define objects $R_{Y}X$ and $L_{X}Y$ in $\calt$ through the following triangles
\begin{equation}\label{equation:right exchange braid}
\xymatrix{
X\ar[r]&\coprod_{\ell\in\mathbb{Z}}D\Hom_{\calt}(X,Y[\ell])\otimes_kY[\ell]
\ar[r]&R_{Y}X[1]\ar[r]&X[1]
}\end{equation}
\begin{equation}\label{equation:left exchange braid}
\xymatrix{
Y[-1]\ar[r]& L_{X}Y[-1]\ar[r]&\coprod_{\ell\in\mathbb{Z}}\Hom_{\calt}(X[\ell],Y)\otimes_kX[\ell]\ar[r]&Y.
}\end{equation}
Then $(Y,R_{Y}X)$ and $(L_{X}Y,X)$ are again exceptional pairs in $\calt$.
We call $R_{Y}X$  the {\em{right mutation}} of $X$ at $Y$, and and $L_{X}Y$ the {\em{left mutation}} of $Y$ at $X$, where the two distinguished triangles \eqref{equation:right exchange braid} and \eqref{equation:left exchange braid} are called \emph{exchange triangles}.

We denote by $\exc\calt$  the set of isomorphism classes of full exceptional sequences in $\calt$, where the length of each full exceptional sequence is $n$ which is the rank of the Grothendieck group $K_0(\calt)$.

Then by \cite{GR87}, $B_n$ acts on $\exc\calt$  as follows:
For a full exceptional sequence ${\mathbf X}:=(X_1,\cdots,X_n)$ and $1\le i<n$, set
$$\sigma_i{\mathbf X} := (X_1,\cdots,X_{i-1},X_{i+1},R_{X_{i+1}}X_{i},X_{i+2},\cdots,X_n)$$
$$\sigma_i^{-1}{\mathbf X} := (X_1,\cdots,X_{i-1},L_{X_i}X_{i+1},X_{i},X_{i+2},\cdots,X_n).$$

Let $(X_1,\cdots,X_n)$ be a full exceptional sequence in $\ka$ where,  for all $1 \leq i \leq n$, $X_i=\P_{(\zg_i,f_{\zg_i})}$ for some graded $\gpoint$-arc $(\zg_i,f_{\zg_i})$ in the surface model $(S,M,\zD_A)$ of $A$. Then by Proposition \ref{prop:tree type} the set $(\zg_1,\cdots,\zg_n)$ is an ordered exceptional dissection. In particular, for any $1 \leq i \leq n-1$, $(\zg_i,\zg_{i+1})$ is an exceptional pair. So the possible relative positions of $\zg_i$ and $\zg_{i+1}$ are shown in Figure \ref{figure:mutation of exceptional dissection}, where we view $\zg_1$ as $\zg_i$ and view $\zg_2$ as $\zg_{i+1}$.

Note that for the right mutation, we need to  consider  the maps from $X_i$ to $X_{i+1}[m]$ for any integer $m$, while for the left mutation, we need to consider all the maps from $X_i[m]$ to $X_{i+1}$, see the exchange triangles \eqref{equation:right exchange braid} and \eqref{equation:left exchange braid} respectively. So in both cases we need to associate gradings to the arcs. For the sake of clarity, we redraw and relabel the picture in Figure \ref{figure:mutation of exceptional dissection} as four pictures in Figures \ref{figure:right action} and  \ref{figure:left action}, where $s$ and $t$ ($s'$ and $t'$ respectively) denote the intersections of the $\gpoint$-arcs $R_{\zg_{i+1}\zg_i}$ and $\zg_i$ ($L_{\zg_i}\zg_{i+1}$ and $\zg_{i+1}$ respectively)
with the initial dual $\rpoint$-dissection $\zD^*_A$ closest to the intersection of $R_{\zg_{i+1}\zg_i}$ and $\zg_i$ ($L_{\zg_i}\zg_{i+1}$ and $\zg_{i+1}$ respectively).

\begin{figure}[H]
\begin{center}
{\begin{tikzpicture}[scale=0.4]
\draw[green] (0,0) circle [radius=0.2];
\draw[green] (0,5) circle [radius=0.2];
\draw[green] (5,0) circle [radius=0.2];

\draw[-] (0,0.2) -- (0,4.8);

\draw[-] (0.2,0) -- (4.8,0);

\draw[-,blue] (0.15,4.85) -- (4.85,0.15);

\draw (3.7,2.6) node {\tiny$R_{\zg_{i+1}\zg_i}$};

\draw (2.3,-.5) node {\tiny$\zg_i$};

\draw (-1.2,2.5) node {\tiny$\zg_{i+1}$};

\draw (-.6,5.5) node {\tiny$q_2$};
\draw (-.5,-.5) node {\tiny$q_1$};
\draw (5.5,-.5) node {\tiny$q_3$};

\draw (0.9,0.9) node {\tiny$a$};
\draw[->] (0.7,0) to [out=90,in=0] (0,0.7);

\draw (2.5,-2) node {case I};

\draw[-,red] (4.2,-0.2) -- (3.8,.3);
\draw[-,red] (3.8,0.8) -- (4.8,1);
\draw (4.3,1.4) node {\tiny$s$};
\draw (3.7,-.6) node {\tiny$t$};

\end{tikzpicture}}
\qquad\qquad
{\begin{tikzpicture}[scale=0.4]
\draw[green] (0,0) circle [radius=0.2];
\draw[green] (0,5) circle [radius=0.2];
\draw[green] (5,5) circle [radius=0.2];
\draw[green] (5,0) circle [radius=0.2];

\draw[-] (0,0.2) -- (0,4.8);
\draw[-] (5,0.2) -- (5,4.8);

\draw[-] (0.2,0) -- (4.8,0);
\draw[-] (0.2,5) -- (4.8,5);

\draw[-] (0.15,0.15) -- (4.85,4.85);
\draw[-,blue] (0.15,4.85) -- (4.85,0.15);

\draw (1.9,1.3) node {\tiny$\zg_i$};
\draw (2,4.2) node {\tiny$R_{\zg_{i+1}\zg_i}$};

\draw (6.2,2.5) node {\tiny$\zg_{i+1}$};
\draw (-1.2,2.5) node {\tiny$\zg_{i+1}$};

\draw (-.6,5.5) node {\tiny$q_2$};
\draw (-.5,-.5) node {\tiny$q_1$};
\draw (5.5,-.5) node {\tiny$q_1$};
\draw (5.5,5.5) node {\tiny$q_2$};

\draw (0.6,1.3) node {\tiny$a_1$};
\draw (4.5,3.7) node {\tiny$a_2$};
\draw[->] (0.5,0.5) to [out=135,in=0] (0,0.7);
\draw[->] (4.5,4.5) to [out=315,in=180] (5,4.2);

\draw[-,red] (3.5,3) -- (2.5,3);
\draw[-,red] (2.8,1.8) -- (3.8,2);
\draw (3.1,1.4) node {\tiny$s$};
\draw (3.1,3.6) node {\tiny$t$};

\draw (2.5,-2) node {case II};
\end{tikzpicture}}
\end{center}
\begin{center}
\caption{Possible relative positions of $\zg_{i}$ and $\zg_{i+1}$ with regards to the right mutation.}\label{figure:right action}
\end{center}
\end{figure}
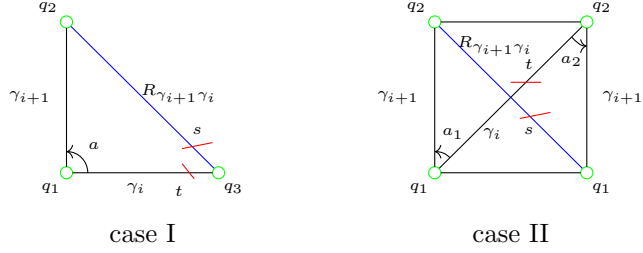

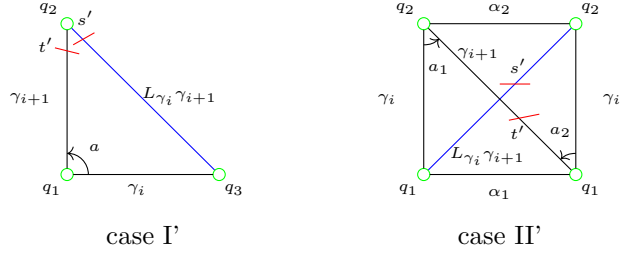
\begin{figure}[H]
\begin{center}
{\begin{tikzpicture}[scale=0.4]
\draw[green] (0,0) circle [radius=0.2];
\draw[green] (0,5) circle [radius=0.2];
\draw[green] (5,0) circle [radius=0.2];

\draw[-] (0,0.2) -- (0,4.8);

\draw[-] (0.2,0) -- (4.8,0);

\draw[-,blue] (0.15,4.85) -- (4.85,0.15);

\draw (3.7,2.6) node {\tiny$L_{\zg_i}\zg_{i+1}$};

\draw (2.3,-.5) node {\tiny$\zg_i$};

\draw (-1.2,2.5) node {\tiny$\zg_{i+1}$};

\draw (-.6,5.5) node {\tiny$q_2$};
\draw (-.5,-.5) node {\tiny$q_1$};
\draw (5.5,-.5) node {\tiny$q_3$};

\draw (0.9,0.9) node {\tiny$a$};
\draw[->] (0.7,0) to [out=90,in=0] (0,0.7);

\draw (2.5,-2) node {case I'};

\draw[-,red] (.2,4.3) -- (.9,4.7);
\draw[-,red] (-.4,4.2) -- (.4,4);
\draw (.6,5.2) node {\tiny$s'$};
\draw (-.7,4.3) node {\tiny$t'$};

\end{tikzpicture}}
\qquad\qquad
{\begin{tikzpicture}[scale=0.4]
\draw[green] (0,0) circle [radius=0.2];
\draw[green] (0,5) circle [radius=0.2];
\draw[green] (5,5) circle [radius=0.2];
\draw[green] (5,0) circle [radius=0.2];

\draw[-] (0,0.2) -- (0,4.8);
\draw[-] (5,0.2) -- (5,4.8);

\draw[-] (0.2,0) -- (4.8,0);
\draw[-] (0.2,5) -- (4.8,5);

\draw[-,blue] (0.15,0.15) -- (4.85,4.85);
\draw[-] (0.15,4.85) -- (4.85,0.15);

\draw (2.1,.6) node {\tiny$L_{\zg_i}\zg_{i+1}$};
\draw (1.8,4) node {\tiny$\zg_{i+1}$};
\draw (2.5,5.5) node {\tiny$\za_2$};
\draw (2.5,-.65) node {\tiny$\za_1$};

\draw (6.2,2.5) node {\tiny$\zg_{i}$};
\draw (-1.2,2.5) node {\tiny$\zg_{i}$};

\draw (-.6,5.5) node {\tiny$q_2$};
\draw (-.5,-.5) node {\tiny$q_1$};
\draw (5.5,-.5) node {\tiny$q_1$};
\draw (5.5,5.5) node {\tiny$q_2$};

\draw (0.5,3.5) node {\tiny$a_1$};
\draw (4.5,1.4) node {\tiny$a_2$};

\draw[<-] (0.5,4.5) to [out=225,in=0] (0,4.3);
\draw[->] (5,0.7) to [out=180,in=45] (4.5,.5);

\draw[-,red] (3.5,3) -- (2.5,3);
\draw[-,red] (2.8,1.8) -- (3.8,2);
\draw (3.1,1.4) node {\tiny$t'$};
\draw (3.1,3.6) node {\tiny$s'$};

\draw (2.5,-2) node {case II'};
\end{tikzpicture}}
\end{center}
\begin{center}
\caption{Possible relative positions of $\zg_{i}$ and $\zg_{i+1}$ with regards to the left mutation.}\label{figure:left action}
\end{center}
\end{figure}

\begin{theorem}\label{thm:braid group actions}
Let $A$ be a gentle algebra with associated  surface model $(S,M,\zD_A)$. Let $\mathbf X=(X_1,\cdots,X_n)$ be a full exceptional sequence in $\ka$, and assume that $X_i=\P_{(\zg_i,f_{\zg_i})}$ for some graded $\gpoint$-arc $(\zg_i,f_{\zg_i})$, for all  $1\leq i\leq n$.

(1) If  $\zg_i$ and $\zg_{i+1}$ do not intersect then $R_{X_{i+1}}X_i=X_i$ and $L_{X_{i}}X_{i+1}=X_{i+1}$. So
 $$\sigma_i{\mathbf X} = \sigma_i^{-1}{\mathbf X} = (X_1,\cdots,X_{i-1},X_{i+1},X_i,X_{i+2},\cdots,X_n),$$
that is, the actions of $\sigma_i$ and $\sigma_i^{-1}$ change the order of $X_i$ and $X_{i+1}$.

(2) Suppose that $\zg_i$ and $\zg_{i+1}$ intersect. Then
\begin{itemize}
\item $R_{X_{i+1}}X_{i}=\P_{({R_{\zg_{i+1}}\zg_{i}},f_{R_{\zg_{i+1}}\zg_{i}})}$ with $f_{R_{\zg_{i+1}}\zg_{i}}(s)=f_{\zg_i}(t)$, where the arc ${R_{\zg_{i+1}}\zg_{i}}$ and the intersections $s$ and $t$ are depicted in Figure \ref{figure:right action}.
\item
 $L_{X_i}X_{i+1}=\P_{({L_{\zg_i}\zg_{i+1}},f_{L_{\zg_i}\zg_{i+1}})}$ with $f_{L_{\zg_i}\zg_{i+1}}(s')=f_{\zg_{i+1}}(t')$, where the arc ${L_{\zg_i}\zg_{i+1}}$ and the intersections $s'$ and $t'$ are depicted in Figure \ref{figure:left action}.
 \end{itemize}
 In this case,
  $$\sigma_i{\mathbf X} = (X_1,\cdots,X_{i-1},X_{i+1},\P_{({R_{\zg_{i+1}}\zg_{i}},f_{R_{\zg_{i+1}}\zg_{i}})},X_{i+2},\cdots,X_n),$$
  $$\sigma_i^{-1}{\mathbf X} = (X_1,\cdots,X_{i-1},\P_{({L_{\zg_i}\zg_{i+1}},f_{L_{\zg_i}\zg_{i+1}})},X_i,X_{i+2},\cdots,X_n).$$
\end{theorem}
\begin{proof}
(1) Since  $\zg_i$ and $\zg_{i+1}$ do not intersect, there exists no map from $X_i$ to $X_{i+1}[\mathbb{Z}]$. Therefore the distinguished triangles
\begin{eqnarray*}
\xymatrix{
X_i\ar[r]&0\ar[r]&X_{i}[1]\ar[r]&X_i[1],
}\\
\xymatrix{
X_{i+1}[-1]\ar[r]&X_{i+1}[-1]\ar[r]&0\ar[r]&X_{i+1}
}\end{eqnarray*}
are the respective exchange triangles of the right and left mutations in  \eqref{equation:right exchange braid} and \eqref{equation:left exchange braid}. So
$ R_{X_{i+1}}X_{i}=X_i$, $L_{X_i}X_{i+1}=X_{i+1}$, and
 $$\sigma_i{\mathbf X} = \sigma_i^{-1}{\mathbf X} = (X_1,\cdots,X_{i-1},X_{i+1},X_i,X_{i+2},\cdots,X_n).$$

(2) We only prove the case for the right braid group action in detail, the proof for left braid group action is similar.
For this we consider the two cases in Figure \ref{figure:right action}.

For case I, note that there exists a unique $m\in \mathbb{Z}$ such that the intersection $q_1$ induces a map
$$a: \P_{(\zg_i,f_{\zg_i})} \longrightarrow \P_{(\zg_{i+1},f_{\zg_{i+1}})}[m].$$

Moreover, note that up to isomorphism this is the unique map from $\P_{(\zg_i,f_{\zg_i})}$ to $\P_{(\zg_{i+1},f_{\zg_{i+1}})}[\mathbb{Z}]$, since there are no other intersections of $\zg_i$ and $\zg_{i+1}$. Thus for some proper grading $g$ of ${R_{\zg_{i+1}}\zg_{i}}$, the distinguished triangle
\begin{equation*}
X_i\s{a}\longrightarrow X_{i+1}[m]\longrightarrow \P_{({R_{\zg_{i+1}}\zg_{i}},g)}\longrightarrow X_i[1]
\end{equation*}
is exactly the triangle in \eqref{equation:right exchange braid} defining the right mutation. So $R_{X_{i+1}}X_{i}[1]=\P_{({R_{\zg_{i+1}}\zg_{i}},g)}$. Note that by \cite[Lemma 2.3]{CS20}, we have $g(s)=f_{\zg_i}(t)-1$. Set $f_{R_{\zg_{i+1}}\zg_{i}}=g[-1]$, then  $R_{X_{i+1}}X_{i}=\P_{({R_{\zg_{i+1}}\zg_{i}},f_{{R_{\zg_{i+1}}\zg_{i}}})}$ with $f_{R_{\zg_{i+1}}\zg_{i}}(s)=f_{\zg_i}(t)$.

For case II, there are unique $m, n \in \mathbb{Z}$ such that we have  homomorphisms
$$a_1: \P_{(\zg_i,f_{\zg_i})} \longrightarrow \P_{(\zg_{i+1},f_{\zg_{i+1}})}[m],$$
$$a_2: \P_{(\zg_i,f_{\zg_i})} \longrightarrow \P_{(\zg_{i+1},f_{\zg_{i+1}})}[n]$$
arising from the intersections at $q_1$ and $q_2$ respectively.
Then by \cite[Proposition 2.5]{CS20}, there is a grading $g$ over ${R_{\zg_{i+1}}\zg_{i}}$ such that
\begin{equation*}
X_i\s{(a_1,a_2)}\longrightarrow X_{i+1}[m]\oplus X_{i+1}[n]\longrightarrow \P_{({R_{\zg_{i+1}}\zg_{i}},g)}\longrightarrow X_i[1]
\end{equation*}
is a distinguished triangle.   On the other hand, note that
$a_1$ and $a_2$ form a basis of the homomorphism space from $\P_{(\zg_i,f_{\zg_i})}$ to $\P_{(\zg_{i+1},f_{\zg_{i+1}})}[\mathbb{Z}]$ since there are no other intersections of $\zg_i$ and $\zg_{i+1}$. So the above triangle is again  exactly the exchange triangle  in \eqref{equation:right exchange braid} defining the right mutation. Therefore $R_{X_{i+1}}X_{i}[1]=\P_{({R_{\zg_{i+1}}\zg_{i}},g)}$ with $g(s)=f_{\zg_i}(t)-1$. Then  $R_{X_{i+1}}X_{i}=\P_{({R_{\zg_{i+1}}\zg_{i}},f_{{R_{\zg_{i+1}}\zg_{i}}})}$ if we set $f_{R_{\zg_{i+1}}\zg_{i}}=g[-1]$, where $f_{R_{\zg_{i+1}}\zg_{i}}(s)=f_{\zg_i}(t)$.

To sum up, for both cases, we have $R_{X_{i+1}}X_{i}=\P_{({R_{\zg_{i+1}}\zg_{i}},f_{{R_{\zg_{i+1}}\zg_{i}}})}$ with $f_{R_{\zg_{i+1}}\zg_{i}}=g[-1]$, and thus
$$\sigma_i{\mathbf X} = (X_1,\cdots,X_{i-1},X_{i+1},\P_{({R_{\zg_{i+1}}\zg_{i}},f_{{R_{\zg_{i+1}}\zg_{i}}})},X_{i+2},\cdots,X_n).$$

\end{proof}

The following proposition can be directly derived
from above theorem and Proposition \ref{proposition:transitivity}.

\begin{proposition}\label{proposition:transitivity2}
If the {\bf RCEA} condition is satisfied for any marked surface, then for any gentle algebra $\A$ such that there are full exceptional sequences in $\ka$, the action of $\mathbb{Z}^n\ltimes B_n$ on the set of full exceptional sequences in $\ka$ is transitive.
\end{proposition}

Then Theorem \ref{theorem:transitivity2} follows from above proposition and Theorem \ref{theorem:transitivity}.
We have the following straightforward corollary from Theorem \ref{theorem:transitivity2}, noticing that for those derived-discrete algebra $\A$ which are hereditary, \cite{C93} has proven the transitivity of the action of $\mathbb{Z}^n\ltimes B_n$ on the set of ordered exceptional dissections in $\ka$ and that the remaining derived-discrete algebras are gentle algebras associated to surfaces of genus zero.

\begin{corollary}\label{corollary:transitivity}
Let $\A$ be a derived-discrete algebra. The action of $\mathbb{Z}^n\ltimes B_n$ on the set of ordered exceptional dissections in $\ka$ is transitive.
\end{corollary}

\section{Three dualities: Exceptional, Koszul and Serre dualities}\label{subsection:three dualities}

In this section we study the relation between the duality of a full exceptional sequence introduced in \cite[Section 8]{R90} and Koszul duality for  gentle algebras. We further show that the duality of exceptional sequences can also be explained using the Auslander-Reiten translation in $\ka$, thereby giving a new proof of this result by Bondal~\cite{Bon90} in the case of the perfect derived categories of gentle algebras.

For the rest of this section, we fix the following notation. Denote by $A$ a gentle algebra associated to a marked surface $(S,M,\zD_A)$ such that there exists a full exceptional sequence in $\ka$ and let $\mathbf X=(X_1,\cdots,X_n)$ be such a full exceptional sequence associated to an ordered exceptional dissection $\zG=(\zg_1,\cdots,\zg_n)$ on $(S,M)$. We define the ordered  dual exceptional dissection of $\Gamma$ to be the dual dissection $\Gamma^*$  with the following order $\zG^*=(\zg^*_n,\cdots,\zg^*_1)$. Denote by $\A(\zG)$ and $\A(\zG^*)$ the gentle algebras arising from $\zG$ and $\zG^*$ respectively. Then  $\A(\zG^*)$ is the Koszul dual of $\A(\zG)$, see Remark \ref{remark: dissection and quiver}. Furthermore, $\zG^*$ is exceptional, this follows directly from the fact that the quiver of $\A(\zG^*)$ is the opposite quiver of $\A(\zG)$.

For any $i \leq j \leq n$, set $$R^{j-i}X_i=R_{X_j}\cdots R_{X_{i+1}} X_{i},$$ and for any $1 \leq k \leq i$, set
$$L^{i-k}X_i=L_{X_k}\cdots L_{X_{i-1}} X_{i}.$$ In particular, we have $R^{j-i}X_i=X_i$ and $L^{i-k}X_i=X_i$ if $j=i$ and $k=i$ respectively.
We will denote by $R^{j-i}\zg_i$ the arc associated to $R^{j-i}X_i$, and denote by $L^{i-k}\zg_i$ the arc associated to $L^{i-k}X_i$.

\begin{definition}\cite{R90}\label{definition:exceptional duality}
Let $\mathbf X=(X_1,\cdots,X_n)$ be a full exceptional sequence in $\ka$.
We call
$\mathbf{RX}=(X_n,R^1X_{n-1}\cdots,R^{n-1}X_1)$ the \emph{right dual}  of $\mathbf X$
 and
$\mathbf{LX}=(L^{n-1}X_n,L^{n-2}X_{n-1}\cdots,X_1)$   \emph{left dual} of $\mathbf X$.
\end{definition}

The following lemma directly follows from the definitions.

\begin{lemma}\label{lemma:exceptional duality}
Let  $\omega_0=\sigma_1(\sigma_2\sigma_1)
\cdots(\sigma_{n-1}\sigma_{n-2}\cdots\sigma_1)$ be the element in $ B_n$ corresponding to the longest element in the associated symmetric group $S_n$. Then
$\mathbf{RX}=\omega_0\mathbf X$ and
$\mathbf{LX}=\omega_0^{-1}\mathbf X$. So we have $\mathbf{R}=\mathbf{L}^{-1}$ on the set of full exceptional sequences in $\ka$.
\end{lemma}

  We now give an interpretation of the left and right dual of a full exceptional sequence in terms of the geometric surface model of $A$. For this,  note that in an exceptional dissection we cannot have polygons such that if the polygon is unfolded one of its edges appears twice since such an arc must be contained in a non-exceptional cycle. An example of such a non-exceptional cycle is given in Figure~\ref{figure:configuration of intersections0} where the arc $\zg_1$ would appear twice in the unfolded polygon.
Therefore, in particular, for any $1\leq i \leq n$, $\zg_i$ belongs to two polygons  of $\zG$, as in Figure \ref{figure:two dualities}, where $k$ and $m$ are allowed to be zero, which means that the $\rpoint$-marked points $q$ and $q'$ in the polygons lie in the same boundary component, otherwise we have $1\leq j \leq k \leq n$,  and $1\leq l \leq m \leq n$. Denote the sets of arcs in Figure  \ref{figure:two dualities} by $\Sigma=\{\za_1,\cdots,\za_k\}$ and $\Sigma'=\{\za'_1,\cdots,\za'_m\}$.

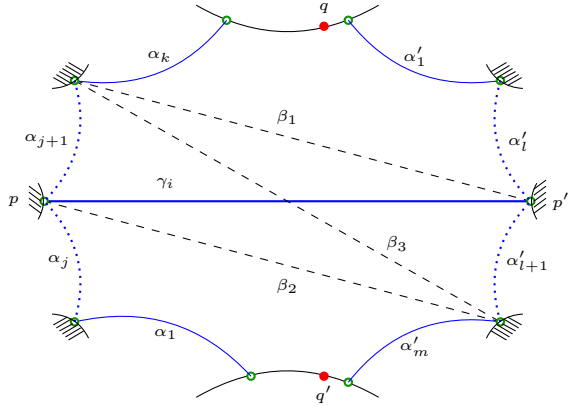
\begin{figure}[H]
 \[\scalebox{1}{
\begin{tikzpicture}[scale=0.8]
\draw[thick,blue](-4,0)--(4,0);
    \draw (-2,0.3) node {\tiny$\zg_i$};
    \draw (-2,-2.2) node {\tiny$\za_1$};
    \draw (-3.75,-1) node {\tiny$\za_j$};
    \draw (-2.1,2.4) node {\tiny$\za_k$};
    \draw (-3.95,1) node {\tiny$\za_{j+1}$};

    \draw (2.1,2.4) node {\tiny$\za'_1$};
    \draw (3.8,1) node {\tiny$\za'_l$};
    \draw (2.1,-2.4) node {\tiny$\za'_m$};
    \draw (3.95,-1) node {\tiny$\za'_{l+1}$};
    \draw (0,1.4) node {\tiny$\zb_1$};
    \draw (0,-1.4) node {\tiny$\zb_2$};
    \draw (1.8,-0.7) node {\tiny$\zb_3$};

\draw[dashed](-3.5,2)--(4,0);
\draw[dashed](3.5,-2)--(-4,0);
\draw[dashed](-3.5,2)--(3.5,-2);

\draw[dark-green,thick,fill=white] (4,0) circle (0.06);
\draw[dark-green,thick,fill=white] (-4,0) circle (0.06);
    \draw (4.5,0) node {\tiny$p'$};
    \draw (-4.5,0) node {\tiny$p$};

	\draw[bend right](4.1,0.3)to(4.1,-0.3);
    \draw[-]  (4.25,0.35)to(4.05,0.2);
    \draw[-]  (4.25,0.35-0.1)to(4.05,0.2-0.1);
    \draw[-]  (4.25,0.35-0.2)to(4.05,0.2-0.2);
    \draw[-]  (4.25,0.35-0.3)to(4.05,0.2-0.3);
    \draw[-]  (4.25,0.35-0.4)to(4.05,0.2-0.4);

	\draw[bend left](-4.1,0.3)to(-4.1,-0.3);
    \draw[-]  (-4.25,0.35-0.1)to(-4.05,0.2-0.1);
    \draw[-]  (-4.25,0.35-0.2)to(-4.05,0.2-0.2);
    \draw[-]  (-4.25,0.35-0.3)to(-4.05,0.2-0.3);
    \draw[-]  (-4.25,0.35-0.4)to(-4.05,0.2-0.4);


	\draw[bend left](1.5,3.25)to(-1.5,3.25);
    \draw[red,thick,fill] (.6,2.9) circle (0.06);
    \draw (.6,3.2) node {\tiny$q$};
    \draw (.6,-3.2) node {\tiny$q'$};

\draw[dark-green,thick,fill=white] (1,3) circle (0.06);
	\draw[bend right,blue](1,3)to(3.5,2);
\draw[dark-green,thick,fill=white] (3.5,2) circle (0.06);
	\draw[bend right,blue,dotted,thick](3.5,2)to(4,0);

\draw[dark-green,thick,fill=white] (-1,3) circle (0.06);
	\draw[bend left,blue](-1,3)to(-3.5,2);
\draw[dark-green,thick,fill=white] (-3.5,2) circle (0.06);
	\draw[bend left,blue,dotted,thick](-3.5,2)to(-4,0);

    \draw[bend right](3.27,2.27)to(3.87,1.87);
    \draw[-]  (3.6,2.3)to(3.35,2.15);
    \draw[-]  (3.6+0.05,2.3-0.05)to(3.35+0.05,2.15-0.05);
    \draw[-]  (3.6+0.1,2.3-0.1)to(3.35+0.1,2.15-0.1);
    \draw[-]  (3.6+0.15,2.3-0.15)to(3.35+0.18,2.15-0.12);
    \draw[-]  (3.6+0.2,2.3-0.2)to(3.35+0.2,2.15-0.2);
    \draw[-]  (3.6+0.25,2.3-0.25)to(3.35+0.27,2.15-0.23);

    \draw[bend left](-3.27,2.27)to(-3.87,1.87);
    \draw[-]  (-3.6,2.3)to(-3.35,2.15);
    \draw[-]  (-3.6-0.05,2.3-0.05)to(-3.35-0.05,2.15-0.05);
    \draw[-]  (-3.6-0.1,2.3-0.1)to(-3.35-0.1,2.15-0.1);
    \draw[-]  (-3.6-0.15,2.3-0.15)to(-3.35-0.18,2.15-0.12);
    \draw[-]  (-3.6-0.2,2.3-0.2)to(-3.35-0.2,2.15-0.2);
    \draw[-]  (-3.6-0.25,2.3-0.25)to(-3.35-0.27,2.15-0.23);


	\draw[bend right](1.5,-3.25)to(-1.5,-3.25);
    \draw[red,thick,fill] (.6,-2.9) circle (-0.06);

\draw[dark-green,thick,fill=white] (1,-3) circle (0.06);
	\draw[bend left,blue](1,-3)to(3.5,-2);
\draw[dark-green,thick,fill=white] (3.5,-2) circle (0.06);
	\draw[bend left,blue,dotted,thick](3.5,-2)to(4,0);

\draw[dark-green,thick,fill=white] (-.6,-2.9) circle (0.06);
 	\draw[bend right,blue](-.6,-2.9)to(-3.5,-2);
\draw[dark-green,thick,fill=white] (-3.5,-2) circle (0.06);
    \draw[bend right,blue,dotted,thick](-3.5,-2)to(-4,0);

    \draw[bend left](3.27,-2.27)to(3.87,-1.87);
    \draw[-]  (3.6,-2.3)to(3.35,-2.15);
    \draw[-]  (3.6+0.05,-2.3+0.05)to(3.35+0.05,-2.15+0.05);
    \draw[-]  (3.6+0.1,-2.3+0.1)to(3.35+0.1,-2.15+0.1);
    \draw[-]  (3.6+0.15,-2.3+0.15)to(3.35+0.18,-2.15+0.12);
    \draw[-]  (3.6+0.2,-2.3+0.2)to(3.35+0.2,-2.15+0.2);
    \draw[-]  (3.6+0.25,-2.3+0.25)to(3.35+0.27,-2.15+0.23);

    \draw[bend right](-3.27,-2.27)to(-3.87,-1.87);
    \draw[-]  (-3.6,-2.3)to(-3.35,-2.15);
    \draw[-]  (-3.6-0.05,-2.3+0.05)to(-3.35-0.05,-2.15+0.05);
    \draw[-]  (-3.6-0.1,-2.3+0.1)to(-3.35-0.1,-2.15+0.1);
    \draw[-]  (-3.6-0.15,-2.3+0.15)to(-3.35-0.18,-2.15+0.12);
    \draw[-]  (-3.6-0.2,-2.3+0.2)to(-3.35-0.2,-2.15+0.2);
    \draw[-]  (-3.6-0.25,-2.3+0.25)to(-3.35-0.27,-2.15+0.23);
 \end{tikzpicture}}\]
\begin{center}
\caption{An $\gpoint$-arc $\zg_i$ of an exceptional dissection belongs to two polygons, where two marked points and two $\gpoint$-arcs may coincide.}\label{figure:two dualities}
\end{center}
\end{figure}

Let $s$ be the intersection $\zg_i\cap\zg^*_i$, noting that following our conventions, this is the unique intersection of $\zg_i$ and  $\zg^*_i$.
Recall that $D\zg_i^*$ and $D^{-1}\zg^*_i$ are the direct  and  inverse twists of $\zg^*_i$ to the next $\gpoint$-point, respectively. We also have the following  intersections, which are in fact the unique intersections between the two arcs in consideration,
\[u=\zg_i\cap D\zg^*_i,~ u'=\zg_i\cap D^{-1}\zg^*_i,~ u''=D\zg^*_i\cap D^{-1}\zg^*_i,\]
where $u''$ always is in the interior of $S$, while $u$ and $u'$ may be on the boundary, see Figure \ref{figure:two dualities2}.

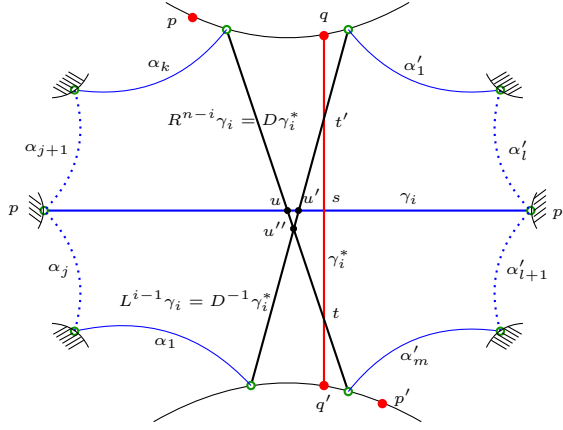
\begin{figure}[H]
 \[\scalebox{1}{
\begin{tikzpicture}[scale=0.8]

    \draw (-2,-2.2) node {\tiny$\za_1$};
    \draw (-3.75,-1) node {\tiny$\za_j$};
    \draw (-2.1,2.4) node {\tiny$\za_k$};
    \draw (-3.95,1) node {\tiny$\za_{j+1}$};

    \draw (2.1,2.4) node {\tiny$\za'_1$};
    \draw (3.8,1) node {\tiny$\za'_l$};
    \draw (2.1,-2.4) node {\tiny$\za'_m$};
    \draw (3.95,-1) node {\tiny$\za'_{l+1}$};

\draw[thick,blue](-4,0)--(4,0);
\draw[thick,red](0.6,-2.9)--(0.6,2.9);
\draw[thick](-1,3)--(1,-3);
\draw[thick](-.6,-2.9)--(1,3);
\draw[thick,fill] (0,0) circle (0.04);
\draw[thick,fill] (.18,0) circle (0.04);
\draw[thick,fill] (0.1,-.3) circle (0.04);
    \draw (.8,.2) node {\tiny$s$};
    \draw (.9,1.5) node {\tiny$t'$};
    \draw (.8,-1.7) node {\tiny$t$};

    \draw (-.2,.15) node {\tiny$u$};
    \draw (-.2,-.3) node {\tiny$u''$};
    \draw (.4,.25) node {\tiny$u'$};
    \draw (2,0.2) node {\tiny$\zg_i$};
    \draw (0.85,-0.8) node {\tiny$\zg_i^*$};
    \draw (-0.88,1.5) node {\tiny$R^{n-i}\zg_i=D\zg_i^*$};
    \draw (-1.5,-1.5) node {\tiny$L^{i-1}\zg_i=D^{-1}\zg_i^*$};
%

\draw[dark-green,thick,fill=white] (4,0) circle (0.06);
\draw[dark-green,thick,fill=white] (-4,0) circle (0.06);
    \draw (4.5,0) node {\tiny$p'$};
    \draw (-4.5,0) node {\tiny$p$};

	\draw[bend right](4.1,0.3)to(4.1,-0.3);
    \draw[-]  (4.25,0.35)to(4.05,0.2);
    \draw[-]  (4.25,0.35-0.1)to(4.05,0.2-0.1);
    \draw[-]  (4.25,0.35-0.2)to(4.05,0.2-0.2);
    \draw[-]  (4.25,0.35-0.3)to(4.05,0.2-0.3);
    \draw[-]  (4.25,0.35-0.4)to(4.05,0.2-0.4);

	\draw[bend left](-4.1,0.3)to(-4.1,-0.3);
    \draw[-]  (-4.25,0.35-0.1)to(-4.05,0.2-0.1);
    \draw[-]  (-4.25,0.35-0.2)to(-4.05,0.2-0.2);
    \draw[-]  (-4.25,0.35-0.3)to(-4.05,0.2-0.3);
    \draw[-]  (-4.25,0.35-0.4)to(-4.05,0.2-0.4);


	\draw[bend left](2,3.45)to(-2,3.45);
    \draw[red,thick,fill] (0.6,2.9) circle (0.06);
    \draw (.6,3.2) node {\tiny$q$};
    \draw (.6,-3.2) node {\tiny$q'$};

    \draw[red,thick,fill] (-1.56,3.2) circle (0.06);
    \draw (-1.9,3.1) node {\tiny$p$};

\draw[dark-green,thick,fill=white] (1,3) circle (0.06);
	\draw[bend right,blue](1,3)to(3.5,2);
\draw[dark-green,thick,fill=white] (3.5,2) circle (0.06);
	\draw[bend right,blue,dotted,thick](3.5,2)to(4,0);

\draw[dark-green,thick,fill=white] (-1,3) circle (0.06);
	\draw[bend left,blue](-1,3)to(-3.5,2);
\draw[dark-green,thick,fill=white] (-3.5,2) circle (0.06);
	\draw[bend left,blue,dotted,thick](-3.5,2)to(-4,0);

    \draw[bend right](3.27,2.27)to(3.87,1.87);
    \draw[-]  (3.6,2.3)to(3.35,2.15);
    \draw[-]  (3.6+0.05,2.3-0.05)to(3.35+0.05,2.15-0.05);
    \draw[-]  (3.6+0.1,2.3-0.1)to(3.35+0.1,2.15-0.1);
    \draw[-]  (3.6+0.15,2.3-0.15)to(3.35+0.18,2.15-0.12);
    \draw[-]  (3.6+0.2,2.3-0.2)to(3.35+0.2,2.15-0.2);
    \draw[-]  (3.6+0.25,2.3-0.25)to(3.35+0.27,2.15-0.23);

    \draw[bend left](-3.27,2.27)to(-3.87,1.87);
    \draw[-]  (-3.6,2.3)to(-3.35,2.15);
    \draw[-]  (-3.6-0.05,2.3-0.05)to(-3.35-0.05,2.15-0.05);
    \draw[-]  (-3.6-0.1,2.3-0.1)to(-3.35-0.1,2.15-0.1);
    \draw[-]  (-3.6-0.15,2.3-0.15)to(-3.35-0.18,2.15-0.12);
    \draw[-]  (-3.6-0.2,2.3-0.2)to(-3.35-0.2,2.15-0.2);
    \draw[-]  (-3.6-0.25,2.3-0.25)to(-3.35-0.27,2.15-0.23);


	\draw[bend right](2.2,-3.5)to(-2.2,-3.5);
    \draw[red,thick,fill] (0.6,-2.9) circle (-0.06);
    \draw[red,thick,fill] (1.56,-3.2) circle (0.06);
    \draw (1.9,-3.1) node {\tiny$p'$};

\draw[dark-green,thick,fill=white] (1,-3) circle (0.06);
	\draw[bend left,blue](1,-3)to(3.5,-2);
\draw[dark-green,thick,fill=white] (3.5,-2) circle (0.06);
	\draw[bend left,blue,dotted,thick](3.5,-2)to(4,0);

\draw[dark-green,thick,fill=white] (-.6,-2.9) circle (0.06);
 	\draw[bend right,blue](-.6,-2.9)to(-3.5,-2);
\draw[dark-green,thick,fill=white] (-3.5,-2) circle (0.06);
	\draw[bend right,blue,dotted,thick](-3.5,-2)to(-4,0);

    \draw[bend left](3.27,-2.27)to(3.87,-1.87);
    \draw[-]  (3.6,-2.3)to(3.35,-2.15);
    \draw[-]  (3.6+0.05,-2.3+0.05)to(3.35+0.05,-2.15+0.05);
    \draw[-]  (3.6+0.1,-2.3+0.1)to(3.35+0.1,-2.15+0.1);
    \draw[-]  (3.6+0.15,-2.3+0.15)to(3.35+0.18,-2.15+0.12);
    \draw[-]  (3.6+0.2,-2.3+0.2)to(3.35+0.2,-2.15+0.2);
    \draw[-]  (3.6+0.25,-2.3+0.25)to(3.35+0.27,-2.15+0.23);

    \draw[bend right](-3.27,-2.27)to(-3.87,-1.87);
    \draw[-]  (-3.6,-2.3)to(-3.35,-2.15);
    \draw[-]  (-3.6-0.05,-2.3+0.05)to(-3.35-0.05,-2.15+0.05);
    \draw[-]  (-3.6-0.1,-2.3+0.1)to(-3.35-0.1,-2.15+0.1);
    \draw[-]  (-3.6-0.15,-2.3+0.15)to(-3.35-0.18,-2.15+0.12);
    \draw[-]  (-3.6-0.2,-2.3+0.2)to(-3.35-0.2,-2.15+0.2);
    \draw[-]  (-3.6-0.25,-2.3+0.25)to(-3.35-0.27,-2.15+0.23);
 \end{tikzpicture}}\]
\begin{center}
\caption{$R^{n-i}\zg_i=D\zg_i^*$ and $L^{i-1}\zg_i=D^{-1}\zg_i^*$ for an $\gpoint$-arc $\zg_i$ in an exceptional dissection.}\label{figure:two dualities2}
\end{center}
\end{figure}

The following proposition gives an explicit description of the object $R^{n-i}X_i$ in terms of the twist of the dual $\zg^*_i$ of the arc $\zg_i$ corresponding to $X_i$.

\begin{proposition}\label{proposition:two dualities} With the notation above, we have $R^{n-i}X_i=\P_{(D\zg_i^*,f)}$, where $f$ is the grading of $D\zg_i^*$ such that there is a map from $\P_{(D\zg_i^*,f)}$ to $X_i$ arising from the intersection $u$ of $D\zg^*_i$ and $\zg_i$.
\end{proposition}
\begin{proof}
At first we show that $R^{n-i}\zg_i=D\zg_i^*$. Following our conventions, see Figure~\ref{figure:two dualities},
we have $$\za_1\preceq \cdots \preceq \za_j\preceq \zg_i \preceq \za_{j+1} \preceq \cdots \preceq \za_k,$$
$$\za'_1\preceq \cdots \preceq \za'_l\preceq \zg_i \preceq \za'_{l+1} \preceq \cdots \preceq \za'_m.$$
Then either
 $\zg_{i+1}=\za_{j+1}$ or $\za'_{l+1}$ if $\zg_{i+1}\in \Sigma\sqcup \Sigma'$, or  $\zg_{i+1}\in\zG\setminus (\Sigma\sqcup \Sigma')$. In either case,  $\zg_{i+1}\in \{\za_{j+1},\za'_{l+1}\}\sqcup(\zG\setminus (\Sigma\sqcup \Sigma'))$.

\emph{Case 1.} If $\zg_{i+1}=\za_{j+1}\neq\za'_{l+1}$, then by Theorem \ref{thm:braid group actions} (2), $R\zg_i=\zb_1$ as shown in Figure \ref{figure:two dualities} (compare also to the left picture in Figure \ref{figure:right action}).

\emph{Case 2.} If $\zg_{i+1}=\za'_{l+1}\neq\za_{j+1}$, then similar to the first case, $R\zg_i=\zb_2$ as shown in Figure \ref{figure:two dualities}.

\emph{Case 3.} If $\zg_{i+1}=\za_{j+1}=\za'_{l+1}$, then by Theorem \ref{thm:braid group actions} (2), $R\zg_i=\zb_3$ as in Figure \ref{figure:two dualities} (compare to the right picture in Figure \ref{figure:right action}).

\emph{Case 4.} If $\zg_{i+1}\in \zG\setminus (\Sigma\sqcup \Sigma')$ then  we show that $\Hom(X_i,X_{i+1})=0$ and by Theorem \ref{thm:braid group actions} (1), $R\zg_i=\zg_i$. Suppose for contradiction that   $\Hom(X_i,X_{i+1}) \neq 0$. Then $\zg_i\cap\zg_{i+1}\in \{p,p'\}$. If $p\in\zg_i\cap\zg_{i+1}$, then $\zg_{i+1}$ follows $\zg_i$ anticlockwise at $p$.  On the other hand, $\za_{j+1}$ directly  follows $\zg_i$ at $p$ in the anticlockwise direction, thus $\za_{j+1}\preceq \zg_{i+1}$. This is a contradiction since $\zg_{i+1}\preceq \za_{j+1}$ ($\zg_{i+1}$ is a minimal element which is larger than $\zg_i$) and $\zg_{i+1}\neq \za_{j+1}$. The argument for the case  $p'\in\zg_i\cap\zg_{i+1}$ is similar.

To sum up, for all the cases, $R\zg_i=\zg_i$ or $R\zg_i$ is obtained by rotating $\zg_i$ clockwise at one endpoint or at both endpoints.
Now let
$$\sigma_i\mathbf X=\{X_1,\cdots,X_{i-1},X_{i+1},RX_i,X_{i+2},\cdots,X_n\}$$
be the new full exceptional sequence.
Using a similar discussion for $R\zg_i$ as  for $\zg_i$ above, we know that $R^2\zg_i=R\zg_i$ or $R^2\zg_i$ is obtained by rotating $R\zg_i$ clockwise at one endpoint or at both endpoints. Continue this process until we obtain a full exceptional dissection
\begin{equation*}\label{equ:excp}
  \sigma_{i+r-1}\cdots\sigma_i(\mathbf X)=(X_1,\cdots,X_{i-1},X_{i+1},\cdots,X_{i+r},R^rX_i,X_{i+r+1},\cdots,X_n)
\end{equation*}
such that we can not continue the process any more, that is, $q$ is the $\rpoint$-point which directly  follows (anticlockwise) an endpoint of $R^r\zg_i$ and $q'$ is the $\rpoint$-point which  directly  follows (anticlockwise) the other endpoint of $R^r\zg_i$.  In this case, there is no intersection between $\zg_v$ and $R^r\zg_i$, for any $i+r+1\leq v \leq n$. Otherwise, they intersect at the endpoints of $R^r\zg_i$ which implies $\zg_v\preceq R^r\zg_i$, see Figure \ref{figure:two dualities2}. But this contradicts the fact that $\sigma_{i+r-1}\cdots\sigma_i(\mathbf X)$ is an exceptional sequence.
Therefore $\Hom(R^rX_i,X_v)=0$ for any $i+r+1\leq v \leq n$ and by inductively using Theorem \ref{thm:braid group actions} (1), we obtain $R^rX_i=R^{n-i}X_i$ and thus $R^r\zg_i=R^{n-i}\zg_i$. On the other hand, since $R^r\zg_i=D\zg_i^*$, so $R^{n-i}\zg_i=D\zg_i^*$.

Now we consider the grading $f$.
By the definition of $\sigma_{j}, i \leq j \leq i+r-1$, we have the following exchange triangle

\[\xymatrix{
R^{j-i}X_i\ar[r]&\coprod_{\ell\in\mathbb{Z}}
D\Hom(R^{j-i}X_i,R^{j-i+1}X_i[\ell])\otimes_kR^{j-i+1}X_i[\ell]\ar[r]&\\
}\]
\[\xymatrix{R^{j-i+1}X_i[1]\ar[r]^{a_j}&R^{j-i}X_i[1],}\]
where the last map $a_j$ is nonzero.
Furthermore, the composition $a_{i} \cdots a_{n}$ is again a nonzero map from $R^{n-i}X_i$ to $X_i$, which corresponds to the intersection $u$ of $R^{n-i}\zg_i=D\zg^*$ and $\zg_i$. This implies that $f$ is the grading such that there is a non-zero map from $R^{n-i}X_i=\P_{(R^{n-i}\zg_i,f)}$ to $X_i$.
\end{proof}

The following proposition gives an explicit description of the object $L^{i-1}X_i$. The proof is similar to the proof of Proposition \ref{proposition:two dualities}.

\begin{proposition}\label{proposition:two dualities2}
With the notation above we have $L^{i-1}X_i=\P_{(L^{i-1}\zg_i,g)}$, where $L^{i-1}\zg_i=D^{-1}(\zg_i^*)$ and $g$ is the grading such that there is a map from $X_i$ to $L^{i-1}X_i$ arising from the intersection $u'$.
\end{proposition}

The following theorem shows that the right (resp. left) duality of an exceptional sequence corresponds to Koszul duality   followed by  twisting (resp. inversely twisting) the  arcs in the corresponding dual exceptional dissection.

\begin{theorem}\label{theorem:two dualities1}
Let $\mathbf X=(X_1,\cdots,X_n)$ be a full exceptional sequence in $\ka$ associated to an ordered exceptional dissection $\zG=(\zg_1,\cdots,\zg_n)$ of $(S,M)$, that is  $X_i=\P_{(\zg_i,f_{\zg_i})}$ for some grading $f_{\zg_i}$. Then the right dual
$\mathbf{RX}$ is induced  by  the equality
 $$\mathbf{R}\zG=D(\zG^*)$$ where $D(\zG^*)$ is the twist of the ordered dual exceptional dissection $\zG^*$.

The left dual $\mathbf{LX}$ is induced  by  the equality
$$\mathbf{L}\zG=D^{-1}(\zG^*)$$ where $D^{-1}(\zG^*)$ is the inverse twist of the ordered dual exceptional dissection $\zG^*$.
\end{theorem}

\begin{proof}
The proof directly follows from Lemma \ref{lemma:exceptional duality} and Proposition \ref{proposition:two dualities}.
\end{proof}

Note that in Theorem~\ref{theorem:two dualities1} when we say that $\mathbf{RX}$ is induced  by  the equality
 $\mathbf{R}\zG=D(\zG^*)$ we mean that the grading of $\mathbf{RX}$ is  induced by the gradings as described in Proposition~\ref{proposition:two dualities} and similarly for $\mathbf{LX}$ the grading is the one induced by the gradings in Proposition~\ref{proposition:two dualities2}.

It is shown in \cite{BM90} that if a triangulated category has full exceptional sequences, then it has a Serre functor, and in \cite[Proposition 4.2]{Bon90} that the duality of the exceptional sequences can be obtained using the Serre functor. In the following we give a new proof of this result for $\ka$.

\begin{theorem}\label{theorem:two dualities}
Let $\mathbf X=(X_1,\cdots,X_n)$ be a full exceptional sequence in $\ka$ and let $\mathbb{S}$ be the Serre functor in $\ka$. Then $\tau(\mathbf X) = {\bf L}^{2} [-1](\mathbf X)$ and thus $\mathbb{S} (\mathbf X) = {\mathbf L}^2(\mathbf X)$.
\end{theorem}

\begin{proof}
We only need to show $\tau(\mathbf X) = {\bf L}^{2} [-1](\mathbf X)$, since $\tau = \mathbb{S}[-1]$.
 Let $X_i=\P_{(\zg_i,f_i)}$ so that  $\zG=(\zg_1, \zg_2,\cdots,\zg_n)$ is the ordered exceptional dissection corresponding to $\bf X$ and set ${\bf L}^{2}(\mathbf X)=(Y_1,Y_2,\cdots,Y_n)$. Then by Theorem \ref{theorem:two dualities1}, $\mathbf{L}\zG=D^{-1}(\zG^*)$. Furthermore, it is easy to see that the twist operator $D$ (as well as $D^{-1}$) and the dual operator $*$ commute, that is, $D^{-1}(\zG^*)=(D^{-1}(\zG))^*$. So we have $\mathbf{L}^2(\zG)=D^{-2}(\zG)$.
On the other hand, we know from Lemma \ref{lemma:tau and Dehn} that the arc of $\tau X_i$ is exactly $D^{-2}\zg_i$. Thus $\tau X_i$ and $Y_i$ correspond to the same arc $D^{-2}\zg_i$.

Now we consider the grading.
On the one hand, there is a non-zero map $F$ from $X_i[-1]$ to $\tau X_i$ arising from the  Auslander-Reiten triangle for  $X_i$.
On the other hand, by Proposition \ref{proposition:two dualities2}, there is a non-zero map $G$ from $X_i$ to $Y_i$. Note that both maps $F$ and $G$ arise from the unique intersection between $D^{-2}\zg_i$ and $\zg_i$. Thus we have $\tau X_i=Y_i[-1]$ and $\tau(\mathbf X) = {\bf L}^{2} [-1](\mathbf X)$.
\end{proof}

The following corollary directly follows from Theorem \ref{theorem:two dualities}.

\begin{corollary}
Let $\mathbf X=(X_1,\cdots,X_n)$ be a full exceptional sequence in $\ka$. Then $L^{n-1}X_n=\mathbb{S} X_n=\tau X_n [1]$ and $R^{n-1}X_1=\mathbb{S}^{-1} X_1=\tau^{-1} X_1 [-1]$.
\end{corollary}

\end{document}